\documentclass{amsart}

\usepackage{amsmath}
\usepackage{amsfonts}
\usepackage{amssymb}
\usepackage{amsthm}
\usepackage{comment}
\usepackage{epsfig}
\usepackage{psfrag}
\usepackage{mathrsfs}
\usepackage{amscd}
\usepackage[all]{xy}
\usepackage{rotating}
\usepackage{lscape}
\usepackage{amsbsy}
\usepackage{verbatim}
\usepackage{moreverb}
\usepackage{color}
\usepackage{bbm}
\usepackage{eucal}
\usepackage{stmaryrd}
\usepackage{quiver}

\usepackage{tikz-cd}
\usetikzlibrary{patterns,shapes.geometric,arrows,decorations.markings}
\usepackage{tikz-3dplot}

\usepackage{caption}
\usepackage{subcaption}

\usepackage[colorlinks=true,linkcolor=blue,citecolor=red]{hyperref}

\colorlet{lightgray}{black!15}

\tikzset{->-/.style={decoration={
  markings,
  mark=at position .5 with {\arrow{>}}},postaction={decorate}}}
\tikzset{midarrow/.style={decoration={
    markings,
    mark=at position {#1} with {\arrow{>}}},postaction={decorate}}}

\newtheorem{theorem}{Theorem}[section]
\newtheorem{prop}[theorem]{Proposition}
\newtheorem{lemma}[theorem]{Lemma}
\newtheorem{cor}[theorem]{Corollary}
\newtheorem{conj}[theorem]{Conjecture}

\theoremstyle{definition}
\newtheorem{definition}[theorem]{Definition}

\newtheorem{remark}[theorem]{Remark}
\newtheorem{example}[theorem]{Example}

\newtheorem{l.notation}[theorem]{Local Notation}

\theoremstyle{remark}

\definecolor{orange}{rgb}{.95,0.5,0}
\definecolor{light-gray}{gray}{0.75}
\definecolor{brown}{cmyk}{0, 0.8, 1, 0.6}
\definecolor{plum}{rgb}{.5,0,1}

\DeclareMathOperator{\ev}{\mathsf{ev}}

\DeclareMathOperator{\Alg}{\sf Alg}

\DeclareMathOperator{\Mod}{\sf Mod}

\DeclareMathOperator{\CAlg}{\sf CAlg}

\DeclareMathOperator{\colim}{{\sf colim}}

\DeclareMathOperator{\limit}{{\sf lim}}

\DeclareMathOperator{\Der}{\sf Der}

\DeclareMathOperator{\Hom}{\sf Hom}

\DeclareMathOperator{\Fun}{{\sf Fun}}

\DeclareMathOperator{\Map}{{\sf Map}}

\DeclareMathOperator{\Spec}{\sf Spec}

\DeclareMathOperator{\ind}{\sf Ind}

\DeclareMathOperator{\pt}{\Spec(\Bbbk)}
\DeclareMathOperator{\Sym}{\sf Sym}

\DeclareMathOperator{\Cat}{{\sf Cat}}

\DeclareMathOperator{\qc}{\sf QC}

\DeclareMathOperator{\triv}{\mathsf{triv}}

\DeclareMathOperator{\op}{\mathsf{op}}

\DeclareMathOperator{\Loc}{\sf LocSys}

\DeclareMathOperator{\spaces}{\mathsf{Spaces}}
\DeclareMathOperator{\Spaces}{\mathsf{Spaces}}

\DeclareMathOperator{\Disk}{\sf Disk}

\DeclareMathOperator{\Lie}{\sf Lie}

\newcommand{\ra}{\rightarrow}
\newcommand{\la}{\leftarrow}
\newcommand{\xra}{\xrightarrow}
\newcommand{\xla}{\xleftarrow}

\def\cA{\mathcal A}\def\cC{\mathcal C}\def\cD{\mathcal D}

\def\cL{\mathcal L}
\def\cO{\mathcal O}\def\cP{\mathcal P}

\def\cV{\mathcal V}

\def\AA{\mathbb A}\def\BB{\mathbb B}\def\DD{\mathbb D}
\def\EE{\mathbb E}\def\GG{\mathbb G}
\def\LL{\mathbb L}
\def\PP{\mathbb P}
\def\RR{\mathbb R}\def\TT{\mathbb T}

\def\ZZ{\mathbb Z}

\def\sB{\mathsf B}\def\sC{\mathsf C}

\def\sO{\mathsf O}
\def\sS{\mathsf S}

\def\gg{\frak G}

\DeclareMathOperator{\DR}{\mathsf{DR}}

\DeclareMathOperator{\dSt}{\mathsf{Stacks}}

\DeclareMathOperator{\id}{\sf id}

\DeclareMathOperator{\Rep}{\sf Rep}

\DeclareMathOperator{\Perf}{\sf Perf}

\def\gg{\frak g}

\DeclareMathOperator{\gr}{{\sf gr}}

\DeclareMathOperator{\LR}{\mathbb{L}_{\mathit{R}_{\sf red}| \mathit{R}}}

\begin{document}

\title{Quantization of Betti Mapping Stacks}

\author{John Francis, Owen Gwilliam, and Rajiv Kaipa}

\address{Department of Mathematics\\Northwestern University\\Evanston, IL 60208}
\email{jnkf@math.northwestern.edu}
\address{Department of Mathematics\\ University of Massaschusetts, Amherst\\ MA}
\email{ogwilliam@umass.edu}
\address{Department of Mathematics\\Northwestern University\\Evanston, IL 60208}
\email{rajivkaipa2027@u.northwestern.edu }
\thanks{JF was supported by the National Science Foundation under awards 1812057 and 2604507. 
This material is based upon work supported by the National Science Foundation under Grant No. DMS-1440140, while JF and OG were in residence at the Mathematical Sciences Research Institute in Berkeley, California, during the Spring 2020 semester, and by Grant No. DMS-2424139, while JF was in residence at the Simons Laufer Mathematical Sciences Institute in Berkeley, California, during the Fall 2026 semester.}

\begin{abstract}
    We describe a quantization procedure for Betti mapping stacks with shifted Poisson targets by combining factorization homology and formal geometry. We quantize derived stacks of local systems of low-dimensional manifolds. 
\end{abstract}

\keywords{}

\subjclass[2020]{Primary 14A30, 13D10. Secondary 55N35, 53D55, 14M35, 57R56.}

\maketitle

\tableofcontents
\section{Introduction}
We quantize Betti mapping stacks by combining the techniques of factorization homology and Gelfand--Kazhdan formal geometry ~\cite{Gelfand:1971jg}. This is a version of the AKSZ construction. In particular, we quantize the derived stack of local systems on low-dimensional manifolds. This approach globalizes the perturbative quantization of Chern--Simons theory by Axelrod--Singer and Kontsevich over moduli spaces of classical solutions,
and is a topological analogue of Beilinson--Drinfeld's construction of differential operators on the moduli space of bundles on an algebraic curve.

We now articulate our results, after recollecting essential background from \cite{CPTVV} and~\cite{oldfact}. 

Let $\Bbbk$ be a field of characteristic~$0$. 
Let $X$ be an Artin derived stack, locally of finite presentation over $\Bbbk$. 
Its functions $\cO_X$ can be encoded as flat sections of the infinite jet bundle (with its canonical connection) on~$X$, 
as we now describe, following~\cite{CPTVV}.

Recall that $X$ has an associated prestack $X_{\DR}$, known as its {\it de Rham stack},
whose functor of points assigns to any connective dg commutative algebra $R$, the space (or $\infty$-groupoid)
\[
X_{\DR}(R) = X(R_{\sf red})
\]
with $R_{\sf red} \simeq \pi_0(R)/{\sf nilpotents}$ the {\em reduced classical ring} of~$R$.
In particular, $X_{\DR}$ is insensitive to the derived or nilpotent directions of~$X$.

Now view $X$ as a family over its de Rham stack $q\colon X \rightarrow X_{\DR}$. 
A $\Bbbk$-point 
\[\Spec(\Bbbk) \xrightarrow[]{x} X_{\DR}\]
corresponds to a reduced point
\[
\Spec(\Bbbk) \xrightarrow[]{\overline{x}} X ~,
\]
and so the fiber of $q$ at $x$ is the formal completion $\hat{X}_x$. Viewing $\hat{X}_x$ as a formal moduli problem over and under $\Spec(\Bbbk)$, 
one can describe functions on $\hat{X}_x$ as the Lie algebra cohomology
\[
\cO(\hat{X}_x) \simeq \sC^{\ast}_{\Lie}\bigl(\TT\hat{X}_x[-1]\bigr)
\]
where $\TT\hat{X}_x[-1]$ is the dg Lie algebra which is the tangent space of the based loop space at $x$ of the formal moduli problem~$\hat{X}_x$. 
The algebra of functions $\cO(\hat{X}_x)$ admits extra structure of a graded mixed algebra, as we discuss later. 

As the point $x$ varies, the various Chevalley--Eilenberg complexes $\sC^{\ast}_{\Lie}(\TT\hat{X}_x[-1])$ glue to form a global object over $X_{\DR}$ known as the presheaf $\cP_X$ of {\em mixed principal parts}, 
a presheaf of graded mixed algebras on $X_{\DR}$ constructed by~\cite{CPTVV}. 
We now sketch its definition.

Given a point $f\colon \Spec(R) \to X_{\DR}$,
corresponding to a reduced point $\overline{f}\colon \Spec(R_{\sf red}) \to X$, the value of the presheaf $\cP_X$ on the point $f$ is roughly given by the relative derived de Rham cohomology of the map $\Spec(R_{\sf red}) \rightarrow \Spec(R)\times_{X_{\DR}}X$, 
completed with respect to its Hodge filtration. There is a realization of the object $\cP_X$, denoted $|\cP_X|$, that recovers what is usually called principal parts. More precisely, there is an equivalence of sheaves over $X_{\DR}$
\[|\cP_X| \simeq q_*\cO_X~,\]
where $q\colon X \ra X_{\DR}$ is the canonical map.

This presheaf of mixed principal parts $\cP_X$ is a commutative algebra over another presheaf $\cD_{X_{\DR}}$, 
which can be thought of as a graded mixed enhancement of the structure sheaf of~$X_{\DR}$. Again, the realization of $\cD_{X_{\DR}}$ recovers the structure sheaf of $X_{\DR}$:
\[|\cD_{X_{\DR}}| \simeq \cO_{X_{\DR}}~.\]
In terms of its $\cD_{X_{\DR}}$-linear structure, the presheaf $\cP_X$ can be understood as sections of the infinite jet bundle with its canonical connection.

We now summarize our key technical result: 
\begin{quote}
principal parts on $X^M$ can be computed via factorization homology over~$M$.
\end{quote}  
Here $M$ is a smooth $d$-manifold, viewed as a derived stack (i.e., its Betti stack),
while factorization homology over $M$ is defined as
\[
\int_M\cA := \underset{U\in\Disk_M}\colim \cA(U)
\]
for $\cA$ a symmetric monoidal functor from the symmetric monoidal category $\Disk_M$ generated by framed $n$-disks embedded in $M$ into a nice enough symmetric monoidal category~$\cV$. 
(See \cite{oldfact} for a complete development of this notion.)

To be more precise, 
the presheaf of mixed principal parts $\cP_X$ can be enhanced to a presheaf of $\Disk_M$-algebras~$\cA_X$, called the algebra of {\em local principal parts}.  
Then we prove the following.

\begin{prop}
Let $M$ be a smooth $n$-manifold with the homotopy type of a finite CW complex, 
and let $X$ be a derived Artin stack, 
locally of finite presentation over~$\Bbbk$. 
There is an equivalence
\[\cP_{X^M} \simeq \int_M\cA_X\]
of modules over~$\cD_{X^M_{\DR}}$.
\end{prop}

\begin{remark}
    One must work with the graded mixed objects of \cite{CPTVV} to obtain the proposition above. See Remark \ref{failureforordinaryprincipalparts} for a more detailed discussion. 
~\hfill$\Diamond$\end{remark}

Therefore one could deform the presheaf of mixed principal parts $\cP_{X^M}$ --- so to speak, quantize $X^M$ --- by finding a deformation $\cA_X^{\pi}[\![\hbar]\!]$ of $\cA_X$ and defining
\[
\cP^{\pi}_{X^M}[\![\hbar]\!] := \int_M\cA_X^{\pi}[\![\hbar]\!]\]
as the \emph{quantum} principal parts on the mapping space. The theory of shifted Poisson structures as developed in \cite{CPTVV} provides a convenient framework for producing deformations of $\cA_X$. We have the following result.

\begin{prop} 
Let $X$ be an Artin derived stack locally of finite presentation over $\Bbbk$. 
Given an $n$-shifted Poisson structure $\pi \in {\sf Poiss}(X, n)$ with $n\geq1$, a $d$-manifold $M$ where $d \leq n+1$, and a framing of $M\times \RR^{n+1-d}$
we obtain an $\EE_{n+1-d}$-quantization of the mapping stack $\Map(\Pi_\infty M, X)$ by 
\[
\cP_{X^M}^{\pi}[\![\hbar]\!] := \int_M \cA_X^{\pi}[\![\hbar]\!]  
\]
\end{prop}

The key inputs are results of \cite{CPTVV} and \cite{Tom} providing equivalences between the space of shifted Poisson structures on a derived stack $X$ and shifted Poisson algebra structures on $\cP_X$. 

\begin{remark}
    \cite{CPTVV} relate shifted Poisson structures on $X$ with shifted Poisson algebra structures on a certain ind-object called $\cP_X(\infty)$. This framework is not compatible with our calculations. While we were writing this paper, \cite{Tom} provided a different approach to producing shifted Poisson algebra structures on $\cP_X$, which ended up being more convenient for our work. See Remark \ref{rmk on CPTVV poisson issues} for more details. 
~\hfill$\Diamond$\end{remark}

\begin{remark}
    The proposition above should be thought of as a version of the AKSZ construction. Given an $n$-shifted Poisson structure on $X$, we quantize the mapping stack $X^M$. \cite{Tom} produces a similar result: given an $n$-shifted Poisson structure on $X$, one obtains a shifted Poisson structure on the mapping stack $X^M$. It would be interesting to compare the quantizations obtained via the two approaches above. 
~\hfill$\Diamond$\end{remark}
Passing to realizations and taking global sections over $X^M_{\DR}$, we also obtain quantizations of global functions on mapping stacks, which is a main result of the present work.

\begin{theorem}
Let $X$ be an Artin derived stack locally of finite presentation over $\Bbbk$. 
Given an $n$-shifted Poisson structure $\pi \in {\sf Poiss}(X, n)$ with $n\geq1$, and a framed $d$-manifold $M$ where $d \leq n+1$, 
we obtain an $\EE_{n+1-d}$-quantization of the algebra $\Gamma(X^M, \cO_{X^M})$: 
\[Z^\hbar_\pi(M):= \Gamma\Bigl(X^M_{\DR}, \Bigl\vert\int_M\cA_X^{\pi}[\![\hbar]\!]\Bigr\vert \Bigr)~.\]
\end{theorem}

Furthermore, as \cite{CPTVV} identify the symmetric monoidal category $\Perf(X)$ with a certain full subcategory of modules over $\cP_X$, 
which we will denote $\Perf^{\epsilon, \gr}_{\cP_X}$, 
an $\EE_n$-deformation of $\cP_X$ gives rise to an $\EE_{n-1}$-deformation of the category~$\Perf(X)$. 
Applying factorization homology,
our result quantizes $\Perf(X^M)$ to a category of modules over~$\cP^\pi_{X^M}[\![\hbar]\!]$. 

As a special case of strong interest,
consider the derived character stack 
\[
\Loc_{\GG}(M):= \Map(M, \sB\GG)
\]
for $n =\dim(M) \leq 3$ and $\GG$ an affine algebraic group over $\Bbbk$. 
Let $\gg$ denote its Lie algebra. 
Given a level $\lambda \in \Sym^2(\gg)^{\GG}$, 
we construct a presheaf of $\Disk_M$-algebras $\cA_{\sB\GG}^{\lambda}$ on $\Loc_{\GG}(M)_{\DR}$ deforming $\cA_{\sB\GG}$, 
and we obtain the following result. 

\begin{prop}
Given a level $\lambda \in \Sym^2(\gg)^{\GG}$, there is an $\EE_{3-d}$-deformation  
\[\cP^{\lambda}_{\Loc_{\GG}(M)}[\![\hbar]\!] := \int_M \cA^{\lambda}_{\sB\GG}[\![\hbar]\!]\]
of the presheaf of mixed principal parts.
\end{prop}

The space $\Sym^2(\gg)^{\GG}$ is precisely the space of $2$-shifted Poisson structures on $\sB\GG$ in the sense of~\cite{CPTVV}. 
Therefore given a $2$-shifted Poisson structure on $\sB\GG$, we obtain quantizations of~$\Loc_{\GG}(M)$. 
These algebras might be interpreted as the quantum observables of Chern--Simons theory at level~$\lambda$.

If $n\leq 2$, we also obtain quantizations of perfect complexes on the derived character stack. 

\begin{prop}
Given $\lambda \in \Sym^2(\gg)^{\GG}$, there is an $\EE_{2-d}$-deformation
\[\Perf^{\lambda}(\Loc_{\GG}(M)):= \Perf^{\epsilon, \gr}_{\cP^{\lambda}_{\Loc_{\GG}(M)}}\]
of the category~$\Perf(\Loc_{\GG}(M))$.
\end{prop}

Note that this result is similar to the results in~\cite{BBJ} and \cite{Karlsson2024DeformationQV}. 
They construct deformations of $\qc(\Loc_{\GG}(M))$ as 
\[\qc_q(\Loc_{\GG}(M)) := \int_M \Rep_q(\GG),\]
where $\Rep_q(\GG)$ is an $\EE_2$-deformation of the category of representations of $\GG$, i.e., a quantum group. 

\subsection{Overview of the paper}
\begin{itemize}
    \item In {\bf Section~\ref{formal localization}}, we review background about de Rham stacks and principal parts from \cite{CPTVV}. 
    \item In {\bf Section~\ref{Deformation Quantization}}, we study the problem of deformation quantization of derived Artin stacks. We collect results on the various quantum objects that one can associate to a derived stack equipped with a shifted Poisson structure. 
    \item In {\bf Section~\ref{sec.principal.parts}}, we analyze principal parts on Betti mapping stacks. This is the technical heart of the paper. 
    \item In {\bf Section~\ref{factorization homology}}, we review the theory of factorization homology and perform calculations that are useful within the context of this article. 
    \item In {\bf Section~\ref{Principal Parts on Mapping Stacks as Factorization Homology}}, we show that principal parts on Betti mapping stacks can be computed as factorization homology of a $\Disk_M$-algebra.
    \item In {\bf Section~\ref{Quantizing Mapping Stacks}}, we prove our main quantization result: we quantize Betti mapping stacks with shifted Poisson target. This is a version of the AKSZ construction. 
    \item In {\bf Section~\ref{Quantizing BG and Loc}}, we exhibit quantizations of the derived stack of local systems $\Loc_{\GG}(M)$ on low-dimensional manifolds using a more explicit approach.
\end{itemize}
\subsection{Conventions}

\begin{itemize}
    \item $\Bbbk$ is a field of characteristic~0. 
    \item The word category will mean $(\infty, 1)$-category. E.g., see \cite{HTT}. 
    \item $\Spaces$ is the category of spaces (or $\infty$-groupoids or anima). 
    \item $\Mod_\Bbbk$ is the category of modules over $\Bbbk$. 
    \item $\CAlg(\cC)$ is the category of commutative algebras in a symmetric monoidal category~$\cC$.
    \item $\CAlg^{\geq 0}$ is the category of connective almost finitely presented $\Bbbk$-algebras. 
    \item ${\sf Aff}$ is the opposite category of $\CAlg^{\geq 0}$, endowed with the \'etale topology (\cite{TV} Definition~2.2.2.3). 
    Whenever we write $\Spec(R)$, it is assumed that $\Spec(R)$ is an object of ${\sf Aff}$ (i.e., $R$ is an almost finitely presented connective $\Bbbk$-algebra).
    \item A derived stack is an \'etale sheaf on ${\sf Aff}$. 
    $\dSt$ is the full subcategory of $\Fun(\CAlg^{\geq 0},\Spaces)$ consisting of derived stacks.
\end{itemize}
Unless an exception is stated locally, we use the following symbols as follows:
 the text
\begin{itemize}
\item $X$ is a derived Artin stack locally of finite presentation over $\Bbbk$;
\item $M$ is a smooth $d$-manifold with the homotopy type of a finite CW complex;
\item $\Pi_\infty M$ to denote the underlying $\infty$-groupoid of $M$ or the derived stack it determines (known as its Betti stack); and
\item $X^M$ is the derived stack~$\Map_{\dSt}(\Pi_\infty M, X)$.
\end{itemize}

\subsection{Acknowledgements}
RK thanks Aleksandr Karapetyan and Preston Cranford for several helpful discussions. We thank Damien Calaque, Sam Gunningham, Pavel Safronov, Bertrand Toën, and Gabriele Vezzosi for their comments and encouragement on a preliminary draft of this paper.

\subsection{AI statement} After a draft of this paper was completed, AI was used to generate a referee report, on the basis of which we made a number of revisions. None of the paper was written by AI.

\section{Formal Localization}\label{formal localization}

In this section we briefly review some constructions and results from~\cite{CPTVV}.

\subsection{Graded Mixed Complexes}

In this subsection, we will collect facts about graded mixed complexes. 
See \cite{CPTVV} section 1.1 or \cite{pav1} Section 1.1 for a more thorough treatment.

\begin{definition} \label{gradedcomplexes}
Let 
\[\Mod_{\Bbbk}^{\gr}:= \Fun(\ZZ^{\sf discrete}, \Mod_{\Bbbk})\]
denote the category of {\em graded complexes}, 
which we regard as a stable symmetric monoidal category via the Day convolution symmetric monoidal structure. 
The grading is called the {\em weight grading}.
\end{definition}

\begin{definition}
For each $n \in \ZZ$, the functor
\[(-)(n)\colon \Mod_{\Bbbk} \rightarrow \Mod_{\Bbbk}^{\gr}\]
sends an object $E\in \Mod_{\Bbbk}$ to the graded complex $E(n)$ consisting of $E$ in weight~$n$. 
\end{definition}

\begin{definition}
For each $n \in \ZZ$, the functor 
\[(-)_{n}\colon \Mod_{\Bbbk}^{\gr} \rightarrow \Mod_{\Bbbk}\]
sends a graded complex $E$ to its weight $n$ piece~$E_n$. 
\end{definition}

The forgetful functor
\[\CAlg(\Mod_{\Bbbk}^{\gr}) \longrightarrow \Mod_{\Bbbk}^{\gr}\]
has a left adjoint given by the free commutative algebra construction, 
denoted
\[\Sym(-)\colon \Mod_{\Bbbk}^{\gr} \longrightarrow \CAlg(\Mod_{\Bbbk}^{\gr}).\]
Given an object $V$ of $\Mod_{\Bbbk}$, there is an equivalence of objects in~$\Mod_{\Bbbk}$
\[\Sym\bigl(V(1)\bigr)_k \simeq \Sym^k(V) \simeq (V^{\otimes k})_{\Sigma_k}\]
where $\Sym\bigl(V(1)\bigr)_k $ denotes the weight $k$-piece of the graded object~$\Sym\bigl(V(1)\bigr)$.

\begin{definition}
The category of {\em graded mixed complexes} is 
\[\Mod^{\epsilon, \gr}_{\Bbbk} := \Mod_{\Bbbk[\epsilon]}(\Mod_{\Bbbk}^{\gr}),\]
where $\epsilon$ has cohomological degree $1$ and weight $1$. The algebra $\Bbbk[\epsilon]$ is a Hopf algebra: its couinit is given by $\epsilon \mapsto 0$ and the comultiplication is given by
\[\epsilon \mapsto \epsilon\otimes1 + 1\otimes \epsilon.\]
We regard $\Mod^{\epsilon, \gr}_{\Bbbk}$ as a symmetric monoidal category where the monoidal structure comes from the Hopf algebra structure on $\Bbbk[\epsilon]$. 
\end{definition}

An object $E$ of $\Mod^{\epsilon, \gr}_{\Bbbk}$ can be described as a collection of $\Bbbk$-modules $\{E_n\}_{n \in \ZZ}$, 
along with a family of morphisms $\{\epsilon\colon E_n \rightarrow E_{n+1}[1]\}_{n}$ such that $\epsilon^2 =0$. The collection of maps $\epsilon$ is called the {\em mixed structure} on~$E$. 

\begin{definition}
The functor
\[{\sf fgt}_{\epsilon}\colon \Mod_{\Bbbk}^{\epsilon, \gr} \rightarrow \Mod_{\Bbbk}^{\gr}\]
forgets the mixed structure on a graded mixed complex. 
\end{definition}
Note that the functor ${\sf fgt}_\epsilon$ is symmetric monoidal, colimit preserving, and creates equivalences.

\begin{remark}
Section 3.2 of \cite{Rak} constructs a symmetric monoidal equivalence of categories 
\[\widehat{{\sf Fil}}(\Mod_{\Bbbk}) \simeq \Mod_{\Bbbk}^{\epsilon, \gr}\]
where $\widehat{{\sf Fil}}(\Mod_{\Bbbk})$ denotes the category of complete filtered complexes. 
Under this equivalence, forgetting the mixed structure amounts to taking the associated graded object of a filtered object. 
All of the results in this paper can be stated in terms of filtered objects instead of graded mixed objects. 
We opt to use graded mixed objects to stay consistent with~\cite{CPTVV}. 
~\hfill$\Diamond$\end{remark}

\begin{definition}
    The functor 
    \[(-)((n))\colon \Mod^{\epsilon, \gr}_{\Bbbk} \ra \Mod^{\epsilon, \gr}_{\Bbbk}\]
    sends a graded mixed complex $E\in \Mod^{\epsilon, \gr}_{\Bbbk}$ to $E((n)):=E\otimes\Bbbk(-n)$. 
\end{definition}

 There is an adjunction given by the enrichment of $\Mod_{\Bbbk}^{\epsilon, \gr}$ over $\Mod_{\Bbbk}$
 \[\Bbbk(0)\otimes - \colon\Mod_{\Bbbk} \rightleftarrows \Mod_{\Bbbk}^{\epsilon, \gr} \colon \underline{\Hom}(\Bbbk(0), -).\]
 \begin{definition} \label{stdrel}
     The \emph{standard realization functor}
     \[|-|: \Mod_{\Bbbk}^{\epsilon, \gr} \longrightarrow \Mod_{\Bbbk}\]
     is defined to be the right adjoint in the adjunction above, i.e., $|-|:= \underline{\Hom}(\Bbbk(0), -)$.
 
 \end{definition}
One also obtains a standard realization functor on functor categories, which we also denote by the same symbol:

\[|- |: \Fun(\cC, \Mod^{\epsilon, \gr}_{\Bbbk}) \longrightarrow \Fun(\cC, \Mod_{\Bbbk})~.\]
In our examples of interest, $\cC$ will be the over-category $\bigl({\sf Aff}_{/X}\bigr)^{\sf op}$, where $X$ is a derived stack. 

The standard realization functor has the structure of a lax symmetric monoidal functor (as it is right adjoint to a symmetric monoidal functor), and therefore sends algebra objects to algebra objects. 

 \begin{remark}
     By \cite{CPTVV} Proposition 1.5.1, given an object $E \in \Mod_{\Bbbk}^{\epsilon, \gr}$, there is an equivalence in $\Mod_{\Bbbk}$
     \[|E| \simeq \prod_{n \geq 0}E_n\]
     where $\prod_{n \geq 0}E_n$ is given the total differential, i.e., the sum of the mixed structure and the internal differential. If one views $E$ as a double complex, then taking standard realization amounts to taking {\em half} of the total complex associated to the double complex (observe that only the non-negative weights appear). 
 ~\hfill$\Diamond$\end{remark}
\begin{remark}
    Given an object $E\in \Mod_{\Bbbk}^{\epsilon, \gr}$ that lives in nonnegative weights, i.e., $E_k\simeq 0$ for $k < 0$, we have 
    \[|E((n))|\simeq |E|~\]
    for all $n\leq 0$
~\hfill$\Diamond$\end{remark}

\subsection{De Rham Theory}

In this subsection, we will review de Rham theory as developed in \cite{CPTVV}. 
For $R$ a connective $\Bbbk$-algebra, 
let $\Mod_{R}^{\epsilon, \gr}$ denote the category of $R$-modules in $\Mod_{\Bbbk}^{\epsilon, \gr}$. 
Consider the functor 
\[(-)_0: \CAlg(\Mod_{R}^{\epsilon,\gr}) \rightarrow \CAlg(\Mod_R)\]
sending a graded mixed $R$-algebra $S$ to its weight $0$ part. 
This functor admits a left adjoint (see \cite{CPTVV} Proposition 1.3.8), 
which leads to the following definition. 

\begin{definition}
The \emph{de Rham algebra} functor  
\[\DR(-/R)\colon \CAlg(\Mod_R) \rightarrow \CAlg(\Mod_{R}^{\epsilon,\gr})\]
is the left adjoint to $(-)_0: \CAlg(\Mod_{R}^{\epsilon,\gr}) \rightarrow \CAlg(\Mod_R)$.
\end{definition}

For an $R$-algebra $S$, the de Rham algebra $\DR(S/R)$ can be thought of as a version of Illusie's Hodge-completed derived de Rham cohomology of the $R$-algebra $S$ (see \cite{Bh}).
There is an equivalence of graded algebras 
\[{\sf fgt}_{\epsilon}(\DR(S/R))\simeq \Sym_S(\LL_{S| R}[-1](1))~.\] 
The mixed structure can be thought of as the de Rham differential. 

In what follows, let $\CAlg^{\geq 0}$ denote the category of connective almost finitely presented $\Bbbk$-algebras. 

\begin{definition} \label{gradedmixedfunctions}
The functor
\[\DD: \CAlg^{\geq 0} \rightarrow \CAlg(\Mod_{\Bbbk}^{\epsilon, \gr})\]
sends a connective $\Bbbk$-algebra $R$ to the de Rham algebra $\DR(R_{\sf red}/R)$, 
where $R_{\sf red} \simeq \pi_0(R)/{\sf nilpotents}$. The functor 
\[\DD: \dSt^{\op} \longrightarrow \CAlg(\Mod_{\Bbbk}^{\epsilon,gr})\]
is the right Kan extension, taking values
\[\DD(X) \simeq \lim_{\Spec(R) \to X}\DD(R)\]
for each derived stack~$X$.
\end{definition}
Given a derived stack $X$, one could think of the algebra $\DD(X)$ as the relative de Rham cohomology of the map $X_{\sf red} \ra X$. Note that the relative de Rham cohomology of a closed immersion is functions on the corresponding formal completion (see \cite{Bh}). As the formal completion of $X$ along $X_{\sf red} \ra X$ is $X$, one could think of $\DD(X)$ as a graded mixed enhancement of functions on $X$. The source of the weight grading on $\DD(X)$ could be thought of as the Hodge filtration on relative de Rham cohomology: interesting weights in $\DD(X)$ appear when the relative cotangent complex of the embedding $X_{\sf red}\ra X$ is nontrivial. 
\begin{example}
Consider the $\Bbbk$-algebra $\Bbbk[x]$ where $|x| = n$. 
For $n \neq 0$, then $(\Bbbk[x])_{\sf red} \simeq \Bbbk$, 
and we have an equivalence of (non-mixed) graded algebras
\[{\sf fgt}_{\epsilon}\bigl(\DD(\Bbbk[x])\bigr) \simeq \Sym_{\Bbbk}\bigl(\LL_{\Bbbk|\Bbbk[x]}[-1](1)\bigr) \simeq \Sym_{\Bbbk}\bigl(\Bbbk[-n](1)\bigr).\]
Therefore $\DD(\Bbbk[x])$ looks like a graded polynomial algebra in one variable, where the variable has cohomological degree $n$ and weight~$1$. Note that we have an equivalence of commutative algebras 
\[\Bbbk[x] \simeq \prod_{k\geq 0}\Sym_{\Bbbk}(\Bbbk[-n](1))_k~,\]
i.e., $\DD(\Bbbk[x])$ looks like $\Bbbk[x]$ with an additional grading.

For $n = 0$, $(\Bbbk[x])_{\sf red} \simeq \Bbbk[x]$, 
and therefore we have 
\[{\sf fgt}_{\epsilon}\bigl(\DD(\Bbbk[x])\bigr) \simeq \Sym_{\Bbbk[x]}\bigl(\LL_{\Bbbk[x]|\Bbbk[x]}[-1](1)\bigr) \simeq \Sym_{\Bbbk[x]}(0) \simeq \Bbbk[x].\]
Observe that the weight grading is {\em trivial} when $\Bbbk[x]$ is reduced.~\hfill$\Diamond$ 
\end{example}

The notion of $\DD(X)$ being a graded mixed enhancement of global functions on $X$ can be made precise in the following manner. Assume $X = \Spec(R)$ is a derived affine. \cite{CPTVV} Lemma 2.2.4 shows that there is an equivalence in $\CAlg_{\Bbbk}$:
\[R \simeq |\DD(X)|\]
where $|\DD(X)|$ denotes the standard realization functor of Definition \ref{stdrel}. Now for $X$ a general derived stack (not necessarily affine), we have 
\[\begin{split}
    |\DD(X)| &\simeq |\lim_{\Spec(R) \ra X}\DD(\Spec(R))| \\
    & \simeq \lim_{\Spec(R) \ra X}|\DD(\Spec(R))| \\
    & \simeq \lim_{\Spec(R)\ra X} R \\
    & \simeq \Gamma(X, \cO_X)~,
\end{split}\]
where the second equivalence uses the fact that standard realization is a right adjoint.

\subsection{Principal Parts and Formal Localization}

We will now review the formal localization machinery of \cite{CPTVV} Section~2. 
Let $X$ be a derived Artin stack, locally of finite presentation over $\Bbbk$. The hypotheses on $X$ ensure that the results of \cite{CPTVV} apply. Consider the de Rham stack $X_{\DR}$ (see \cite{CPTVV} Definition 2.1.3). (Also see \cite{Lu5} for an account of de Rham stacks and their relation to the theory of $D$-modules.) As a functor of points the de Rham stack is given by $X_{\DR}(R) \simeq X(R_{\sf red})$ for $R\in \CAlg_{\Bbbk}^{\geq 0}$. There is a canonical map $q\colon X \ra X_{\DR}$, exhibiting $X$ as a family over $X_{\DR}$. 

Given a map $f\colon\Spec(R) \ra X_{\DR}$, we have a diagram 
\[\begin{tikzcd}
{\Spec(R_{\sf red})} \ar[rrd, "\overline{f}"] \ar[rd] \ar[ddr] && \\
& {X_R} \arrow[rd, phantom, "\lrcorner", very near start] \ar[r] \ar[d] & X \ar[d, "q"] \\
& {\Spec(R)} \ar[r, "f"] & {X_{\DR}}
\end{tikzcd}\]
where the square is a pullback square so
$X_R:= X \times_{X_{\DR}}\Spec(R)$,
a stack that which lives between $\Spec(R_{\sf red})$ and~$\Spec(R)$.
We will use items from this diagram repeatedly below.

Consider the following presheaves of graded mixed algebras on $X_{\DR}$. 

\begin{definition}
The \emph{crystalline structure sheaf} is the presheaf 
\[
\begin{array}{ccccc}
\cD_{X_{\DR}}\colon & ({\sf Aff}_{/X_{\DR}})^{\sf{op}} & \longrightarrow & \CAlg(\Mod_{\Bbbk}^{\epsilon, \gr})\\
& (\Spec(R) \xrightarrow[]{f} X_{\DR}) & \longmapsto & \DD(R)
\end{array}
\]
The \emph{presheaf of mixed principal parts} is the presheaf 
\[
\begin{array}{ccccc}
\cP_{X}\colon & ({\sf Aff}_{/X_{\DR}})^{\sf{op}} &\longrightarrow& \CAlg(\Mod_{\Bbbk}^{\epsilon, \gr})\\
&(\Spec(R) \xrightarrow[]{f} X_{\DR}) & \longmapsto &\cP_X(f)\simeq\DD(X_R)
\end{array}
\]
where $X_R:= X \times_{X_{\DR}}\Spec(R)$.
\end{definition}

Note that for every map $f\colon\Spec(R) \ra X_{\DR}$, 
the canonical map $X_R \to \Spec(R)$ induces a map $\DD(R) \to \DD(X_R)$ 
and hence equips $\cP_X$ with a $\cD_{X_{\DR}}$-algebra structure.  

Section 2.4.2 of \cite{CPTVV} shows that there are equivalences of presheaves on $X_{\DR}$:
\[|\cD_{X_{\DR}}| \simeq \cO_{X_{\DR}}~,\]
and 
\[|\cP_X| \simeq q_*\cO_X~.\]
The notation $|-|$ denotes the standard realization functor of Definition \ref{stdrel}, but now applied to the functor category $\Fun\bigl(({\sf Aff}_{/X_{\DR}})^{\sf{op}}, \Mod^{\epsilon, \gr}_{\Bbbk}\bigr)$. The sheaf $q_*\cO_X$ is the usual sheaf of principal parts on $X$ along with its natural crystalline connection. Its value at a point $f\colon\Spec(R) \ra X_{\DR}$ is functions on $X_R$.

One can describe the underlying graded algebra of mixed principal parts rather explicitly. 

\begin{theorem}[\cite{CPTVV} Proposition 2.2.7]    \label{underlyinggradedofpp}
Given a map $f \colon \Spec(R) \ra X_{\DR}$, there is a natural equivalence  
\[
{\sf fgt}_{\epsilon}\bigl(\cP_X(f)\bigr) \simeq \Sym_{R_{\sf red}}\bigl(\LL_{R_{\sf red}| X_R}[-1](1)\bigr) 
\]
of graded $\Sym_{R_{\sf red}}\bigl(\LL_{R_{\sf red}| R}[-1](1)\bigr) $-algebras. 
\end{theorem} 

\begin{example}
When $R = \Bbbk$ and $X = \sB\GG$ is the classifying stack of an algebraic group over $\Bbbk$, 
this construction gives familiar results: 
given a map $x\colon \Spec (\Bbbk) \ra \sB\GG_{\DR}$, the fiber $\sB\GG_{\Bbbk} \simeq \Spec(x)\times_{\sB\GG_{\DR}}\sB\GG$ is just the formal completion of $\sB\GG$ along the induced map $x\colon \Spec(\Bbbk) \rightarrow \sB\GG$ 
and the relative cotangent complex is given by $\LL_{\Bbbk| \sB\GG_{\Bbbk}} \simeq {\sf cofib}(\gg^{\vee}[-1] \rightarrow 0) \simeq \gg^{\vee}$. 
Therefore, the underlying (non-mixed) graded algebra of $\cP_{\sB\GG}(x)$ is given by $\Sym_\Bbbk\bigl(\gg^{\vee}[-1](1)\bigr)$. 
In fact, as a graded mixed algebra, $\cP_{\sB\GG}(x)$ is a graded mixed version of the usual Chevalley--Eilenberg cochain complex of $\gg$ (see \cite{CPTVV} section~3.6.2).   
~\hfill$\Diamond$\end{example}

\begin{remark} 
The case of a $\Bbbk$-point for general $X$ is similar. 
Given a map $x\colon \Spec(\Bbbk) \ra X_{\DR}$, the fiber $X_x$ is the formal completion of $X$ along the induced map $x\colon \Spec(\Bbbk) \to X$, 
and the relative cotangent complex is now given by $\LL_{\Spec(\Bbbk)| X_x} \simeq \LL_{X,x}[1]$. 
Therefore, the underlying graded algebra of $\cP_X(x)$ is given by $\Sym_\Bbbk\bigl(\LL_{X,x}(1)\bigr)$. 
Let $\TT_{X,x} \simeq \LL_{X,x}^{\vee}$ denote the tangent complex at $x$. 
Viewing $X_x$ as a formal moduli problem under $\Spec(\Bbbk)$, one obtains a Lie algebra structure on $\TT_{X, x}[-1]$. As a graded mixed algebra $\cP_X(x)$ is a graded mixed version of the Chevalley--Eilenberg cochain complex of $\TT_{X,x}[-1]$. 

If we have $R\simeq R_{\sf red}$, then principal parts at an $R$-point $f\colon\Spec(R) \ra X_{\DR}$ can again be thought of as Chevalley--Eilenberg cochains of the Lie algebra 
(now in $\Mod_R$) 
corresponding to the formal moduli problem $\Spec(R) \rightarrow X_R \rightarrow \Spec(R)$. 
(Note that we have a map $\Spec(R) \rightarrow X_R$ as $R$ is reduced.)

In the general case, a map $f \colon \Spec(R) \ra X_{\DR}$ has a fiber $X_R$ that is a formal moduli problem under $\Spec(R_{\sf red})$ and over $\Spec(R)$. 
There is a corresponding Lie algebroid $\TT_{R_{\sf red}| X_R} \to \TT_{R_{\sf red}| R}$, 
and the mixed principal parts $\cP_X(f)$ can be thought of as the Chevalley--Eilenberg complex of this. 
See \cite{CaGri} Section 3.3 for an account of ideas related to this. 
~\hfill$\Diamond$\end{remark}

In some sense, the derived stack $X$ can be thought of as the relative spectrum of the $\cD_{X_{\DR}}$-algebra $\cP_X$ . This is the formal localization result of \cite{CPTVV}: perfect complexes on $X$ are the same data as perfect modules over $\cP_X$. Note that $\cP_X$ is a commutative algebra in the symmetric monoidal category of presheaves of graded mixed modules on ${\sf dAff_{/X_{\DR}}}$, and therefore we obtain a symmetric monoidal category of modules $\Mod_{\cP_X}$. 

\begin{definition}
The category of \emph{perfect $\cP_X$-modules} $\Perf^{\epsilon, \gr}_{\cP_X}$ is the full subcategory of $\Mod^{\epsilon, \gr}_{\cP_X}$ consisting of objects $E$ satisfying:
\begin{itemize}

\item For all $f \colon \Spec(R) \ra X_{\DR}$, there exists an equivalence (noncanonical) of (non-mixed) graded $\cP_X(f)$-modules 
\[
E(f) \simeq \cP_X(f) \underset{R_{\sf red}}\otimes E_0
\]
for some $E_0 \in \Perf_{R_{\sf red}}$;

\item given maps 
\[\Spec(S) \xrightarrow[]{g} \Spec(R) \xrightarrow[]{f} X_{\DR}~,\]
the natural map 
\[
E(f) \underset{\cP_X(f)}\otimes \cP_X(f\circ g) \rightarrow E(f \circ g)
\]
is an equivalence.
\end{itemize}
\end{definition}

\begin{theorem}[Corollary 2.4.7~\cite{CPTVV}]
\label{formloc}
There is a symmetric monoidal equivalence of categories 
\[
\Perf(X) \simeq \Perf^{\epsilon, \gr}_{\cP_X}.
\]
\end{theorem}

\begin{remark}
    Informally, the functor 
    \[\Perf(X) \longrightarrow \Perf^{\epsilon, \gr}_{\cP_X}\]
    could be thought of as taking a vector bundle $\cV$ on $X$ to the jet bundle $J^{\infty}\cV$ equipped with its natural flat connection (thus defining an object over $X_{\DR}$). The definition of the category $\Perf^{\epsilon, \gr}_{\cP_X}$ precisely carves out the full subcategory of $\Mod^{\epsilon, \gr}_{\cP_X}$ spanned by objects that look like jet bundles of a vector bundle, making the functor above an equivalence. The functor in the other direction can be thought of as taking flat sections. Also see \cite{BK}.~\hfill$\Diamond$ 
\end{remark}

\section{Deformation Quantization}\label{Deformation Quantization}

Let $X$ be derived Artin, locally of finite presentation over $\Bbbk$.  In this section we discuss deformation quantizations of derived Artin stacks. For $n \geq 1$, given an $n$-shifted Poisson structure $X$ there are several quantum objects that one can construct:
\begin{itemize}
    \item a deformation of the presheaf of mixed principal parts~$\cP_X$, 
    \item a deformation of $\Gamma(X, \cO_X)$, the global functions on~$X$, 
    \item a deformation of the structure sheaf on the site of formally \'etale derived Artin stacks over~$X$,
    \item a deformation of the category of perfect complexes~$\Perf(X)$. 
\end{itemize}
In this section, we will construct deformations mentioned above. Note that there is similar discussion on deformations of principal parts and perfect complexes in \cite{CPTVV} Section 3.5. The constructions in this section are close cousins of the constructions in \cite{CPTVV}. Remark~\ref{rmk on CPTVV poisson issues} and Definitions~\ref{definitionofdeformation},\ref{definitionofdeformationquantization} highlight the difference between our constructions and those presented in~\cite{CPTVV}. 

\begin{remark}
Consider the mapping stack $\Loc_{\GG}(\Sigma)=\Map(\Sigma, \sB\GG)$ where $\Sigma$ is a compact orientable surface of genus $g\geq 2$. 
If $\GG$ is a connected semisimple algebraic group, then there is a dense open substack of irreducible $\GG$-local systems.  We will later construct quantizations of $\Map(\Sigma, \sB\GG)$, and in particular quantize the structure sheaf.  Therefore, a quantization of $\Loc_{\GG}(\Sigma)=\Map(\Sigma, \sB\GG)$ restricts to a quantization of the stack of irreducible local systems. 
There are various constructions of quantizations of this substack; see~\cite{Jordan} for a survey.
~\hfill$\Diamond$
\end{remark}

We will first describe our notion of deformation quantization. 
\subsection{Deformation Quantization of Graded Mixed Algebras}
Given an graded mixed algebra $P \in \CAlg\bigl(\Mod^{\epsilon, \gr}_{\Bbbk}\bigr)$, one can view it as an $\EE_{n+1}$-algebra and study deformations of it. Let $\hbar$ be a formal variable of weight $2$ and cohomlogical degree $0$.

\begin{definition}\label{definitionofdeformation}
Given a graded mixed algebra $P \in \CAlg(\Mod_{\Bbbk}^{\epsilon, \gr})$ and $n \geq -1$, an {\em $\EE_{n+1}$-deformation of $P$} is an object 
\[A \in \Alg_{\EE_{n+1}}\bigl(\Mod^{\epsilon, \gr}_{\Bbbk[\![\hbar]\!]}\bigr)~,\]
together with an equivalence of $\EE_{n+1}$-algebras 
\[A\otimes_{\Bbbk[\![\hbar]\!]}\Bbbk \simeq P\]
where $\Bbbk = \Bbbk[\![\hbar]\!]/(\hbar)$.
\end{definition}

\begin{remark}
    One could also consider the extension of scalars $-\otimes\Bbbk[\hbar]$ instead of $-\otimes\Bbbk[\![\hbar]\!]$ above. However, as $\hbar$ is of weight 2 (any positive weight works for the argument below), there is essentially no difference as $\Bbbk[\hbar]$ is already $\hbar$-adically complete, i.e., the natural map of graded mixed objects
    \[\Bbbk[\hbar] \longrightarrow \lim_k\Bbbk[\hbar]/(\hbar^k)\]
    is an equivalence. Indeed, as equivalences are detected weightwise, and as limits are also computed weightwise, we have
    \[\Bbbk\cdot\hbar^n \simeq \bigl(\Bbbk[\hbar]\bigr)_{2n} \simeq \bigl(\lim_k\Bbbk[\hbar]/(\hbar^k)\bigr)_{2n} \simeq \lim_k\bigl(\Bbbk[\hbar]/(\hbar^k)\bigr)_{2n} \simeq \Bbbk\cdot\hbar^{n}.\]
This is not true if $\hbar$ is of weight zero. 
    
    We opt to use the power series notation as the realization of $\Bbbk[\hbar]$ is given by 
    \[\begin{split}
        |\Bbbk[\hbar]| &\simeq |\lim_k\Bbbk[\hbar]/(\hbar^k)| \\
        &\simeq \lim_k|\Bbbk[\hbar]/(\hbar^k)| \\
        &\simeq \lim_k \Bbbk[\hbar]/(\hbar^k) \\
        &\simeq \Bbbk[\![\hbar]\!]
    \end{split}\]
    where $\Bbbk[\![\hbar]\!]$ on the right is a complex concentrated in degree~$0$.~\hfill$\Diamond$
\end{remark}

Such a deformation projects a Poisson shadow at the level of cohomology.
For a graded mixed complex $M$, write 
\[H^*(M) \colon = H^*\bigl({\sf fgt}_{\epsilon}(M)\bigr)~,\]
the cohomology with respect only to the internal differential.
Recall that $\EE_{n+1}(2) \simeq S^n$, 
and therefore, for $n\geq 1$, given an $\EE_{n+1}$-deformation $A$ of $P$, 
a choice of fundamental class of $S^n$ induces an operation (a shifted Poisson bracket)
\[\{-,-\}_A\colon H^*(A)^{\otimes 2} \ra H^*({}A)[-n]~.\]
By definition of an $\EE_{n+1}$-deformation, this operation is zero mod $\hbar$ and so we obtain an operation 
\[\{-,-\}_A\colon H^*(P)^{\otimes 2} \ra \hbar H^*(P)[-n]~.\]

\begin{definition}
For $n\geq 1$, 
the {\em semiclassical bracket} of an $\EE_{n+1}$-deformation $A$ of $P$ is the map 
\[\{-,-\}_A\colon H^*(P)^{\otimes 2} \ra \hbar H^*(P)[-n]~.\]
\end{definition}

The bracket encodes the first order information of the deformation,
and we will relate this to the $\PP_{n+1}$-algebra structures appearing via formality. 

Following \cite{Tom}, there is an operad $\PP_{\Bbbk[n+1](-2)}$ in $\Mod^{\epsilon, \gr}_{\Bbbk}$. An algebra over $\PP_{\Bbbk[n+1](-2)}$ is roughly the data of a graded mixed algebra $P \in \CAlg(\Mod_{\Bbbk}^{\epsilon, \gr})$ along with a compatible Poission bracket $\{-,-\}$ of cohomological degree $-n$ and weight $-2$. 

One can define the space of $\PP_{\Bbbk[n+1](-2)}$-algebra structures on it as the fiber of  
\[\Alg_{\PP_{\Bbbk[n+1](-2)}}(\Mod_{\Bbbk}^{\epsilon, \gr}) \rightarrow \CAlg(\Mod_{\Bbbk}^{\epsilon, \gr})\]
over $\{P\}$. We denote this space by $\PP_{n+1}-{\sf str}(P)$: the abuse of notation will be justified below. Given $\pi \in \PP_{n+1}-{\sf str}(P)$, we again obtain a map 
\[\hbar\cdot \{-,-\}_{\pi}\colon H^*(P)^{\otimes 2} \ra \hbar H^*(P)[-n]~.\]

Roughly speaking, we will define a deformation quantization of a Poisson algebra structure $\pi$ to be an $\EE_{n+1}$-deformation whose induced semiclassical bracket agrees with the map above. 

\begin{definition}\label{definitionofdeformationquantization}
    For $n\geq 1$, given a graded mixed algebra $P \in \CAlg(\Mod_{\Bbbk}^{\epsilon, \gr})$, and $\pi \in \PP_{n+1}-{\sf str}(P)$, a {\em deformation quantization of $\pi$} is the data of an $\EE_{n+1}$-deformation $A$ of $P$ satisfying 
    \[\{-,-\}_{A} = \hbar\cdot \{-,-\}_{\pi}\]
    as maps 
    \[H^*(P)^{\otimes 2} \ra \hbar H^*(P)[-n]~.\]
\end{definition}

\begin{remark}
    Note that the definition above merely asks for a comparison on cohomology. Relating an $\EE_{n+1}$-deformation to the full homotopy coherent $\PP_{n+1}$-algebra structure it quantizes (especially when $n=0,-1$) leads one to the notion of the Beilinson-Drinfeld operad. In this article, we will make extensive use of Dunn-Lurie additivity for the $\EE_n$-operads, and therefore we choose not to work with Beilinson-Drinfeld quantizations. We instead resort to various hacks to demonstrate that we produce nontrivial quantizations. See Remark~\ref{beilinsondrinfeld}.
~\hfill$\Diamond$\end{remark}

By formality theorems for the $\EE_{n+1}$-operad (for $n \geq 1$), we immediately obtain deformation quantizations of Poisson structures. 

\begin{lemma}\label{existence of deformation quantization}
    For $n \geq 1$, given given a graded mixed algebra $P \in \CAlg(\Mod_{\Bbbk}^{\epsilon, \gr})$, and $\pi \in \PP_{n+1}-{\sf str}(P)$, a deformation quantization of $\pi$ exists. 
\end{lemma}
\begin{proof}
    Let $\hbar$ be a formal variable of weight $2$ and cohomlogical degree $0$. Given $\pi \in \PP_{n+1}-{\sf str}(P)$, rescaling a bracket of arity $k$ by $\hbar^{k-1}$, we obtain a graded mixed $\PP_{n+1}$-algebra structure on $P[\![\hbar]\!]$. Now, by formality, $P[\![\hbar]\!]$ can be viewed as an $\EE_{n+1}$-algebra. Observe that this is an $\EE_{n+1}$-deformation of $P$ as the brackets on $P[\![\hbar]\!]$ are zero mod $\hbar$. Formality can be normalized so that the fundamental class of $S^1 \simeq \EE_{n+1}(2)$ acts by $\{-,-\}_{\pi}$. Therefore $P[\![\hbar]\!]$ is a deformation quantization of $\pi$. 
\end{proof}

\begin{remark}
    The Lemma above establishes the nontriviality of an $\EE_{n+1}$-deformation obtained via formality of the $\EE_{n+1}$-operad for $n\geq 1$. Later in this paper, we will construct deformations using factorization homology, sometime producing $\EE_0$ and $\EE_1$-deformations. We will establish nontriviality of the deformation using other methods. 
~\hfill$\Diamond$\end{remark}

We will now examine the semi-classical structure on $\EE_1$ and $\EE_0$-deformations. Note that these structures will only be used to detect {\em nontriviality} of the deformation: we do not introduce the notion of quantizing a Poisson algebra structure. 

The homology of the $\EE_1$-operad is the associative operad ${\sf Ass}$. The map of operads $\Lie \ra {\sf Ass}$, gives us the following definition. 
\begin{definition}
    Given an $\EE_1$-deformation $A$ of $P$, the {\em semiclassical structure} is the induced commutator operation on cohomology
    \[[-,-]_A\colon H^*(A)^{\otimes 2} \ra H^*(A)~.\]
\end{definition}
Again, for the trivial deformation, this is the zero bracket. 

We now treat the $\EE_0$-case. 
\begin{definition}
    Given an $\EE_0$-deformation $A$ of $P$, consider the $\hbar$-adic filtration
    \[F^{p}A \colon = \hbar^pA~.\]
    The {\em semiclassical structure} associated to the deformation $A$ is the $E_1$-page of the associated spectral sequence.
\end{definition}
Note that a nontrivial differential on the $E_1$-page detects nontriviality of the deformation. 

The results of this section also apply to presheaves of graded mixed complexes. 
\subsection{Deformation Quantization of Artin Stacks}

\cite{CPTVV} define the notion of a shifted Poisson structure on $X$. In particular, for $n \geq 0$, \cite{CPTVV} Definition 3.1.1 introduces a {\em space} of $n$-shifted Poisson structures ${\sf Poiss}(X, n)$. A point $\pi \in {\sf Poiss}(X, n)$ is roughly a derived version of a Poisson bivector. 
\begin{example}
(\cite{CPTVV} Section 3.1) 
Let $\GG$ be a reductive group over $\Bbbk$, and let $\gg$ be its Lie algebra. We have 
    \[
    \begin{split}
        \pi_0\bigl({\sf Poiss}(\sB\GG, 2)\bigr) &\simeq \Sym^2(\gg)^{\GG} \\
        \pi_0\bigl({\sf Poiss}(\sB\GG, 1)\bigr) &\simeq \wedge^3(\gg)^{\GG} \\
        \pi_0\bigl({\sf Poiss}(\sB\GG, n)\bigr) &\simeq \ast, n\geq 3.
    \end{split}
    \]
so that $B\GG$ has interesting Poisson structures in only a narrow range of dimensions.~\hfill$\Diamond$
\end{example}

Shifted Poisson structures on $X$ are closely related to {\em shifted Poisson algebra} structures on the algebra of principal parts~$\cP_X$. Following \cite{Tom}, the data of an $n$-shifted Poisson structure $\pi \in {\sf Poiss}(X,n)$ is roughly equivalent to the data of a $\cD_{X_{\DR}}$-linear $\PP_{n+1}$-algebra structure on $\cP_X$. 

\begin{theorem}[\cite{Tom} Proposition 3.8] 
Let $X$ be a derived Artin stack, locally of finite presentation over $\Bbbk$. There is an equivalence  
\[{\sf Poiss}(X, n) \simeq \Alg_{\PP_{\Bbbk[n+1](-2)}}(\Mod^{\epsilon, \gr}_{\cD_{X_{\DR}}})\underset{\CAlg(\Mod^{\epsilon, \gr}_{\cD_{X_{\DR}}})}\times\{\cP_X\}\]
in~$\spaces$.
\end{theorem}

\begin{remark}
\label{rmk on CPTVV poisson issues}
Note that the approach presented in \cite{Tom} is different from the approach presented in \cite{CPTVV}. \cite{CPTVV} relate shifted Poisson structures on $X$ with Poisson algebra structures on certain algebras living in $\ind(\Mod^{\epsilon, \gr}_{\cD_{X_{\DR}}})$. 
There is a functor 
\[
-\otimes\Bbbk(\infty):\Mod^{\epsilon, \gr}_{\cD_{X_{\DR}}} \longrightarrow \ind(\Mod^{\epsilon, \gr}_{\cD_{X_{\DR}}})~.
\]
With $\cP_X(\infty):= \cP_X\otimes\Bbbk(\infty)$,
\cite{CPTVV} show that there is an equivalence of spaces 
\[{\sf Poiss}(X, n) \simeq \PP_{n+1}{\rm -}{\sf structures}(\cP_X(\infty))~.
\]
In Section~\ref{sec.principal.parts}, we will express the algebra of principal parts $\cP_{X^M}$ on a mapping stack as a colimit. 
As the map to the ind-completion $-\otimes\Bbbk(\infty)$ is not colimit preserving in general, our approach is incompatible with the approach presented in \cite{CPTVV}. Therefore we opt to use the comparison theorem from \cite{Tom} instead, as it allows us to perform our computations in the category $\Mod^{\epsilon, \gr}_{\cD_{X_{\DR}}}$~\hfill$\Diamond$
\end{remark}

Therefore, given an $n$-shifted Poisson structure $\pi \in {\sf Poiss}(X, n)$, we obtain a $\PP_{\Bbbk[n+1](-2)}$-algebra structure on $\cP_X \in \CAlg(\Mod_{\cD_{X_{\DR}}}^{\epsilon, \gr})$.  We have the following lemma. 

\begin{lemma} \label{poissbivecttodef}
Let $X$ be a derived Artin stack, locally of finite presentation over $\Bbbk$. 
For each $n \geq 1$, given an $n$-shifted Poisson structure $\pi \in {\sf Poiss}(X, n)$, 
there is an $\EE_{n+1}$-deformation quantization (in the sense of Definition~\ref{definitionofdeformationquantization}) of principal parts of $X$
\[
\cP_X^{\pi}[\![\hbar]\!] \in \Alg_{\EE_{n+1}}(\Mod_{\cD_{X_{\DR}}[\![\hbar]\!]}^{\epsilon, \gr})
\]
\end{lemma}
\begin{proof}
    This is immediate from Lemma~\ref{existence of deformation quantization}. 
\end{proof}

\begin{remark}
    The $\EE_{n+1}$-algebra $\cP_X^{\pi}[\![\hbar]\!] \in \Alg_{\EE_{n+1}}(\Mod_{\cD_{X_{\DR}}[\![\hbar]\!]}^{\epsilon, \gr})$ appearing in the statement of Lemma~\ref{poissbivecttodef} has underlying $\cD_{X_{\DR}}$-module $\cP_X[\![\hbar]\!]$. 
~\hfill$\Diamond$\end{remark}

\begin{remark}\label{beilinsondrinfeld}
    Given a $\cD_{X_{\DR}}[\![\hbar]\!]$-linear $\PP_{n+1}$-algebra structure on $\cP_{X}[\![\hbar]\!]$, one could also consider the notion of a $\BB\DD_{n+1}$-quantization of principal parts. In our notation, this would mean a $\cD_{X_{\DR}}[\![\hbar]\!][\![u]\!]$-linear $\BB\DD_{n+1}$-algebra structure on $\cP_X[\![\hbar]\!][\![u]\!]$ such that 
    \[\cP_X[\![\hbar]\!][\![u]\!] \underset{\cD_{X_{\DR}}[\![\hbar]\!][\![u]\!]}{\otimes} \cD_{X_{\DR}}[\![\hbar]\!] \simeq \cP_X[\![\hbar]\!]\]
    as $\PP_{n+1}$-algebras, i.e., we keep track of the Poisson structure that is being quantized. For $n\geq 1$, there is a formality isomorphism $\BB\DD_{n+1}\simeq \PP_{n+1}$, therefore one also immediately obtains $\BB\DD_{n+1}$-quantizations starting from an $n$-shifted Poisson structure on $X$. We choose not to deal with $\BB\DD_{n+1}$-quantizations in this article, as the focus is on constructing quantizations of mapping stacks $X^M$ using factorization homology. Given an $n$-shifted Poisson structure on $X$ ($n\geq1$), and a $d$-manifold $M$ along with a framing of $M\times \RR^{n+1-d}$, our methods produce an $\EE_{n+1-d}$-quantization (for $n+1-d\geq 2$, this is in the sense of Definition~\ref{definitionofdeformationquantization}) of the Betti mapping stack $X^M$. Given an additivity theorem of the form 
    \[\Alg_{\EE_d}\bigl(\Alg_{\BB\DD_{n+1-d}}\bigr) \simeq \Alg_{\BB\DD_{n+1}}\]
    our methods would produce $\BB\DD_{n+1-d}$-quantizations of mapping stacks. 
~\hfill$\Diamond$\end{remark}

Consider the canonical map $q: X \ra X_{\DR}$. 
Recall that we have
\[|\cD_{X_{\DR}}| \simeq \cO_{X_{\DR}}\]
and 
\[|\cP_X| \simeq q_*\cO_X~.\]
The realization functor $|-|$ is, in fact, lax symmetric monoidal as it is right adjoint to a symmetric monoidal functor given by sending an ordinary complex to a graded mixed complex concentrated in weight~$0$.
Hence $\cD_{X_{\DR}}$-linear deformations of $\cP_X$ yield deformation quantizations of $q_*\cO_X$,
and taking sections of $q_*\cO_X$ over $X_{\DR}$, we also obtain deformations of the global functions~$\Gamma(X, \cO_X)$. 
We now pursue this line of thought.

\begin{lemma} \label{classicallimits}
Let $\cC$ and $\cD$ be presentably symmetric monoidal stable categories. 
Given an adjunction 
\[F\colon \cC \rightleftarrows \cD\colon G\]
where the left adjoint $F$ has the structure of a symmetric monoidal functor, $G$ admits the structure of a lax monoidal functor. Given a commutative algebra $R \in \CAlg(\cD)$, we obtain an adjunction 
\[F_R: \Mod_{G(R)}(\cC) \rightleftarrows \Mod_{R}(\cD): G_R~\]
where the left adjoint $F_R$ is symmetric monoidal and $G_R$ is lax monoidal. 
This further induces an adjunction
\[F_R: \Alg_{\EE_{n+1}}(\Mod_{G(R)}(\cC)) \rightleftarrows \Alg_{\EE_{n+1}}(\Mod_{R}(\cD)): G_R~.\]
Given $P, Q \in \Mod_{R}(\cD)$, we have the canonical maps (given by the lax monoidal structure)
\[\mu_{P,Q}: G_R(P)\underset{G_R(R)}\otimes G_R(Q) \longrightarrow G_R(P\otimes_RQ)~.\]
For a cofiber sequence in $\Mod_{R}(\cD)$
\[Q' \ra R \ra Q~,\]
if $\mu_{P, Q'}$ and  $\mu_{P, R}$ are equivalences, then so is $\mu_{P,Q}$. Further, if $P,Q \in \Alg_{\EE_{n+1}}(\Mod_{R}(\cD))$, then $\mu_{P,Q}$ is an equivalence of $\EE_{n+1}$-algebras. 
\end{lemma}

We postpone the proof of this lemma and deduce the deformation statements instead. The following asserts that an $\EE_{n+1}$-deformation of principal parts $\cP_X$ defines a deformation of the pushforward $q_\ast\cO_X$, where $q\colon X \ra X_{\DR}$ is the canonical map.

\begin{lemma} \label{quantizingordinaryprincipalparts}
   Consider the map $q\colon X \ra X_{\DR}$ for $X$ a derived Artin stack, locally of finite presentation over $\Bbbk$. Given an $\EE_{n+1}$-algebra
   \[\cP^{\pi}_X[\![\hbar]\!] \in \Alg_{\EE_{n+1}}(\Mod_{\cD_{X_{\DR}}[\![\hbar]\!]}^{\epsilon, \gr})\qquad {\textrm with }\qquad\cP_X^{\pi}[\![\hbar]\!] \otimes_{\cD_{X_{\DR}}[\![\hbar]\!]}{\cD_{X_{\DR}}}|_{\hbar=0} \simeq \cP_{X}\]
   one obtains an $\EE_{n+1}$-algebra 
   \[q_*\cO_X^{\pi}[\![\hbar]\!] \in \Mod_{\cO_{X_{\DR}}[\![\hbar]\!]}\]
   deforming the sheaf $q_*\cO_X$: i.e., we have an equivalence
   \[q_*\cO_X^{\pi}[\![\hbar]\!]\otimes_{\cO_{X_{\DR}}[\![\hbar]\!]}\cO_{X_{\DR}} \simeq q_*\cO_X ~.\]
   
\end{lemma}

\begin{proof}
We begin with some computations. We will first show that 
\[|\cD_{X_{\DR}}[\![\hbar]\!]| \simeq \cO_{X_{\DR}}[\![\hbar]\!]~.\]
as objects of $\Fun\bigl(({\sf Aff}_{/X_{\DR}})^{\op}, \Mod_{\Bbbk}\bigr)$. We compute
\[
\begin{split}
   |\cD_{X_{\DR}}[\![\hbar]\!]| &\simeq |\lim_k \cD_{X_{\DR}}[\hbar]/(\hbar^k)| \\
   &\simeq \lim_k|\cD_{X_{\DR}}[\hbar]/(\hbar^k)| \\
   &\simeq \lim_k|\cD_{X_{\DR}}|\otimes\Bbbk[\hbar]/(\hbar^k) \\
   &\simeq |\cD_{X_{\DR}}|[\![\hbar]\!] \simeq \cO_{X_{\DR}}[\![\hbar]\!]~.
\end{split}
\]
The second equivalence is due to the fact that realization is a right adjoint. The third is due to the fact that right adjoints between stable categories are exact, and the final equivalence follows from the fact that $|\cD_{X_{\DR}}|\simeq \cO_{X_{\DR}}$. 
We will now show that the lax monoidal structure provides equivalences
\[\mu_{\cP_X[\![\hbar]\!], \cD_{X_{\DR}}[\![\hbar]\!]}: |\cP^{\pi}_X[\![\hbar]\!]\otimes_{\cD_{X_{\DR}}[\![\hbar]\!]}\cD_{X_{\DR}}[\![\hbar]\!]| \xra{\simeq} q_*\cO_X[\![\hbar]\!]\]
       and 
       \[\mu_{\cP_X[\![\hbar]\!], \cD_{X_{\DR}}[\![\hbar]\!]((-2))}\colon |\cP^{\pi}_X[\![\hbar]\!]\otimes_{\cD_{X_{\DR}}[\![\hbar]\!]}\cD_{X_{\DR}}[\![\hbar]\!]((-2))| \simeq q_*\cO_X[\![\hbar]\!]\]
       in $\Mod_{\cO_{X_{\DR}}[\![\hbar]\!]}$.

For the first equivalence, note that the underlying module of $\cP^{\pi}_X[\![\hbar]\!]$ is merely $\cP_X[\![\hbar]\!]$. Therefore, as objects of $\Mod_{\cO_{X_{\DR}}[\![\hbar]\!]}$, we have 
\[\begin{split}
   |\cP^{\pi}_X[\![\hbar]\!]\otimes_{\cD_{X_{\DR}}[\![\hbar]\!]}\cD_{X_{\DR}}[\![\hbar]\!]|  &\simeq |\cP^{\pi}_X[\![\hbar]\!]| \\
   &\simeq |\cP_X[\![\hbar]\!]| \\
   &\simeq |\lim_k \cP_{X}[\hbar]/(\hbar^k)| \\
   &\simeq \lim_k|\cP_{X}[\hbar]/(\hbar^k)| \\
   &\simeq \lim_k|\cP_{X}|\otimes\Bbbk[\hbar]/(\hbar^k) \\
   &\simeq |\cP_{X}|[\![\hbar]\!] \simeq q_*\cO_{X_{\DR}}[\![\hbar]\!]
\end{split}
 \]
and therefore we have the first equivalence. 
For the second, note that the weight-shifted object 
\[\cP^{\pi}_X[\![\hbar]\!]\otimes_{\cD_{X_{\DR}}[\![\hbar]\!]}\cD_{X_{\DR}}[\![\hbar]\!]((-2)) \simeq \cP^{\pi}_X[\![\hbar]\!]((-2))\]
still lives in non-negative weights, and that shifting the weights up has no effect on the realization. Therefore we have the second equivalence by the same computation. 

Now the lemma follows immediately from Lemma \ref{classicallimits}.
   \begin{itemize}
       \item Set $\cC = \Fun\bigl(({\sf Aff}_{/X_{\DR}})^{\op}, \Mod_{\Bbbk}\bigr)$ and $\cD = \Fun\bigl(({\sf Aff}_{/X_{\DR}})^{\op}, \Mod_{\Bbbk}^{\epsilon, \gr}\bigr)$. The adjunction is given by 
       \[ F=(-)(0)\colon\Fun\bigl(({\sf Aff}_{/X_{\DR}})^{\op}, \Mod_{\Bbbk}\bigr) \rightleftarrows\Fun\bigl(({\sf Aff}_{/X_{\DR}})^{\op}, \Mod_{\Bbbk}^{\epsilon, \gr}\bigr)\colon |-|=G~.\]
       \item Set $R = \cD_{X_{\DR}}[\![\hbar]\!]$. We verified $G_R(R) \simeq |\cD_{X_{\DR}}[\![\hbar]\!]|\simeq  \cO_{X_{\DR}}[\![\hbar]\!]$.
       \item Set $P = \cP_X^{\pi}[\![\hbar]\!] \in \Alg_{\EE_{n+1}}\bigl(\Mod^{\epsilon, \gr}_{\cD_{X_{\DR}}[\![\hbar]\!]})$.
       \item Set $Q' = \cD_{X_{\DR}}[\![\hbar]\!]((-2))$ and $Q = \cD_{X_{\DR}}$. We have a cofiber sequence in $\Mod^{\epsilon, \gr}_{\cD_{X_{\DR}}[\![\hbar]\!]}$:
       \[\cD_{X_{\DR}}[\![\hbar]\!]((-2)) \xra{\cdot \hbar} \cD_{X_{\DR}}[\![\hbar]\!] \ra \cD_{X_{\DR}}~.\]
       \item We have verified that there are equivalences 
       \[\mu_{\cP_X[\![\hbar]\!], \cD_{X_{\DR}}[\![\hbar]\!]}: |\cP^{\pi}_X[\![\hbar]\!]\otimes_{\cD_{X_{\DR}}[\![\hbar]\!]}\cD_{X_{\DR}}[\![\hbar]\!]| \xra{\simeq} q_*\cO_X[\![\hbar]\!]\]
       and 
       \[\mu_{\cP_X[\![\hbar]\!], \cD_{X_{\DR}}[\![\hbar]\!]((-2))}\colon |\cP^{\pi}_X[\![\hbar]\!]\otimes_{\cD_{X_{\DR}}[\![\hbar]\!]}\cD_{X_{\DR}}[\![\hbar]\!]((-2))| \simeq q_*\cO_X[\![\hbar]\!]\]
       in $\Mod_{\cO_{X_{\DR}}[\![\hbar]\!]}$
       
   \end{itemize}
   Therefore by Lemma \ref{classicallimits}, we have an $\EE_{n+1}$-algebra 
       \[q_*\cO_X^{\pi}[\![\hbar]\!]:= |\cP_X^{\pi}[\![\hbar]\!]| \in \Alg_{\EE_{n+1}}\bigl(\Mod_{\cO_{{X_{\DR}}}[\![\hbar]\!]}\bigr)\]
       deforming $q_*\cO_X$: 
       \[q_*\cO_X^{\pi}[\![\hbar]\!]\otimes_{\cO_{X_{\DR}}[\![\hbar]\!]}\cO_{X_{\DR}} \simeq |\cP_X^{\pi}[\![\hbar]\!]|\otimes_{|\cD_{X_{\DR}}[\![\hbar]\!]|}|\cD_{X_{\DR}}| \simeq |\cP_X^{\pi}[\![\hbar]\!] \otimes_{\cD_{X_{\DR}}[\![\hbar]\!]}\cD_{X_{\DR}}| \simeq |\cP_X| \simeq q_*\cO_X~.\]

\end{proof}

The functor of taking global sections 
\[\Gamma(X_{\DR}, -)\colon \Fun\bigl(({\sf Aff}/X_{\DR})^{\op}, \Mod_{\Bbbk}\bigr) \longrightarrow \Mod_{\Bbbk}\]
again has the structure of a lax symmetric monoidal functor. It is right adjoint to the symmetric monoidal functor given by taking the constant presheaf. 
Note that $\Gamma(X_{\DR}, q_{*}\cO_X) \simeq \Gamma(X, \cO_X)$. Therefore, taking global sections $\Gamma(X_{\DR}, q_*\cO_X^{\pi}[\![\hbar]\!])$ provide deformations of global functions on $X$. 

\begin{lemma} \label{quantizingglobalfunctions}
 Let $X$ be a derived Artin stack, locally of finite presentation over $\Bbbk$. Consider the family $q\colon X \ra X_{\DR}$ Given an $\EE_{n+1}$-deformation of principal parts $\cP_X$, i.e., an object
   \[\cP^{\pi}_X[\![\hbar]\!] \in \Alg_{\EE_{n+1}}(\Mod_{\cD_{X_{\DR}}[\![\hbar]\!]}^{\epsilon, \gr})\]
   together with an equivalence 
   \[\cP^{\pi}_X[\![\hbar]\!] \otimes_{\cD_{X_{\DR}}[\![\hbar]\!]}{\cD_{X_{\DR}}}|_{\hbar=0} \simeq \cP_{X}~,\]
   one obtains an $\EE_{n+1}$-deformations of $\Gamma(X, \cO_X)$, i.e., 
   $\EE_{n+1}$-algebra $\Gamma(X, \cO_X)^{\pi}[\![\hbar]\!] \in \Mod_{\Gamma(X_{\DR}, \cO_{X_{\DR}}[\![\hbar]\!])}$ together with an equivalence
   \[\Gamma(X, \cO_X)^{\pi}[\![\hbar]\!] \otimes_{\Gamma(X_{\DR}, \cO_{X_{\DR}}[\![\hbar]\!])}{\Gamma(X_{\DR}, \cO_{X_{\DR}})}_{\hbar =0} \simeq \Gamma(X, \cO_X).\]   
\end{lemma}
\begin{proof}
The argument is very similar to Lemma \ref{quantizingordinaryprincipalparts}. We begin with some observations. We have 
\[\begin{split}
    \Gamma(X_{\DR}, \cO_{X_{\DR}}[\![\hbar]\!]) &\simeq \Gamma(X_{\DR}, \lim_k \cO_{X_{\DR}}[\hbar]/\hbar^k) \\
    &\simeq \lim_k\Gamma(X_{\DR}, \cO_{X_{\DR}}[\hbar]/\hbar^k) \\
    &\simeq \lim_k\Gamma(X_{\DR}, \cO_{X_{\DR}})\otimes\Bbbk[\hbar]/\hbar^k \\
    &\simeq \Gamma(X_{\DR}, \cO_{X_{\DR}})[\![\hbar]\!]~.
    \end{split}\]
The second equivalence is due to the fact that global sections is a right adjoint and the third is due to the fact that right adjoints between stable categories are exact.
Now we will show that we have an equivalence of objects of $\Mod_{\Gamma(X_{\DR}, \cO_{X_{\DR}})[\![\hbar]\!]}$
\[\Gamma(X_{\DR}, q_*\cO_X^{\pi}[\![\hbar]\!]) \simeq \Gamma(X_{\DR}, q_*\cO_X)[\![\hbar]\!]~.\]
Note that this is merely a computation involving the underlying module of $\Gamma(X_{\DR}, q_*\cO_X^{\pi}[\![\hbar]\!])$. We again have 
\[\begin{split}
    \Gamma(X_{\DR}, q_*\cO_X^{\pi}[\![\hbar]\!]) &\simeq \Gamma(X_{\DR},q_*\cO_X[\![\hbar]\!]) \\
    &\simeq \Gamma(X_{\DR}, \lim_k q_*\cO_{X}[\hbar]/\hbar^k) \\
    &\simeq \lim_k\Gamma(X_{\DR}, q_*\cO_{X}[\hbar]/\hbar^k) \\
    &\simeq \lim_k\Gamma(X_{\DR}, q_*\cO_{X})\otimes[\hbar]/\hbar^k \\
    &\simeq \Gamma(X_{\DR}, q_*\cO_{X_{\DR}})[\![\hbar]\!]~.
    \end{split}\]

Now the lemma follows immediately from Lemma \ref{classicallimits}. Indeed, we have the following. 
\begin{itemize}
    \item Set $\cC = \Mod_{\Bbbk}$ and $\cD = \Fun\bigl(({\sf Aff}_{/X_{\DR}})^{\op}, \Mod_{\Bbbk}\bigr)$. Consider the adjunction
    \[F=c\colon \Mod_{\Bbbk} \rightleftarrows \Fun\bigl(({\sf Aff}_{/X_{\DR}})^{\op}, \Mod_{\Bbbk}\bigr)\colon \Gamma(X_{\DR}, -)=G\]
    where $c$ denotes the constant presheaf functor. 
    \item Set $R = \cO_{X_{\DR}}[\![\hbar]\!]$. We verified that $G_R(R) = \Gamma(X_{\DR}, \cO_{X_{\DR}}[\![\hbar]\!]) \simeq \Gamma(X_{\DR}, \cO_{X_{\DR}})[\![\hbar]\!]~.$
    \item Set $P =  q_*\cO_X^{\pi}[\![\hbar]\!] \in \Alg_{\EE_{n+1}}(\Mod_{\cO_{X_{\DR}}[\![\hbar]\!]})$. 
    \item Set $Q' = \cO_{X_{\DR}}[\![\hbar]\!]$ and $Q = \cO_{X_{\DR}}$. We have a cofiber sequence in $\Mod_{\cO_{X_{\DR}}[\![\hbar]\!]}$
    \[\cO_{X_{\DR}}[\![\hbar]\!] \xra{\cdot \hbar}\cO_{X_{\DR}}[\![\hbar]\!] \ra \cO_{X_{\DR}}~.\]
    \item We verified that 
    \[\mu_{q_*\cO_X^{\pi}[\![\hbar]\!], \cO_{X_{\DR}}[\![\hbar]\!]}: \Gamma(X_{\DR}, q_*\cO_X^{\pi}[\![\hbar]\!]) \otimes_{\Gamma(X_{\DR}, \cO_{X_{\DR}}[\![\hbar]\!])}\Gamma(X_{\DR},\cO_{X_{\DR}}[\![\hbar]\!]) \xra{\simeq} \Gamma(X_{\DR}, q_*\cO_X)[\![\hbar]\!]\]
    is an equivalence in $\Mod_{\Gamma(X_{\DR}, \cO_{X_{\DR}}[\![\hbar]\!])}$. 
\end{itemize}
Therefore by Lemma \ref{classicallimits}, we have an $\EE_{n+1}$-algebra 
\[\Gamma(X, \cO_X)^{\pi}[\![\hbar]\!]:= \Gamma(X_{\DR}, q_*\cO_X^{\pi}[\![\hbar]\!]) \in \Alg_{\EE_{n+1}}\bigl(\Mod_{\Gamma(X_{\DR}, \cO_{X_{\DR}}[\![\hbar]\!])}\bigr)\]
deforming $\Gamma(X, \cO_X)$:
\[\begin{split}
    \Gamma(X, \cO_X)^{\pi}[\![\hbar]\!] \otimes_{\Gamma(X_{\DR}, \cO_{X_{\DR}}[\![\hbar]\!])}\Gamma(X_{\DR},\cO_{X_{\DR}}) &\simeq  \Gamma(X_{\DR}, q_*\cO_X^{\pi}[\![\hbar]\!]) \otimes_{\Gamma(X_{\DR}, \cO_{X_{\DR}}[\![\hbar]\!])}\Gamma(X_{\DR},\cO_{X_{\DR}}) \\
    &\simeq \Gamma(X_{\DR}, q_*\cO_X) \simeq \Gamma(X, \cO_X)~.
\end{split}\]
Therefore we have the result. 
\end{proof}

\begin{proof}[Proof of Lemma~\ref{classicallimits}]
    Given $P \in \Mod_R(\cD)$, consider the functors
    \[G_R(P)\otimes_{G_R(R)}G_R(-)\colon \Mod_R(\cD) \ra \Mod_{G_R(R)}(\cC)~,\]
    and
    \[G_R(P\otimes_R(-))\colon \Mod_R(\cD) \ra \Mod_{G_R(R)}(\cC)~.\]
    Both of these functors are exact. Indeed, since $G_R$ is a right adjoint, it preserves limits and is therefore exact (note that $\Mod_R(\cD)$ is stable). Tensor products preserve colimits and are again exact. 
    Now consider the cofiber sequence 
    \[Q' \ra R \ra Q\]
    from the statement of the lemma. Applying the natural map 
    \[\mu_{P, -}\colon G_R(P)\otimes_{G_R(R)}G_R(-) \ra G_R(P\otimes_R(-))\]
    to the cofiber sequence above, we obtain a diagram where the rows are cofiber sequences 
    \[
    \begin{tikzcd}
        G_R(P)\otimes_{G_R(R)}G_R(Q') \arrow[r] \arrow[d, "\mu_{P,Q'}"] &G_R(P)\otimes_{G_R(R)}G_R(R) \arrow[r] \arrow[d, "\mu_{P,R}"] &G_R(P)\otimes_{G_R(R)}G_R(Q) \arrow[d, "\mu_{P,Q}"] \\
        G_R(P\otimes_R(Q')) \arrow[r] &G_R(P\otimes_R(R)) \arrow[r] &G_R(P\otimes_R(Q))
    \end{tikzcd}
    \]
    As two of the downward arrows are equivalences by assumption, so is the third, i.e.
    \[\mu_{P,Q}: G_R(P)\otimes_{G_R(R)}G_R(Q) \ra G_R(P\otimes_RQ)\]
    is an equivalence. Now if $P,Q \in \Alg_{\EE_{n+1}}(\Mod_{R}(\cD))$, then $\mu_{P,Q}$ is a map of $\EE_{n+1}$-algebras. The forgetful functor 
    \[\Alg_{\EE_{n+1}}(\Mod_{R}(\cD)) \ra \Mod_R(\cD)\]
    detects equivalences. Therefore by the argument above, $\mu_{P,Q}$ is an equivalence of $\EE_{n+1}$-algebras, and we have the result. 
\end{proof}

Now we show that the deformations obtained above restrict along \'etale maps. Given a map of derived Artin stacks 
\[f\colon Y \ra X~,\]
we have a map of $\cD_{Y_{\DR}}$-algebras 
\[f^{-1}\cP_X \ra \cP_Y~.\]
This map is an equivalence if the map $f$ is formally \'etale.
\begin{lemma} \label{restriction}
    Given a formally \'etale map 
    \[f\colon Y \ra X~,\]
    the natural map 
    \[f^{-1}\cP_X \ra \cP_Y~\]
    is an equivalence of $\cD_{Y_{\DR}}$-algebras.
\end{lemma}
\begin{proof}
    Take a map $g: \Spec(R) \ra Y_{\DR}$. We will show that the map 
    \[f^{-1}\cP_X(g) \ra \cP_Y(g)~\]
    is an equivalence of $\cD_{Y_{\DR}}(g)$-algebras. it suffices to check that this map is an equivalence after forgetting the mixed structure. 
    We have the map of $\Sym_{R_{\sf red}}(\LL_{R_{\sf red}\vert}R[-1](1))\simeq {\sf fgt}_{\epsilon}(\cD_{Y_{\DR}}(g))$-algebras
    \[\Sym_{R_{\sf red}}(\LL_{R_{\sf red}\vert X_R}[-1](1))\simeq {\sf fgt}_{\epsilon}(f^{-1}\cP_X(g)) \ra {\sf fgt}_{\epsilon}(\cP_Y(g))\simeq \Sym_{R_{\sf red}}(\LL_{R_{\sf red}\vert Y_R}[-1](1))\]
    
    This map is induced by the canonical map of stacks over $\Spec(R)$ and under $\Spec(R_{\sf red})$
    \[Y_R=\Spec(R)\underset{Y_{\DR}}{\times} Y \ra \Spec(R)\underset{X_{\DR}}{\times}X~ = X_R\]
    by taking relative cotangent complexes and then taking symmetric algebras. 
    Consider the map of cotangent complexes 
    \[\LL_{R_{\sf red}\vert X_R} \ra \LL_{R_{\sf red}\vert Y_R}\]
    in $\Mod_{R_{\sf red}}^{\LL_{R_{\sf red}\vert R}}$. It fits into a cofiber sequence 
    \[h^*\LL_{Y_R\vert X_R} \longrightarrow \LL_{R_{\sf red}\vert X_R} \longrightarrow \LL_{R_{\sf red}\vert Y_R}~\]
    where $h$ denotes the canonical map $h\colon \Spec(R_{\sf red}) \ra Y_R$. However, note that $\LL_{Y_R\vert X_R}$ is the pullback of $\LL_{Y\vert X}$. As $f\colon Y \ra X$ is formally \'etale we have $\LL_{Y\vert X}\simeq 0$. Therefore we have a cofiber sequence 
    \[0\simeq h^*\LL_{Y_R\vert X_R} \ra \LL_{R_{\sf red}\vert X_R} \ra \LL_{R_{\sf red}\vert Y_R}~\]
    and we have the result. 
\end{proof}
Therefore given a formally \'etale map $Y \ra X$, a quantization of principal parts on $X$ restricts to a deformation of principal parts on $Y$. 
\begin{lemma}
    Assume $f\colon Y \ra X$ is a formally \'etale map between derived Artin stacks, locally of finite presentation over $\Bbbk$. Given an $\EE_{n+1}$-deformation of principal parts $\cP_X$, i.e., an object
   \[\cP^{\pi}_X[\![\hbar]\!] \in \Alg_{\EE_{n+1}}(\Mod_{\cD_{X_{\DR}}[\![\hbar]\!]}^{\epsilon, \gr})\]
   together with an equivalence 
   \[\cP_X^{\pi}[\![\hbar]\!] \otimes_{\cD_{X_{\DR}}[\![\hbar]\!]}{\cD_{X_{\DR}}}|_{\hbar=0} \simeq \cP_{X}~,\]
   we obtain a $\EE_{n+1}$-deformation of $\cP_Y$, i.e., an object 
   \[\cP^{\pi}_Y[\![\hbar]\!] \in \Alg_{\EE_{n+1}}(\Mod_{\cD_{Y_{\DR}}[\![\hbar]\!]}^{\epsilon, \gr})\]
   together with an equivalence 
   \[\cP_Y^{\pi}[\![\hbar]\!] \otimes_{\cD_{Y_{\DR}}[\![\hbar]\!]}{\cD_{Y_{\DR}}}|_{\hbar=0} \simeq \cP_{Y}~.\]
\end{lemma}
In particular, the constructions of Lemma \ref{quantizingordinaryprincipalparts} and \ref{quantizingglobalfunctions} also restrict appropriately along formally \'etale maps. 
Let $X_{\sf f\acute{e}t}$ denote the full subcategory of $(\dSt_{/X})^{\op}$ spanned by formally \'etale maps $Y \ra X$, and covers are given by jointly surjective formally \'etale maps. 
Therefore we have the following lemma. 

\begin{lemma} \label{quantizingthestructuresheaf}
 Let $X$ be a derived Artin stack, locally of finite presentation over $\Bbbk$. Consider the family $q\colon X \ra X_{\DR}$ Given an $\EE_{n+1}$-deformation of principal parts $\cP_X$, i.e., an object
   \[\cP^{\pi}_X[\![\hbar]\!] \in \Alg_{\EE_{n+1}}(\Mod_{\cD_{X_{\DR}}[\![\hbar]\!]}^{\epsilon, \gr})\]
   together with an equivalence 
   \[\cP^{\pi}_X[\![\hbar]\!] \otimes_{\cD_{X_{\DR}}[\![\hbar]\!]}{\cD_{X_{\DR}}}|_{\hbar=0} \simeq \cP_{X}~,\]
   we obtain a sheaf
   \[\cO_X^{\pi}[\![\hbar]\!] \in \Alg_{\EE_{n+1}}\bigl({\sf Shv}(X_{\sf f\acute{e}t})\bigr)\]
   deforming the structure sheaf $\cO_X$ on $X_{\sf f\acute{e}t}$. 
\end{lemma}

\begin{proof}
As deformations restrict along formally \'etale maps, by Lemma \ref{quantizingglobalfunctions}, we immediately have a quantization of the presheaf 
\[
\begin{array}{ccccc}
\cO^{\sf f\acute{e}t}_X: & X_{\sf f\acute{e}t}& \ra&  \Mod_{\Bbbk}\\
&(Y \xra{f} X) &\mapsto &\Gamma(Y, \cO_Y)
\end{array}
\]
given by 
\[
\begin{array}{ccccc}
\cO_X^{\pi, {\sf f\acute{e}t}}[\![\hbar]\!]:& X_{\sf f\acute{e}t}& \ra &\Alg_{\EE_{n+1}}(\Mod_{\Bbbk[\![\hbar]\!]})\\
&(Y \xra{f} X) &\mapsto &\Gamma(Y_{\DR}, \vert\cP^{\pi}_Y[\![\hbar]\!]\vert).
\end{array}
\]
Descent is agnostic to the algebra structures. 
 he underlying module of $\cO_X^{\pi, {\sf f\acute{e}t}}[\![\hbar]\!]$ is given by $\cO_X^{\sf f\acute{e}t}[\![\hbar]\!]$, which satisfies descent. Therefore we obtain a sheaf $\cO_X^{\pi, {\sf f\acute{e}t}}[\![\hbar]\!]$ deforming the sheaf $\cO_X$. 
\end{proof}

\begin{remark}
    In the statements of Lemmas  \ref{quantizingordinaryprincipalparts}, \ref{quantizingglobalfunctions}, and \ref{quantizingthestructuresheaf}, we are not claiming to quantize a given Poisson algebra structure. Nevertheless, nontriviality of the deformations can be detected by defining semi-classical structures on cohomology in a similar manner to the previous subsection.  
~\hfill$\Diamond$\end{remark}

We now discuss deformations of the category of perfect complexes $\Perf(X) \simeq \Perf^{\epsilon, \gr}_{\cP_X}$. Given an $\EE_{n+1}$-quantization $\cP_X^{\pi}[\![\hbar]\!]$ of $\cP_X$, we have the $\EE_n$-category $\Mod_{\cP^{\pi}_{X}[\![\hbar]\!]}\bigl(\Mod^{\epsilon, \gr}_{\cD_{X_{\DR}}[\![\hbar]\!]}\bigr)$. 

\begin{definition}[\cite{CPTVV} Definition 3.5.5]\label{quantumperfect}
    The category of \em{quantum perfect modules} $\Perf(X, \pi):= \Perf_{\cP_X^{\pi}[\![\hbar]\!]}^{\epsilon, \gr}$ is the full subcategory of $\Mod_{\cP_X^{\pi}[\![\hbar]\!]}\bigl(\Mod^{\epsilon, \gr}_{\cD_{X_{\DR}}[\![\hbar]\!]}\bigr)$ spanned by objects $E\in \Mod_{\cP_X^{\pi}[\![\hbar]\!]}\bigl(\Mod^{\epsilon, \gr}_{\cD_{X_{\DR}}[\![\hbar]\!]}\bigr)$ satisfying the following conditions:
    \begin{itemize}
        \item for every $f\colon\Spec(R) \ra X_{\DR}$, we have 
        \[{\sf fgt}_{\epsilon}(E(f)) \simeq {\sf fgt}_{\epsilon}(\cP_X^{\pi}[\![\hbar]\!](f))\underset{A_{\sf red}}{\otimes} L_f\]
        where $L_f \in \Perf_{A_{\sf red}}$
        \item for every $\Spec(S) \xra{g} \Spec(R) \xra{f} X_{\DR}$, the canonical map 
        \[E(f) \underset{\cP_X^{\pi}[\![\hbar]\!](f)}{\otimes} \cP_X^{\pi}[\![\hbar]\!](f\circ g) \xra{\simeq} E(f\circ g)\]
        is an equivalence of graded mixed $\cP_X^{\pi}[\![\hbar]\!](g\circ f)$-modules. 
    \end{itemize}
\end{definition}

The category of quantum perfect modules should be thought of as a one parameter family of $\EE_n$-categories whose zero fiber is the category 
\[\Perf_{\cP_X^{\pi}}^{\epsilon, \gr} \simeq \Perf(X)~.\]

\begin{remark}
    Note that the lemma above is not claiming that $\Perf(X)$ is a deformation quantization of the induced $\PP_n$-category structure on $\Perf(X)$. See \cite{CPTVV} Conjecture 3.5.7. 
~\hfill$\Diamond$\end{remark}

\section{Principal Parts on Mapping Stacks with Betti Source}\label{sec.principal.parts}

Let $X$ be a derived Artin stack locally of finite presentation over $\Bbbk$. 
Let $M$ be a smooth $d$-manifold with the homotopy type of a finite CW complex,
and let $\Pi_\infty M$ denote its underlying $\infty$-groupoid,
which will often be viewed as a category below to describe local systems on $M$.
Let $X^M$ denote the mapping stack $\Map_{\dSt}(\Pi_\infty M, X)$ whose functor of points~is 
\[
\begin{array}{ccccc}
X^M\colon & \CAlg^{\geq 0}_{\Bbbk} & \longrightarrow & \spaces\\
& R & \longmapsto & \Map_{\spaces}(\Pi_\infty M, X(R)).
\end{array}
\]
The finiteness condition on the homotopy type of $M$ ensures that the mapping stack $X^M$ is a derived Artin stack locally of finite presentation over $\Bbbk$~\cite[\S2]{PTVV}. 
(Note that $X^M$ only cares about $M$ as a homotopy type and not as a smooth manifold.)

In this section, we will unwind the formal localization machinery to the family $X^M \rightarrow X^M_{\DR}$. 

\begin{remark}
The functor $(-)_{\DR}: \dSt \to \dSt$ is a right adjoint (see \cite{CPTVV} section 2.1). 
Therefore, we have $(X^M)_{\DR} \simeq (X_{\DR})^M$ and so our notation $X^M_{\DR}$ is unambiguous. 
~\hfill$\Diamond$\end{remark}

We will describe the presheaf of principal parts $\cP_{X^M}$ as the homology of $M$ with coefficients in a local system~$\cA_X \in \Fun\Bigl(\Pi_\infty M, \CAlg\bigl(\Mod_{\cD_{X^M_{\DR}}}^{\epsilon, \gr}\bigr)\Bigr)$. Informally, the local system $\cA_X$ is given by  
\[
\begin{array}{cccc}
    \cA_X\colon & \Pi_\infty M & \longrightarrow & \CAlg\bigl(\Mod_{\cD_{X^M_{\DR}}}^{\epsilon, \gr}\bigr)  \\
     & x & \mapsto & {\sf ev}_x^{\ast}\cP_{X}
\end{array}
\]
where $\ev_x\colon X^M_{\DR} \ra X_{\DR}$ is evaluation at a point~$x\in M$. 
The local system $\cA_X$ will be obtained from a local system
\[\cA_X^{\sf geom}\colon \Pi_\infty M \ra \Fun\bigl({\sf Aff}_{/X^M_{\DR}}, \dSt\bigr)\]
whose fiber at $x\in M$ is the functor 
\[
\begin{array}{ccccc}
     \cA_X^{\sf geom}(x)\colon & {\sf Aff}_{/X^M_{\DR}} & \longrightarrow & \dSt \\
     &(\Spec(R)\xra{f} X^M_{\DR}) & \mapsto & \Spec(R)\times_{X_{\DR}} X 
\end{array}
\]
where the map $\Spec(R) \ra X_{\DR}$ is the composite 
\[\Spec(R) \xra{f} X^M_{\DR} \xra{{\sf ev}_x} X_{\DR}~.\]
Observe that the composition 
\[\Pi_\infty M \xra{\cA_X^{\sf geom}} \Fun({\sf Aff}_{/X^M_{\DR}}, \dSt) \xra{\DD_{\ast}} \Fun\bigl({(\sf dAff}_{/X^M_{\DR}})^{\op}, \CAlg(\Mod_{\Bbbk}^{\epsilon, \gr})\bigr)\]
is the local system $\cA_X$. 
We will endow the local system $\cA_X^{\sf geom}$ with some additional structure. 
We have two constant local systems 
\[\cP_{X^M}^{\sf geom}\colon \Pi_\infty M \ra \ast \ra\Fun({\sf Aff}_{/X^M_{\DR}}, \dSt)\]
\[\cD_{X^M_{\DR}}^{\sf geom}\colon \Pi_\infty M \ra \ast \ra \Fun({\sf Aff}_{/X^M_{\DR}}, \dSt)\]
encoding geometric versions of the presheaf of principal parts and the crystalline structure sheaf associated to the family $X^M \ra X^M_{\DR}$. We will construct a commutative triangle. 
\[
\begin{tikzcd}
  {\cP^{\sf geom}_{X^M}} \arrow[rr] \arrow[dr]
    && {\cA_X^{\sf geom}} \arrow[dl] \\
  & {\cD^{\sf geom}_{X^M_{\DR}}}
\end{tikzcd}
\]

Post-composition with $\DD \colon \dSt^{\op} \ra \CAlg(\Mod^{\epsilon, \gr}_{\Bbbk})$ will produce a commutative triangle in the category $\Fun\Bigl(\Pi_\infty M, \Fun\bigl(({\sf Aff}_{/X^M_{\DR}})^{\op}, \CAlg(\Mod_{\Bbbk}^{\epsilon, \gr})\bigr)\Bigr)$: 
\[
\begin{tikzcd}
    & {\cD_{X^M_{\DR}}} \arrow[dl] \arrow[dr] \\
  {\cA_X} \arrow[rr] 
    && {\cP_{X^M}}  
\end{tikzcd}
\]
where $\cP_{X^M}$ and $\cD_{X^M_{\DR}}$ are viewed as constant local systems. By adjunction,  we will obtain a canonical map 
\[\underset{\Pi_\infty M}\colim \ \cA_X \ra \cP_{X^M}\]
in the category $\Mod^{\epsilon, \gr}_{\cD_{X^M_{\DR}}}$. This map is an equivalence, describing $\cP_{X^M}$ as the homology of $\cA_X$. 
\begin{itemize}
    \item In Subsection~\ref{subsectiongeometriclocalsystem}, we will construct the local system $\cA_X^{\sf geom}$. It is rather technical and can be skipped if the reader is satisfied with the pointwise definition provided above. 
    
    \item In Subsection~\ref{subsectioncotangentcomplexes}, we will collect some useful computations involving cotangent complexes on mapping stacks with Betti source.

    \item In Subsection~\ref{subsectionalgebraiclocalsystem}, we will construct the local system $\cA_X$ and describe $\cP_{X^M}$ as principal parts as the homology of $\cA_X$. 
\end{itemize}

\subsection{The Local System $\cA_X^{\sf geom}$} \label{subsectiongeometriclocalsystem}

\begin{definition}
Let $\cC$ denote the category that is the (nerve of the) $3$-object category 
    \[
    \begin{tikzcd}
        & \ast \arrow[d]\\
        \ast \arrow[r] & \ast
    \end{tikzcd}
    \]
\end{definition}

Consider the functor category $\Fun(\cC, \dSt)$. 
An object of this category is a diagram 
\[X \ra Y \la Z~.\]
We have a functor 
\[
\begin{array}{ccc}
     {\sf Aff}_{/X^M_{\DR}} & \longrightarrow & \Fun(\cC, \dSt)  \\
      (\Spec(R) \ra X^M_{\DR}) & \mapsto & (\Spec(R) \ra X^M_{\DR} \la X^M)
\end{array}
\]
appearing naturally in our setting.

\begin{definition}
The functor 
\[
{\sf fib}_{X^M}\colon 
\Pi_\infty M \times {\sf Aff}_{/X^M_{\DR}} 
\ra 
\Fun(\cC, \dSt)
\]
is the composite 
\[
\Pi_\infty M \times{\sf Aff}_{/X^M_{\DR}} 
\xra{p_2} 
{\sf Aff}_{/X^M_{\DR}} 
\ra \Fun(\cC, \dSt)
\]
where $p_2$ is projection onto the second factor. 
\end{definition}

Similarly, consider the functor 
\[
\begin{array}{ccc}
{\sf Aff}_{/X_{\DR}} & \longrightarrow & \Fun(\cC, \dSt) \\
(\Spec(R) \ra X_{\DR})&  \mapsto &  (\Spec(R) \ra X_{\DR} \la X)~.
\end{array}
\]
Recall the canonical evaluation map of derived stacks
\[
{\sf ev}\colon 
\Pi_\infty M \times X^M \ra X
\]
and the induced functor 
\[
M\times{\sf Aff}_{/X^M_{\DR}} \xra{\sf ev} {\sf Aff}_{/X_{\DR}}~.
\]

\begin{definition}
The functor 
\[
{\sf fib}_{X^M}^{\sf loc}\colon 
\Pi_\infty M \times {\sf Aff}_{/X^M_{\DR}} 
\ra 
\Fun(\cC, \dSt)\]
is the composite 
\[
\Pi_\infty M\times {\sf Aff}_{/X^M_{\DR}} 
\xra{\sf ev} 
{\sf Aff}_{X_{\DR}} 
\ra 
\Fun(\cC, \dSt)~.
\]
\end{definition}

\begin{remark}
    The local system $M \ra \Fun({\sf Aff}_{/X^M_{\DR}}, \Fun(\cC, \dSt))$ adjoint to ${\sf fib}_{X^M}^{\sf loc} \in \Fun(\Pi_\infty M \times {\sf Aff}_{/X^M_{\DR}}, \Fun(\cC, \dSt))$ has fibers
    \[x \in M \mapsto ({\sf fib}_{X^M}^{\sf loc})_x\colon {\sf Aff}_{/X^M_{\DR}} \ra \Fun(\cC, \dSt)\]
    where 
    \[
    \begin{array}{cccc}
         ({\sf fib}_{X^M}^{\sf loc})_x\colon & {\sf Aff}_{/X^M_{\DR}} & \longrightarrow & \Fun(\cC, \dSt)\\
         &(\Spec(R)\xra{f}X^M_{\DR}) &\mapsto & (\Spec(R) \xra{{\sf ev}_x\circ f}X_{\DR} \la X)~.
    \end{array}
    \]
    ~\hfill$\Diamond$
\end{remark}

We will construct a natural transformation of functors 
\[\eta\colon {\sf fib}_{X^M} \ra {\sf fib}_{X^M}^{\sf loc}~.\]
The natural transformation $\eta$ gives us the following: for each $x\in M$ and $f\colon \Spec(R) \ra X^M_{\DR}$, we obtain a commutative diagram in $\dSt$
\[
\begin{tikzcd}
   &\Spec(R) \ar[r, "f"] \ar[d, "{\sf id}"] &{X^M_{\DR}} \ar[d, "{\sf ev}_x"]  &{X^M} \ar[l] \ar[d] \\
   &\Spec(R) \ar[r, "{\sf ev}_x \circ f"] &{X_{\DR}}  & X \ar[l]
\end{tikzcd}
\]
Composing with 
\[\lim\colon \Fun(\cC, \dSt) \ra (\dSt)~,\]
for each $x\in M$ and $f\colon \Spec(R) \ra X^M_{\DR}$, we obtain a commutative diagram 

\[\begin{tikzcd}[row sep={70,between origins}, column sep={70,between origins}]
      & X^M_R \ar{rr}\ar{dd}\ar[dashed, dl] & & X^M \ar{dd}\ar[dl, "{\sf ev}_x"] \\
    \Spec(R)\underset{X_{\DR}}\times X \ar[crossing over]{rr} \ar{dd} & & X \\
      & \Spec(R)  \ar[rr, "f", {xshift=10pt}] \ar[dl, "{\sf id}"] & &  X^M_{\DR} \ar[dl, "{\sf ev}_x"] \\
    \Spec(R) \ar[rr, "{\sf ev}_x \circ f"] && X_{\DR} \ar[from=uu,crossing over]
\end{tikzcd}\]
We need to introduce some notation. 

\begin{definition}
The functor 
\[
c_{(-)}\colon \dSt \ra \Fun(\Pi_\infty M, \dSt)
\]
takes a derived stack $Y$ to the constant local system 
\[
c_Y: \Pi_\infty M \ra \ast \xra{Y} \dSt~.
\]
\end{definition}

We have identifications of mapping spaces 
\[\Map_{\Fun(\Pi_\infty M, \dSt)}(c_{X^M}, c_{X}) \simeq \Map_{\dSt}(\colim c_{X^M}, X) \simeq \Map_{\dSt}(\Pi_\infty M\times X^M, X)~.\]
Therefore, the evaluation map of derived stacks ${\sf ev}: M\times X^M \ra X$ induces a map of constant local systems 
\[c_{X^M} \xra{{\sf ev}} c_X~.\]

\begin{definition}
The functor 
\[
{\sf ev}: \Pi_\infty M \times [1] \ra \dSt
\]
is the adjoint of 
\[
[1] \xra{\sf ev} \Fun(\Pi_\infty M, \dSt)
\]
picking out the morphism of local systems
\[c_{X^M} \xra{{\sf ev}} c_X\]
described above. 
\end{definition}

Recall that we have a natural transformation of endofunctors on $\dSt$:
\[{\sf Crys}\colon{\sf id} \ra (-)_{\DR} \in \Fun(\dSt \times [1], \dSt)~.\]

\begin{definition}
    The functor 
    \[{\sf Crys}^{\sf loc}\colon  [1] \times [1] \ra  \Fun(\Pi_\infty M \times {\sf Aff}_{/X^M_{\DR}}, \dSt)\]
    is the adjoint to the composite 
    \[ \Pi_\infty M\times {\sf Aff}_{/X^M_{\DR}} \times [1] \times [1] \xra{\sf pr} \Pi_\infty M \times [1] \times [1] \xra{{\sf ev}\times{\sf id}} \dSt \times [1] \xra{\sf Crys} \dSt~.\]
\end{definition}

\begin{remark}
    Given $x \in M$ and $f \colon \Spec(R) \ra X^M_{\DR}$, the composite above is the data of a commutative diagram 
    \[
    \begin{tikzcd}
        X^M \arrow[r] \arrow[d, "{\sf ev}_x"] & X^M_{\DR} \arrow[d, "{\sf ev}_x"] \\
        X \arrow [r] & X_{\DR}
    \end{tikzcd}
    \]
Note that ${\sf Crys}^{\sf loc}$ does not depend on the map~$f \colon \Spec(R) \ra X^M_{\DR}$. Informally, ${\sf Crys}^{\sf loc}$ is half of the data of the desired natural transformation $\eta$. 
~\hfill$\Diamond$
\end{remark}

We now construct the other half of the natural transformation $\eta$. 
\begin{definition}
    Given a derived stack $Y \in \dSt$, $\cD_{Y}^{\sf geom}$ is the functor 
    \[
    \begin{array}{cccc}
        \cD_{Y}^{\sf geom}\colon & {\sf Aff}_{/Y} &\longrightarrow &\dSt \\
         & (\Spec(R)\ra Y) & \mapsto & \Spec(R)~.
    \end{array}
    \]
\end{definition}
Consider the projections 
\[\Pi_\infty M \xla{p_1} \Pi_\infty M \times {\sf Aff}_{/X^M_{\DR}} \xra{p_2} {\sf Aff}_{/X^M_{\DR}}~.\]
Observe that there is a canonical map in~$\Fun(\Pi_\infty M \times {\sf Aff}_{/X^M_{\DR}}, \dSt)$:
\[p_2^{\ast} \cD_{X^M_{\DR}}^{\sf geom} \ra p_1^*c_{X^M_{\DR}}~.\]
Consider the composable maps 
\[p_2^{\ast} \cD_{X^M_{\DR}}^{\sf geom} \ra p_1^*c_{X^M_{\DR}} \xra{p_1^{\ast}{\sf ev}} p_1^{\ast}c_{X_{\DR}}~.\]
\begin{definition}
The functor 
 \[{\sf Test}^{\sf loc}\colon [1] \times [1] \ra \Fun(\Pi_\infty M \times {\sf Aff}_{/X^M_{\DR}}, \dSt)\] 
 is the data of a composition
 \[
 \begin{tikzcd}
     p_2^{\ast} \cD_{X^M_{\DR}}^{\sf geom} \arrow[r] \arrow[d, "{\sf id}"] &  p_1^*c_{X^M_{\DR}} \arrow[d, "p_1^{\ast}{\sf ev}"] \\
     p_2^{\ast} \cD_{X^M_{\DR}}^{\sf geom} \arrow[r] & p_1^{\ast}c_{X_{\DR}}
 \end{tikzcd}.
 \]
\end{definition}

\begin{remark}
 At $x\in M$ and $f \colon\Spec(R) \ra X^M_{\DR}$, ${\sf Test}^{\sf loc}$ is the data of a commutative square 
 \[\begin{tikzcd}
	{\Spec(R)} & {X^M_{\DR}} \\
	{\Spec(R)} & {X_{\DR}}
	\arrow["f", from=1-1, to=1-2]
	\arrow["{{\sf id}}"', from=1-1, to=2-1]
	\arrow["{{\ev_x}}", from=1-2, to=2-2]
	\arrow["{{\sf ev}_x \circ f}"', from=2-1, to=2-2]
\end{tikzcd}\]
This is the other half of the desired natural transformation $\eta$. 
~\hfill$\Diamond$
\end{remark}

\begin{lemma} \label{compgeomspan}
    The functors 
    \[{\sf Crys}^{\sf loc}\colon  [1] \times [1] \ra \Fun(\Pi_\infty M \times {\sf Aff}_{/X^M_{\DR}}, \dSt)\]
    and 
     \[{\sf Test}^{\sf loc}\colon [1] \times [1] \ra \Fun(\Pi_\infty M \times {\sf Aff}_{X^M_{\DR}}, \dSt)\] 
     glue along their common edge to give us a natural transformation 
    \[\eta\colon {\sf fib}_{X^M} \ra {\sf fib}_{X^M}^{\sf loc}~.\]
\end{lemma}
\begin{proof}
Write 
\[[1]:= (0 \ra 1)~.\]
Observe from the constructions of ${\sf Crys}^{\sf loc}$ and ${\sf Test}^{\sf loc}$ that we have an equivalence 
\[{\sf Crys}^{\sf loc}|_{[1]\times\{1\}} \simeq {\sf Test}^{\sf loc}|_{[1]\times \{1\}}.\]
Therefore we have an induced functor 
\[\eta\colon  \bigl([1]\times [1]\bigr)\bigsqcup_{[1]\times\{1\}}\bigl([1]\times [1]\bigr) \ra \Fun(\Pi_\infty M \times {\sf Aff}_{/X^M_{\DR}}, \dSt)~.\]
Note that we have $\bigl([1]\times [1]\bigr)\bigsqcup_{[1]\times\{1\}}\bigl([1]\times [1]\bigr) \simeq [1]\times \cC$, therefore we have a functor 
\[\eta\colon [1]\times \cC \ra \Fun(\Pi_\infty M \times {\sf Aff}_{/X^M_{\DR}}, \dSt)~.\]
Taking adjoints, we obtain a functor 
\[\eta\colon [1] \ra \Fun\bigl(\Pi_\infty M \times {\sf Aff}_{/X^M_{\DR}}, \Fun(\cC, \dSt)\bigr)\]
such that 
\[\eta|_{\{0\}} \simeq {\sf fib}_{X^M}\]
and 
\[\eta|_{\{1\}} \simeq {\sf fib}^{\sf loc}_{X^M}~,\]
producing our desired natural transformation. 
\end{proof}

\begin{definition}
    The functor 
    \[\cP_{X^M}^{\sf geom}: \Pi_\infty M \ra \Fun({\sf Aff}_{/X^M_{\DR}}, \dSt)\]
    is the adjoint of the composite 
    \[\Pi_\infty M \times {\sf Aff}_{X^M_{\DR}} \xra{{\sf fib}_{X^M}} \Fun(\cC, \dSt) \xra{\lim} \dSt~.\]
\end{definition}

\begin{definition}
    The functor 
    \[\cA_X^{\sf geom}\colon \Pi_\infty M \ra \Fun({\sf Aff}_{/X^M_{\DR}}, \dSt)\]
    is the adjoint of the composite 
    \[M \times {\sf Aff}_{X^M_{\DR}} \xra{{\sf fib}^{\sf loc}_{X^M}} \Fun(\cC, \dSt) \xra{\lim} \dSt~.\]    
\end{definition}
Observe that by construction, we have natural maps 
\[\cP^{\sf geom}_{X^M} \ra \cD_{X^M_{\DR}}^{\sf geom}\]
and 
\[\cA_X^{\sf geom} \ra \cD_{X^M_{\DR}}^{\sf geom}~,\]
where we view $\cD_{X^M_{\DR}}^{\sf geom}$ as a constant local system. 
\begin{definition}
    Given a derived stack $Y \in \dSt$, $\cD_Y^{\sf geom, red}$ is the functor 
    \[
    \begin{array}{cccc}
        \cD_Y^{\sf geom, red}\colon & {\sf Aff}_{/Y} & \longrightarrow & \dSt \\
         & \Spec(R) & \mapsto &\Spec(R_{\sf red})
    \end{array}
    \]
\end{definition}
Again, viewing $\cD_{X^M_{\DR}}^{\sf geom, red}$ as a constant local system on $M$, observe that we have natural maps 
\[\cD_{X^M_{\DR}}^{\sf geom, red} \ra \cP_{X^M}^{\sf geom}\]
and 
\[\cD_{X^M_{\DR}}^{\sf geom, red} \ra \cA_X^{\sf geom}~.\]
View $\cP_{X^M}^{\sf geom}$ and $\cA_X^{\sf geom}$ as objects of the category~$^{\cD_{X^M_{\DR}}^{\sf geom, red}/}\Fun(\Pi_\infty M, \Fun({\sf Aff}_{/X^M_{\DR}}, \dSt))_{/\cD_{X^M_{\DR}}^{\sf geom}}$.
The following is immediate from Lemma~\ref{compgeomspan}, and the constructions of $\cP_{X^M}^{\sf geom}$ and~$\cA_X^{\sf geom}$.

\begin{cor} \label{compgeomlocsys}
There is a map 
\[\eta\colon \cP_{X^M}^{\sf geom} \ra \cA_X^{\sf geom}\]
of objects of the category~$^{\cD_{X^M_{\DR}}^{\sf geom, red}/}\Fun(\Pi_\infty M, \Fun({\sf Aff}_{/X^M_{\DR}}, \dSt))_{/\cD_{X^M_{\DR}}^{\sf geom}}$.
\end{cor}

\subsection{Cotangent Complexes} \label{subsectioncotangentcomplexes}
In this subsection, we will collect some useful computations involving cotangent complexes. We will use the description of the cotangent complex as an object representing the functor of derivations. 
We begin with some recollections.

\begin{definition}[Definition 1.3.13~\cite{DAG}]\label{der}
Let $Y \in \dSt$ be a derived stack, and let $R \in \CAlg^{\geq 0}$ be a connective commutative algebra. 
Given a map $f\colon\Spec(R) \ra Y$ and a connective $R$-module~$E$, 
let $R[E]$ denote the trivial square-zero extension of $R$ by~$E$. 
The space of \emph{derivations at $f$ with values in $E$} is the pullback 
\[
\begin{tikzcd}
\Der_Y(f, E) & {Y(R[E])} \\
\pt & {Y(R)}
\arrow[from=1-1, to=1-2]
\arrow[from=1-1, to=2-1]
\arrow[from=1-2, to=2-2]
\arrow["f"', from=2-1, to=2-2]
\end{tikzcd}
\]
in $\spaces$.
\end{definition}

\begin{remark}
    This construction is functorial in $f$ and $E$. See \cite{DAG} Definition 1.3.13. ~\hfill$\Diamond$
\end{remark}

\begin{remark}
For affine $Y$, the definition above agrees with Definition 1.3.3 in \cite{CPTVV}. Indeed, assume $Y = \Spec(R)$. For $E \in \Mod_R$, \cite{CPTVV} define 
\[\Der(R, E):= \Map_{{\CAlg_{\Bbbk}}_{/R}}(R, R[E])\]
where $R[E]$ is the trivial square-zero extension of $R$ by $E$ (note that this is naturally an $R$-algebra with an augmentation $R[E]\rightarrow R$). The mapping space above can be written as the fiber product
\[\begin{tikzcd}
{\Map_{\CAlg_{\Bbbk}/R}(R, R[E])} & {\Map_{\CAlg_{\Bbbk}}(R, R[E])} \\
\pt & {\Map_{\CAlg_{\Bbbk}}(R, R)}
\arrow[from=1-1, to=1-2]
\arrow[from=1-1, to=2-1]
\arrow[from=1-2, to=2-2]
\arrow["{{\sf id}_R}", from=2-1, to=2-2]
\end{tikzcd}\]
giving us 
\[\Der(R, E) \simeq \Der_{\Spec(R)}(\id, E)\]
as wanted.~\hfill$\Diamond$
\end{remark}

For any mapping stack $Y^S$ with Betti source~$S \in \Spaces$, 
there is a local (on $S$) description of these spaces of derivations. 
Take a map $f\colon\Spec(R) \ra Y^S$ and a connective $R$-module~$E$. 
For each point $x$ of $S$, we have a map 
\[
f_x := {\sf ev}_x\circ f\colon \Spec(R) \to Y
\]
and its associated space of derivations~$\Der_Y(f_x, E)$. 
As the construction of the space of derivations is functorial, these spaces of derivations organize into a local system of spaces over $S$. More precisely, the local system is given by the composite 
\[S \xra{f} Y(R) \xra{\Der(-, E)} \Spaces~.\]
$\Der_{Y^S}(f, E)$ is recovered as the global sections of this local system. 

\begin{lemma}
\label{derivations}
Given a map $f\colon \Spec(R) \ra Y^S$ and a connective $R$-module $E$, there is an equivalence 
\[
\Der_{Y^S}(f, E) \simeq \lim_{x \in S}\Der_{Y}(f_x, E)
\]
in $\spaces$.
\end{lemma}

\begin{proof}
We have
\[
\begin{split}
\Der_{Y^S}(f, E) &\simeq {\sf fib}_f\bigl(Y^S(R[E])\to Y^S(R)\bigr) \\
&\simeq {\sf fib}_f\bigl(\Map(S, Y(R[E])) \to \Map(S, Y(R))\bigr) \\
&\simeq \lim_{x\in S}{\sf fib}_{f_x}\bigl(Y(R[E]) \to Y(R)\bigr)\\
&\simeq \lim_{x \in S}\Der_{Y}(f_x, E)
\end{split}
\]
where $\Map$ here denotes the mapping spaces.
\end{proof}

This description of derivations leads to a convenient description of the absolute cotangent complex of mapping stacks with Betti source. 
We now review pertinent facts about the cotangent complex,
which can be described as an object representing the functor of derivations. 

\begin{definition}[Definition 1.3.13~\cite{DAG}]
\label{defcot}
Let $Y \in \dSt$ be a derived stack. 
An object $\LL_Y \in \qc(Y)$ is a \emph{cotangent complex for $Y$} if there is a natural equivalence
\[\Der_Y(f, E) \simeq \Map_{\Mod_R}(f^*\LL_Y, E)\]
for every map $f\colon\Spec(R) \ra Y$ and connective $R$-module~$E$.
\end{definition} 
Assume $Y$ admits a cotangent complex. 
The local to global description of derivations on the Betti mapping stack $Y^S$ will allow us to describe the cotangent complex in a similar manner. 
Take a map $f\colon \Spec(R) \ra Y^S$. 
For each $x \in S$, we obtain an $R$-module $f_x^*\LL_Y$. 
As $x$ varies over $S$, these $R$-modules organize into a local system over $S$. More precisely, taking pullbacks of $\LL_Y$, we have a functor 
\[Y(R) \xra{(-)^*\LL_Y} \Mod_R~,\]
and therefore we obtain a local system of $R$-modules 
\[S \xra{f} Y(R) \xra{(-)^*\LL_Y} \Mod_R~.\]
The pullback $f^*\LL_{Y^M}$ has the following description. 

\begin{lemma}
\label{cot}
Given a map $f\colon \Spec(R) \ra Y^S$, consider the local system of $R$-modules on $S$ whose fiber over $x\in S$ is given by $f_x^*\LL_{Y}$. 
There is an equivalence 
\[
f^*\LL_{Y^S} \simeq \underset{x \in S}\colim \ f_x^*\LL_{Y}
\]
of $R$-modules.
\end{lemma}

\begin{proof}
Given a connective $R$-module $E$, we have
\[
\begin{split}
\Map_{\Mod_R}(\colim_{x \in S}f_x^*\LL_{Y}, E) &\simeq \lim_{x \in S} \Map_{\Mod_R}(f_x^*\LL_{Y}, E)\\
&\simeq \lim_{x\in S}\Der_{Y}(f_x, E)\\
&\simeq \Der_{Y^S}(f, E)
\end{split}
\]
where the last equivalence follows from Lemma \ref{derivations}.
\end{proof}

Given a map of derived stacks $f\colon Y \ra Z$, the relative cotangent complex $\LL_{Y | Z}$ has a description in terms of the absolute cotangent complexes of $Y$ and $Z$.

\begin{definition}
\label{defrelcot}
Let the derived stacks $Y$ and $Z$ both admit cotangent complexes.
Given a map of derived stacks $f \colon Y \ra Z$, the \emph{relative cotangent complex of the map $f$} is the cofiber 
\[
\LL_f := {\sf cofib}(f^*\LL_Z \rightarrow \LL_Y)~.
\]
We also use $\LL_{Y| Z}$ when the map is clear from context.
\end{definition}

We now return to studying the family $X^M \to X^M_{\DR}$. 
Recall the map 
\[\eta\colon \cP^{\sf geom}_{X^M} \ra \cA_X^{\sf geom}\]
of objects of the category~$^{\cD_{X^M_{\DR}}^{\sf geom, red}/}\Fun(\Pi_\infty M, \Fun({\sf Aff}_{/X^M_{\DR}}, \dSt))_{/\cD_{X^M_{\DR}}^{\sf geom}}$. Given a map $f\colon \Spec(R) \ra X^M_{\DR}$, we obtain a map 
\[\eta_f\colon  \cP^{\sf geom}_{X^M}(f) \ra \cA_X^{\sf geom}(f)\]
in the category $\Fun\bigl(\Pi_\infty M, ^{\Spec(R_{\sf red})/}\dSt_{/\Spec(R)}\bigr)$. Note that $\cP_{X^M}^{\sf geom}(f)$ is the constant local system on $M$ with values in  $^{\Spec(R_{\sf red})/}\dSt_{/\Spec(R)}$, whose fiber is given by 
\[\Spec(R_{\sf red}) \ra X^M_R = \Spec(R) \times_{X^M_{\DR}} X^M \ra \Spec(R)~.\]
The local system $\cA_X^{\sf geom}(f)$ has fiber over $x\in M$ given by 
\[\Spec(R_{\sf red}) \ra X_{R,x} = \Spec(R)\times_{X_{\DR}} X \ra \Spec(R)\]
where the map $\Spec(R) \ra X_{\DR}$ is the composite
\[\Spec(R) \xra{f} X^M_{\DR} \xra{{\sf ev}_x} X_{\DR}~.\]
The map of local systems $\eta_f$ at a point $x\in M$ is the data of a commutative diagram 
\[\begin{tikzcd}
	& {\Spec(R_{\sf red})} & \\
	{X^M_R} && {X_{R,x}} \\
	& {\Spec(R)}
	\arrow[from=1-2, to=2-1]
	\arrow[from=1-2, to=2-3]
	\arrow["{\eta_{f,x}}", from=2-1, to=2-3]
	\arrow[from=2-1, to=3-2]
	\arrow[from=2-3, to=3-2]
\end{tikzcd}\]
Forming cotangent complexes, we have the following definitions. 
\begin{definition}
    Given a map $f\colon \Spec(R) \ra X^M_{\DR}$, the local system $\cP^{\sf cot}_{X^M}(f)$ on $M$ with values in $\Mod_{R_{\sf red}}^{\LL_{R_{\sf red}| R}}$ is given by the composite 
    \[\Pi_\infty M \xra{\cP_{X^M}^{\sf geom}(f)} \ ^{\Spec(R_{\sf red})/}\dSt_{/\Spec(R)} \xra{\LL_{R_{\sf red}| (-)}} \Mod_{R_{\sf red}}^{\LL_{R_{\sf red}| R}}~.\]
\end{definition}

\begin{definition}
    Given a map $f\colon \Spec(R) \ra X^M_{\DR}$, the local system $\cA^{\sf cot}_{X}(f)$ on $M$ with values in $\Mod_{R_{\sf red}}^{\LL_{R_{\sf red}| R}}$ is given by the composite 
    \[\Pi_\infty M \xra{\cA_X^{\sf geom}(f)} \ ^{\Spec(R_{\sf red})/}\dSt_{/\Spec(R)} \xra{\LL_{R_{\sf red}| (-)}} \Mod_{R_{\sf red}}^{\LL_{R_{\sf red}| R}}~.\]
\end{definition}
The local system $\cP^{\sf cot}_{X^M}(f)$ is the constant local system with value $\LL_{R_{\sf red}| R} \ra \LL_{R_{\sf red}| X^M_R}$, and the local system $\cA^{\sf cot}_{X}(f)$ has fiber over $x\in M$ given by~$\LL_{R_{\sf red}| R} \ra\LL_{R_{\sf red}| X_{R,x}}$.
The map $\eta_{f, x}$ induces a map of local systems on $M$ with values in the category $\Mod_{R_{\sf red}}^{\LL_{R_{\sf red}| R}}$
\[\eta_{f, x}: \cA^{\sf cot}_{X}(f) \ra \cP^{\sf cot}_{X^M}(f)~.\]
By adjunction, we obtain a natural map of objects of the category~$\Mod_{R_{\sf red}}^{\LL_{R_{\sf red}| R}}$:
\[\underset{\Pi_\infty M}{\colim} \ \cA^{\sf cot}_{X}(f) \ra \cP^{\sf cot}_{X^M}(f)~.\]

\begin{lemma}
\label{relcot}
Given a map $f \colon \Spec(R) \ra X^M_{\DR}$, the natural map of objects of the category~$\Mod_{R_{\sf red}}^{\LL_{R_{\sf red}| R}}$
\[\underset{\Pi_\infty M}{\colim} \ \cA^{\sf cot}_{X}(f) \ra \cP^{\sf cot}_{X^M}(f)\]
is an equivalence. 
\end{lemma}
\begin{proof}
From the constructions of the local systems $\cP^{\sf geom}_{X^M}(f)$ and $\cA_X^{\sf geom}(f)$, we in fact have local systems of diagrams in $\dSt$ whose fiber over $x\in M$ are given by
\[\begin{tikzcd}
	{\Spec(R_{\sf red})} && \\
	& {X^M_R} & {X^M} \\
	& {\Spec(R)} & {X^M_{\DR}}
	\arrow["s"{description}, from=1-1, to=2-2]
	\arrow["r", from=1-1, to=2-3]
	\arrow[from=1-1, to=3-2]
	\arrow["g"', from=2-2, to=2-3]
	\arrow[from=2-2, to=3-2]
	\arrow[from=2-3, to=3-3]
	\arrow["f"', from=3-2, to=3-3]
\end{tikzcd}\]
and  
\[\begin{tikzcd}
	{\Spec(R_{\sf red})} && \\
	& {X_{R,x}} & X \\
	& {\Spec(R)} & {X_{\DR}}
	\arrow["{s_x}"{description}, from=1-1, to=2-2]
	\arrow["{r_x}", from=1-1, to=2-3]
	\arrow[from=1-1, to=3-2]
	\arrow["{g_x}"', from=2-2, to=2-3]
	\arrow[from=2-2, to=3-2]
	\arrow[from=2-3, to=3-3]
	\arrow["{f_x}"', from=3-2, to=3-3]
\end{tikzcd}\]
From the diagrams above, we have cofiber sequences of relative cotangent complexes
\[s^{*}\LL_{X^M_R| R} \rightarrow \LR \rightarrow \LL_{R_{\sf red}| X^M_R}\]
and
\[s_x^*\LL_{X_{R,x}|R} \rightarrow \LR \rightarrow \LL_{R_{\sf red}| X_{R,x}}.\]
Note that we have identifications 
\[
\LL_{X^M_R| R} \simeq g^*\LL_{X^M| X^M_{\DR}} \simeq g^*\LL_{X^M}
\]
and
\[
\LL_{X_{R,x}| R} \simeq g_x^*\LL_{X| X_{\DR}} \simeq g_x^*\LL_{X}~,
\]
where the second equivalence in both rows above follows from \cite{CPTVV} Lemma~2.1.10.
Therefore, we have 
\[\LL_{R_{\sf red}| X^M_R} \simeq {\sf cofib}(r^{*}\LL_{X^M} \rightarrow \LR)\]
and
\[\LL_{R_{\sf red}| X_{R,x}} \simeq {\sf cofib}(r^*_x\LL_{X} \to \LR).\]
Applying Lemma \ref{cot} to the complex $r^{*}\LL_{X^M}$, we have
\[\LL_{R_{\sf red}| X^M_R} \simeq {\sf cofib}\Bigl(\bigl(\underset{\Pi_\infty M}{\colim} \ r^*_x \LL_{X} \bigr)\rightarrow \LR\Bigr).\]
The forgetful functor 
\[(\Mod_{R_{\sf red}})_{/\LR} \rightarrow \Mod_{R_{\sf red}}\]
creates colimits, 
so we can compute $\underset{\Pi_\infty M}{\colim} \ r^*_x \LL_{X}$ in the over-category $(\Mod_{R_{\sf red}})_{/\LR}$. 
Now, note that the functor 
\[{\sf cofib}: (\Mod_{R_{\sf red}})_{/\LR} \rightarrow (\Mod_{R_{\sf red}})^{\LR/}\]
is colimit preserving. 
Therefore we have 
\[
\begin{split}
    \underset{\Pi_\infty M}{\colim} \ \cA^{\sf cot}_{X}(f) &\simeq \underset{\Pi_\infty M}{\colim} \  \LL_{R_{\rm{red}}| X_{R,x}} \\
    & \simeq \underset{\Pi_\infty M}{\colim} \ \sf{cofib}(r^*_x\LL_{X} \to \LR) \\
    & \simeq {\sf cofib}\Bigl(\bigl(\underset{\Pi_\infty M}{\colim} \ r^*_x \LL_{X} \bigr)\rightarrow \LR \Bigr) \\
    & \simeq \LL_{R_{\sf red}| X^M_R} \simeq  \cP^{\sf cot}_{X^M}(f)~.
\end{split}
\]
\end{proof}

\begin{remark}
    When $R$ is reduced, we have $\LR \simeq 0$ and the statement of the lemma simplifies to 
    \[\LL_{R| X^M_R}\simeq \underset{x \in M}\colim \ \LL_{R| X_{R,x}} \]
    where the colimit is computed in $\Mod_R$. The under-category only appears when one is at a non-reduced point.
 ~\hfill$\Diamond$\end{remark}

\begin{remark}
When $R = \Bbbk$, $X = \sB\GG$, and a map $f\colon \Spec(\Bbbk) \rightarrow \Loc_{\GG}(M)_{\DR}$ picks out a local system $L$ on~$M$, 
the fiber $\Loc_{\GG}(M)_{\Bbbk}$ is merely the formal completion of the point $L$ in~$\Loc_{\GG}(M)$. 
The $\Bbbk$-linear dual of relative cotangent complex $\LL_{\Spec(\Bbbk)| \Loc_{\GG}(M)_{k}}^{\vee}$ is then given by 
\[
\begin{split}
\LL_{\Spec(\Bbbk)| \Loc_{\GG}(M)_{k}}^{\vee} &\simeq \Gamma(\Pi_\infty M, \LL_{\Spec(\Bbbk)| \sB\GG_{k, x}}^{\vee}) \\
& \simeq \Gamma(\Pi_\infty M, \gg_{L}) \\
& \simeq \Omega^{\ast}(M, \gg_{L})
\end{split}
\]
where $\gg_{L}$ denotes the adjoint flat bundle associated to the $\GG$-local system $L$. 
This computation recovers the claim that the Lie algebra governing deformations of a $\GG$-local system on a manifold as the de Rham complex of the associated flat vector bundle.      
~\hfill$\Diamond$\end{remark}

\subsection{Local Description of Principal Parts} \label{subsectionalgebraiclocalsystem}
Recall the functor 
\[\DD: \dSt^{\op} \ra \CAlg(\Mod_{\Bbbk}^{\epsilon, \gr})\]
from Definition~\ref{gradedmixedfunctions}.
\begin{definition} \label{mainlocalsystem}
    The local system 
    \[\cA_X\colon \Pi_\infty M \ra \Fun\bigl(({\sf Aff}_{/X^M_{\DR}})^{\op}, \CAlg(\Mod_{\Bbbk}^{\epsilon, \gr})\bigr)\]
    is the composite
    \[\Pi_\infty M \xra{\cA_X^{\sf geom}} \Fun({\sf Aff}_{/X^M_{\DR}}, \dSt) \xra{\DD_{\ast}} \Fun\bigl(({\sf Aff}_{/X^M_{\DR}})^{\op}, \CAlg(\Mod_{\Bbbk}^{\epsilon, \gr})\bigr)~.\]
\end{definition}

Observe that the composites 
\[\Pi_\infty M \xra{\cP_{X^M}^{\sf geom}} \Fun({\sf Aff}_{/X^M_{\DR}}, \dSt) \xra{\DD_{\ast}} \Fun\bigl(({\sf Aff}_{/X^M_{\DR}})^{\op}, \CAlg(\Mod_{\Bbbk}^{\epsilon, \gr})\bigr)\]
and 
\[\Pi_\infty M \xra{\cD_{X^M_{\DR}}^{\sf geom}} \Fun({\sf Aff}_{/X^M_{\DR}}, \dSt) \xra{\DD_{\ast}} \Fun\bigl(({\sf Aff}_{/X^M_{\DR}})^{\op}, \CAlg(\Mod_{\Bbbk}^{\epsilon, \gr})\bigr)\]
recover the presheaf of principal parts and the crystalline structure sheaf of $X^M$, viewed as constant local systems on $M$. Note that we have
\[\Fun\bigl(({\sf Aff}_{/X^M_{\DR}})^{\op}, \CAlg(\Mod_{\Bbbk}^{\epsilon, \gr})\bigr) \simeq \CAlg\Bigl(\Fun\bigl(({\sf Aff}_{/X^M_{\DR}})^{\op}, \Mod_{\Bbbk}^{\epsilon, \gr}\bigr)\Bigr)\]
and we will view $\cA_X$ as an object of the latter category from hereon. The natural maps 
\[\cP^{\sf geom}_{X^M} \ra \cD_{X^M_{\DR}}^{\sf geom}\]
and 
\[\cA_X^{\sf geom} \ra \cD_{X^M_{\DR}}^{\sf geom}~\]
upgrade $\cA_X$ and $\cP_{X^M}$ to objects of the category~$\Fun\Bigl(\Pi_\infty M, \CAlg(\Mod^{\epsilon, \gr}_{\cD_{X^M_{\DR}}})\Bigr)$.
The following corollary is immediate from Corollary~\ref{compgeomlocsys}.
\begin{cor} \label{compalglocsys}
There is a natural map of objects of the category~$\Fun\Bigl(\Pi_\infty M, \CAlg(\Mod^{\epsilon, \gr}_{\cD_{X^M_{\DR}}})\Bigr)$:
\[\cA_X \ra \cP_{X^M}\]
where $\cP_{X^M}$ is viewed as a constant local system on $M$. 
\end{cor}
By adjunction, we obtain a natural map in $\CAlg(\Mod^{\epsilon, \gr}_{\cD_{X^M_{\DR}}})$:
\[\underset{\Pi_\infty M}{\colim} \ \cA_X \ra \cP_{X^M}~.\]

\begin{prop} \label{ppascolim}
    The natural map 
    \[\underset{\Pi_\infty M}{\colim} \ \cA_X \ra \cP_{X^M}\]
    is an equivalence of graded mixed $\cD_{X^M_{\DR}}$-algebras.
\end{prop}
To prove the proposition, we need the following lemma. 
\begin{lemma}
\label{nonmixlocsys}
Given a map $f\colon \Spec(R) \ra X^M_{\DR}$, we have a commutative diagram
\[
\begin{tikzcd}
\Pi_{\infty} M \ar[r, "{\cA_X^{\sf cot}}"] \ar[d, "{\cA_X(f)}"] 
& {(\Mod_{R_{\sf red}})^{\LR[-1]/}} \ar[r,"{(-)(1)}"] 
& {(\Mod^{\gr}_{R_{\sf red}})^{\LR[-1](1)/}} \ar[d,"{\Sym_{R_{\sf red}}}"]\\
\CAlg(\Mod^{\epsilon,\gr}_{\cD_{X^M_{\DR}}(f)}) \ar[rr,"{{\sf fgt}_{\epsilon}}"]
&
& {\CAlg(\Mod^{\gr}_{\Sym_{R_{\sf red}}(\LR[-1](1))})}
\end{tikzcd}
\]

In particular, we have an equivalence 
\[{\sf fgt}_{\epsilon} \circ \cA_X(f) \simeq \Sym_{R_{\sf red}} \circ (-)(1) \circ \cA_X^{\sf cot}(f)\]
of local systems on~$M$.
\end{lemma}

\begin{proof}
    \cite{CPTVV} Proposition 2.2.7 constructs an equivalence between the composites
    \[\dSt' \xra{\DD} \CAlg(\Mod_{\Bbbk}^{\epsilon, \gr}) \xra{\sf fgt_{\epsilon}} \CAlg(\Mod_{\Bbbk}^{\gr})\]
    and 
    \[\dSt' \xra{\LL_{(-)_{\sf red}\vert(-)}} \Mod_{\Bbbk}^{\epsilon, \gr} \xra{\Sym\bigl((-)[-1](1)\bigr)} \CAlg(\Mod^{\gr}_{\Bbbk})~,\]
    where $\dSt'$ denotes the full subcategory of $\dSt$ spanned by algebraisable affine formal derived stacks. Given an affine $\Spec(R)$, the equivalence above upgrades to an equivalence between the composites 
    \[^{\Spec(R_{\sf red})_{/}}(\dSt')_{^{/}\Spec(R)} \ra  \CAlg\Bigl(\Mod^{\epsilon, \gr}_{\DD(R)}\Bigr) \ra \CAlg\Bigl(\Mod^{\gr}_{\Sym_{R_{\sf red}}\bigr(\LL_{R_{\sf red}\vert R}[-1](1)\bigl)}\Bigr)\]
    and 
    \[^{\Spec(R_{\sf red})_{/}}(\dSt')_{^{/}\Spec(R)}\ra  (\Mod_{R_{\sf red}})^{\LR[-1]/} \ra \CAlg\Bigl(\Mod^{\gr}_{\Sym_{R_{\sf red}}\bigr(\LL_{R_{\sf red}\vert R}[-1](1)\bigl)}\Bigr)\] 
    
    By \cite{CPTVV} Proposition 2.1.8, for $X$ derived Artin, locally of finite presentation over $\Bbbk$, fibers of the family $X \ra X_{\DR}$ along maps $\Spec(R) \ra X_{\DR}$ are algebraisable affine formal derived stacks. Also note that $(X_R)_{\sf red} \simeq \Spec(R_{\sf red})$. Therefore the local system 
    \[\cA_X^{\sf geom}(f)\colon \ \Pi_{\infty}M \ra ^{\Spec(R_{\sf red})_{/}}(\dSt)_{^{/}\Spec(R)}\]
    lifts to a local system 
    \[\cA_X^{\sf geom}(f)\colon \ \Pi_{\infty}M \ra ^{\Spec(R_{\sf red})_{/}}(\dSt')_{^{/}\Spec(R)}~,\]
    and we have the result. 
    
\end{proof}

\begin{proof}[Proof of Proposition~\ref{ppascolim}]
It suffices to check that for each $f\colon \Spec(R) \ra X^M_{\DR}$
\[\Bigl(\underset{\Pi_\infty M}{\colim} \ \cA_X\Bigr)(f) \longrightarrow \cP_{X^M}(f)\]
is an equivalence of graded mixed $\cD_{X^M_{\DR}}(f)\simeq \DD(R)$-algebras. 
As colimits are computed pointwise, it suffices to check that 
\[\left(\underset{\Pi_\infty M}{\colim}  \ \cA_X(f)\right) \longrightarrow \cP_{X^M}(f)\]
is an equivalence of $\DD(R)$-algebras. 
The functor
\[\Mod^{\epsilon, \gr}_{\DD(R)} \xrightarrow[]{{\sf fgt}_{\epsilon}} \Mod^{\gr}_{\Sym_{R_{\sf red}}(\LR[-1](1))}\]
is symmetric monoidal, colimit preserving, and creates equivalences. Therefore the induced functor 
\[\CAlg\bigl(\Mod^{\epsilon, \gr}_{\DD(R)}\bigr) \xrightarrow[]{{\sf fgt}_{\epsilon}} \CAlg\bigl(\Mod^{\gr}_{\Sym_{R_{\sf red}}(\LR[-1](1))}\bigr)\]
is also colimit-preserving and creates equivalences. 
Hence it suffices to check that the map 
\[\underset{\Pi_\infty M}{\colim}  \ {\sf fgt}_{\epsilon}(\cA_X(f)) \longrightarrow {\sf fgt}_{\epsilon}(\cP_{X^M}(f)) \simeq \Sym_{R_{\sf red}}(\LL_{R_{\sf red}| X^M_R}[-1](1))\]
is an equivalence. 
By Lemma~\ref{nonmixlocsys}, the map above is equivalent to the map 
\[\underset{\Pi_\infty M}{\colim} \ \Sym_{R_{\sf red}}\bigl(\LL_{R_{\sf red} | X_{R,x}}[-1](1)\bigr) \longrightarrow \Sym_{R_{\sf red}}\bigl(\LL_{R_{\sf red}| X^M_R}[-1](1)\bigr)~.\]
However, as 
\[\Sym_{R_{\sf red}}\colon (\Mod^{\gr}_{R_{\sf red}})^{\LR[-1](1)/} \longrightarrow \CAlg(\Mod^{\gr}_{\Sym_{R_{\sf red}}(\LR[-1](1))})\]
is a left adjoint, we have 
\[\underset{\Pi_\infty M}{\colim}  \ \Sym_{R_{\sf red}}\bigl(\LL_{R_{\sf red} | X_{R,x}}[-1](1)\bigr) \simeq \Sym_{R_{\sf red}}\bigl(\underset{\Pi_\infty M}{\colim} \ \LL_{R_{\sf red} | X_{R,x}}[-1](1)\bigr)\]
where the colimit on the right is computed in the category $(\Mod^{\gr}_{R_{\sf red}})^{\LR[-1](1)/}$.
By Lemma~\ref{relcot}, we have 
\[
\underset{\Pi_\infty M}{\colim} \ \LL_{R_{\sf red}| X_{R,x}}[-1](1) \simeq \LL_{R_{\sf red}| X^M_R}[-1](1).
\]
and we have the result.
\end{proof}

\section{Factorization Homology} \label{factorization homology}

In this section we review factorization homology with a focus on features pertinent to our main technical result. 
Let $\cV$ be a presentably symmetric monoidal category, and $M$ a smooth $d$-manifold. 
Recall the notion of a $\Disk_M$-algebra on a smooth $d$-manifold $M$ with values in $\cV$ from~\cite{oldfact}. 
A local system 
\[M \xra{A} \CAlg(\cV)\]
yields a $\Disk_M$-algebra by the following lemma. 
 \begin{lemma} \label{diskfromcommloc}
     There is a functor 
     \[\Fun(\Pi_{\infty}M, \CAlg(\cV)) \longrightarrow \Alg_{\Disk_M}(\cV)\]
     which for each point of $M$ lies over the forgetful functor $\CAlg(\cV)\ra \Alg_{\EE_{n+1}}(\cV)$.
 \end{lemma}
 \begin{proof}
By~\cite{zp}, there is an equivalence of operads
\[
\underset{U\in\Disk_{M}}\colim \Disk_{U} \overset{\sim}\longrightarrow\Disk_{M}
\]
that gives an equivalence of categories
\[
\Alg_{\Disk_M}(\cV)\overset{\sim}\longrightarrow\underset{U\in \Disk_M}\limit \Alg_{\Disk_U}(\cV)~.
\]
The functor $U\mapsto \Alg_{\Disk_U}(\cV)$ is the right Kan extension of its values on connected opens. Therefore we have 
\[\Alg_{\Disk_M}(\cV) \simeq \underset{U \in\Disk_M^{\sf conn}}{\lim} \Alg_{\Disk_{U}}(\cV)~.\]
For $U \in \Disk^{\sf conn}_M$, the canonical maps of operads
\[\Disk_U \longrightarrow \EE_{\infty}\]
induce canonical maps 
\[\CAlg(\cV) \longrightarrow \Alg_{\Disk_U}(\cV)~.\]
Taking limits, we obtain
\[
    \Fun(\Disk_M^{\sf conn}, \CAlg(\cV)) \simeq \underset{\Disk_M^{\sf conn}}{\lim}\CAlg(\cV) \longrightarrow \underset{U \in \Disk_M^{\sf conn}}{\lim} \Disk_{U}(\cV) \simeq \Alg_{\Disk_M}(\cV)~.
\]
Now, by \cite{oldfact} Lemma 2.12, we have an equivalence $\Pi_{\infty}M\simeq\Disk_M^{\sf conn}$. Therefore, we obtain a functor 
\[\Fun(\Pi_{\infty}M, \CAlg(\cV)) \longrightarrow \Alg_{\Disk_M}(\cV)~.\]
 \end{proof}
As we will see in the next section of the paper, our main $\Disk_M$-algebra of interest will be presented to us as a local system of commutative algebras. The factorization homology of a $\Disk_M$-algebra obtained from a local system of commutative algebras has an easy description. Given a local system 
\[M \xra{A} \CAlg(\cV)\]
let $\cA \in \Alg_{\Disk_M}(\cV)$ also denote the corresponding $\Disk_M$-algebra. Then we have 
\[\int_MA \simeq \underset{\Pi_{\infty}M}{\colim} \ A\]
where the colimit on the right is computed in the category $\CAlg(\cV)$. 
\begin{lemma} \label{colim}
The following diagram commutes 
\[\begin{tikzcd}
	{\Fun(\Pi_{\infty}M, \CAlg(\cV))} & {\Alg_{\Disk_M}(\cV)} \\
	{\CAlg(\cV)} & \cV
	\arrow[from=1-1, to=1-2]
	\arrow["\colim"', from=1-1, to=2-1]
	\arrow["{\int_M}", from=1-2, to=2-2]
	\arrow["{{\sf fgt}}"', from=2-1, to=2-2]
\end{tikzcd}\]
\end{lemma}
\begin{proof}
Observe that $\CAlg(\CAlg(\cV)) \simeq \CAlg(\cV)$. Therefore, consider the functor 
\[\Fun(\Pi_{\infty}M, \CAlg(\cV)) \simeq \Fun(\Pi_{\infty}M, \CAlg(\CAlg(\cV))) \longrightarrow \Alg_{\Disk_M}(\CAlg(\cV))~,\]
where the second arrow is obtained by Lemma \ref{diskfromcommloc}. 
The forgetful functor 
\[{\sf fgt}: \CAlg(\cV)\rightarrow \cV\]
is symmetric monoidal, and therefore induces a functor 
\[{\sf fgt}: \Alg_{\Disk_M}\bigl(\CAlg(\cV)\bigr) \longrightarrow \Alg_{\Disk_M}\bigl(\cV\bigr)~.\]
Observe that the diagram 
\[\begin{tikzcd}
	{\Fun(\Pi_{\infty}M, \CAlg(\cV))} & {\Fun(\Pi_{\infty}M, \CAlg(\CAlg(\cV)))} & {\Alg_{\Disk_M}(\CAlg(\cV))} \\
	&& {\Alg_{\Disk_M}(\cV)}
	\arrow["\simeq", from=1-1, to=1-2]
	\arrow[from=1-1, to=2-3]
	\arrow[from=1-2, to=1-3]
	\arrow["{{\sf fgt}}", from=1-3, to=2-3]
\end{tikzcd}\]
commutes. Note that ${\sf fgt}$ preserves sifted colimits (\cite{HA} Corollary 3.2.3.2), 
and therefore the following diagram commutes
\[\begin{tikzcd}
	{{\Fun(\Pi_{\infty}M, \CAlg(\CAlg(\cV)))}} & {{\Fun(\Pi_{\infty}M, \CAlg(\cV))}} \\
	{\Alg_{\Disk_M}(\CAlg(\cV))} & {\Alg_{\Disk_M}(\cV)} \\
	{\CAlg(\cV)} & \cV
	\arrow["{\simeq }", from=1-1, to=1-2]
	\arrow[from=1-1, to=2-1]
	\arrow[from=1-2, to=2-2]
	\arrow["{\sf fgt}", from=2-1, to=2-2]
	\arrow["{\int_M}", from=2-1, to=3-1]
	\arrow["{\int_M}", from=2-2, to=3-2]
	\arrow["{{\sf fgt}}", from=3-1, to=3-2]
\end{tikzcd}\]
Now, given an object $A \in \Fun(\Pi_{\infty}M, \CAlg(\cV)) \simeq \Fun(\Pi_{\infty}M, \CAlg(\CAlg(\cV)))$, the corresponding object $A \in \Alg_{\Disk_M}(\CAlg(\cV))$ is the left Kan extension of its values on $\Disk_M^{\sf conn}$ (note that the coproduct on $\CAlg(\cV)$ is the symmetric monoidal product). Therefore we have 
\[
\begin{split}
    \int_M A &\simeq \underset{\Disk_M}{\colim} \ A \\
    &\simeq \underset{\Disk_M^{\sf conn}}{\colim} \ A \\
    &\simeq \underset{\Pi_{\infty }M}{\colim} \ A
\end{split}
\]
where the colimits are being comuted in $\CAlg(\cV)$, i.e., the following diagram commutes 
\[\begin{tikzcd}
	{\Fun(\Pi_{\infty}M, \CAlg(\CAlg(\cV)))} & {\Alg_{\Disk_M}(\CAlg(\cV))} \\
	{} & {\CAlg(\cV)}
	\arrow[from=1-1, to=1-2]
	\arrow["\colim"', from=1-1, to=2-2]
	\arrow["{\int_M}", from=1-2, to=2-2]
\end{tikzcd}\]
Therefore the following diagram commutes 
\[\begin{tikzcd}
	{{\Fun(\Pi_{\infty}M, \CAlg(\cV))}} & {{\Fun(\Pi_{\infty}M, \CAlg(\CAlg(\cV)))}} & {{\Fun(\Pi_{\infty}M, \CAlg(\cV))}} \\
	& {\Alg_{\Disk_M}(\CAlg(\cV))} & {\Alg_{\Disk_M}(\cV)} \\
	& {\CAlg(\cV)} & \cV
	\arrow["\colim"', from=1-1, to=3-2]
	\arrow["\simeq"', from=1-2, to=1-1]
	\arrow["{\simeq }", from=1-2, to=1-3]
	\arrow[from=1-2, to=2-2]
	\arrow[from=1-3, to=2-3]
	\arrow["{\sf fgt}", from=2-2, to=2-3]
	\arrow["{\int_M}", from=2-2, to=3-2]
	\arrow["{\int_M}", from=2-3, to=3-3]
	\arrow["{{\sf fgt}}", from=3-2, to=3-3]
\end{tikzcd}\]
and we have the result. 

\end{proof}
In Section~7, we will construct deformations of $\Disk_M$-algebras arising from local systems of commutative algebras. These deformations will arise from a deformation of a local system of commutative algebras as a local system of $\EE_{n+1}$-algebras. A local system of $\EE_{n+1}$-algebras does not give rise to a $\Disk_M$-algebra immediately--a choice of {\em framing} of $M\times\RR^{n+1-d}$ is required. The following lemma makes this precise. 
\begin{lemma}
\label{diskasloc}
For $n+1 \geq d$, a framing of $M\times\RR^{n+1-d}$ determines an equivalence
\[
\Alg_{\Disk_M}\bigl(\Alg_{\EE_{n+1-d}}(\cV)\bigr)\simeq \Fun\bigl(\Pi_\infty M,\Alg_{\EE_{n+1}}(\cV)\bigr)
\]
between $\Disk_M$-algebras and local systems of $\EE_n$-algebras on~$M$.
\end{lemma}

\begin{proof}
By~\cite{zp}, there is an equivalence of operads
\[
\underset{U\in\Disk_{M}}\colim \Disk_{U} \overset{\sim}\longrightarrow\Disk_{M}
\]
that gives an equivalence of categories
\[
\Alg_{\Disk_M}\bigl(\Alg_{\EE_{n+1-d}}(\cV)\bigr)\overset{\sim}\longrightarrow\underset{U\in \Disk_M}\limit \Alg_{\Disk_U}\bigl(\Alg_{\EE_{n+1-d}}(\cV)\bigr)~.
\]
The functor $U\mapsto \Alg_{\Disk_U}\bigl(\Alg_{\EE_{n+1-d}}(\cV)\bigr)$ is the right Kan extension of its values on connected opens. 
Consequently, this limit can be computed after restricting to connected opens, which gives
\[
\underset{U\in \Disk_M}\limit \Alg_{\Disk_U}\bigl(\Alg_{\EE_{n+1-d}}(\cV)\bigr)
\simeq
\underset{U \in \Disk_M^{\sf conn}}\limit \Alg_{\Disk_U}\bigl(\Alg_{\EE_{n+1-d}}(\cV)\bigr)~.
\]
Observe that a choice of identification $\RR^d \cong U$ determines an equivalence of categories 
\[\Alg_{\Disk_U}(\Alg_{\EE_{n+1-d}}(\cV)\bigr) \simeq \Alg_{\EE_{d}}\bigl(\Alg_{\EE_{n+1-d}}(\cV)\bigr) \simeq \Alg_{\EE_{n+1}}(\cV)~.\]
Consider the functor 
\[\Disk_M^{\sf conn} \xra{U \mapsto {\Alg_{\Disk(U)}}\bigl(\Alg_{\EE_{n+1-d}}(\cV)\bigr)} \Cat_{\infty}~.\]

We have an identification $\Pi_{\infty}M\simeq\Disk_M^{\sf conn}$ of $\infty$-groupoids (\cite{oldfact} Lemma 2.12)
and so we have a functor
\[
\Pi_{\infty} M \simeq \Disk_M^{\sf conn} \xra{U \mapsto {\Alg_{\Disk(U)}}\bigl(\Alg_{\EE_{n+1-d}}(\cV)\bigr)} \Cat_{\infty}~,
\]
i.e., we have a local system of categories on $\Pi_{\infty}M$. The local system above is classified by a map 
\[\Pi_{\infty}M \xra{\tau_M} \sB\sO(d) \xra{-\oplus_{{\sf id}_{\RR^{n+1-d}}}} \sB\sO(n+1) \ra \Cat_{\infty} \]
where $\tau_M$ classifies the tangent bundle of $M$, and $\sB\sO(n+1) \ra \Cat_{\infty}$ is the $\sO(n+1)$-action on $\Alg_{\EE_{n+1}}(\cV)$. A framing of $M\times\RR^{n+1-d}$, i.e., a nullhomotopy of $\Pi_{\infty}M \ra \sB\sO(d) \ra \sB\sO(n+1)$, provides a trivialization of this functor. 
Now the limit of this functor is equivalent to the category of sections of the unstraightening,
which is
\[
\Gamma_M(\Pi_\infty M\times \Alg_{\EE_{n+1}}(\cV)) \simeq \Fun(\Pi_\infty M, \Alg_{\EE_{n+1}}(\cV))~,
\]
and so we recover the claim.
\end{proof}

\begin{remark}
    The Lemma above does not require a framing of $M$. For example, given a local system of $\EE_3$-algebras on a closed oriented surface $\Sigma$, a choice of framing on $\Sigma\times \RR$ determines the corresponding $\Disk_{M}$-algebra. The manifold $\Sigma$ does not admit a framing in general, but $\Sigma\times \RR$ always does. 
~\hfill$\Diamond$\end{remark}

\begin{remark}
    Given a local system of $\EE_{n}$-algebras $A$ on $M$, and a framing of $M\times\RR^{n+1-d}$, we obtain a $\Disk_M$-algebra valued in $\EE_{n+1-d}$-algebras, which we also denote by $A$. Therefore one obtains the $\EE_{n+1-d}$-algebra 
    \[\int_M A~.\]
    Unlike what happens in Lemma~\ref{colim}, the object $\int_M A$ is sensitive to a lot more than just the homotopy type of $M$; it is sensitive to the homotopy types of all of the configuration spaces ${\sf Conf}_k(M)$ for $k\geq 0$. Therefore, the various ``quantum" constructions made in this article should know a lot more about the manifold $M$ than just its homotopy type. 
~\hfill$\Diamond$\end{remark}

\section{Principal Parts on Mapping Stacks as Factorization Homology}\label{Principal Parts on Mapping Stacks as Factorization Homology}
Throughout this section, we fix $X$ a derived Artin stack, locally of finite presentation over $\Bbbk$, and fix $M$ a $d$-manifold with the homotopy type of a finite CW complex. No framings are required for the results of this section. 

Recall the local system $\cA_X \in \Fun\Bigl(\Pi_\infty M, \CAlg(\Mod^{\epsilon, \gr}_{\cD_{X^M_{\DR}}})\Bigr)$ from Definition \ref{mainlocalsystem}. 
Its image under the composite 
\[\Fun\Bigl(\Pi_\infty M, \CAlg(\Mod^{\epsilon, \gr}_{\cD_{X^M_{\DR}}})\Bigr) \underset{{\sf Lemma \ } \ref{diskfromcommloc}}{\longrightarrow}  \Alg_{\Disk_M}\Bigl(\Mod^{\epsilon, \gr}_{\cD_{X^M_{\DR}}}\Bigr)\]
gives a $\Disk_M$-algebra also denoted by $\cA_X$.

\begin{remark} \label{unpackingdiskalgebra}
We now describe the value of $\cA_X$ on a disk $\RR^d \hookrightarrow M$. Recall the proof of Lemma~\ref{diskasloc}. As the $\Disk_M$-algebra $\cA_X$ is constructed from the local system of commutative algebras $\cA_X$, its value on a disk $\iota \colon \RR^d \hookrightarrow M$ is identified with the fiber of the local system $\cA_X$ at $\iota(0) \in M$, i.e., we have an equivalence in $\Mod^{\epsilon, \gr}_{\cD_{X^M_{\DR}}}$:
\[\cA_X(\iota\colon \RR^d \hookrightarrow M) \simeq \cA_X(\iota(0))~.\]

Given an embedding of several disks $\iota\colon \bigsqcup_{i=1}^k\RR^d \hookrightarrow M$, let $x_i$ denote $0$ in the $i$th copy of $\RR^d$. We have an equivalence in $\Mod^{\epsilon, \gr}_{\cD_{X^M_{\DR}}}$ (note that the tensor product below is over $\cD_{X^M_{\DR}}$):
\[\cA_X\Bigl(\iota\colon \bigsqcup_{i=1}^k\RR^d \hookrightarrow M\Bigr) \simeq \bigotimes_{i=1}^k \cA_X(\iota(x_i))~.\]
Note that the objects above are presheaves of graded mixed complexes on ${\sf Aff}_{/X^M_{\DR}}$ (in fact, algebras over the crystalline structure sheaf $\cD_{X^M_{\DR}}$). Therefore, given a map $f\colon \Spec(R) \rightarrow X^M_{\DR}$, we obtain equivalences of graded mixed $\cD_{X^M_{\DR}}(f)$-modules  
\begin{align*}
    \cA_X\Bigl(\iota\colon \bigsqcup_{i=1}^k\RR^d \hookrightarrow M\Bigr)(f) & \simeq \Bigl(\bigotimes_{i=1}^k \cA_X(\iota(x_i))\Bigr)(f) \\
    & \simeq \bigotimes_{i=1}^{k} \bigl(\cA_X(\iota(x_i))\bigr)(f) \\
    & \simeq \bigotimes_{i=1}^{k} \cP_{X}(\ev_{\iota(x_i)}\circ f)~,
\end{align*}
where the second equivalence follows from the fact that tensor products in $\Mod^{\epsilon, \gr}_{\cD_{X^M_{\DR}}}$ are computed pointwise.
~\hfill$\Diamond$
\end{remark}

Its factorization homology 
\[
\int_M\cA_X \in \Mod^{\epsilon, \gr}_{\cD_{X^M_{\DR}}}
\]
admits a rich structure in the setting of derived geometry.
As tensor products and colimits in $\Mod^{\epsilon, \gr}_{\cD_{X^M_{\DR}}}$ are computed pointwise, 
for any point $(f \colon\Spec(R) \rightarrow X^M_{\DR})$,
we compute $\int_M\cA_X(f)$ in the category $\Mod^{\epsilon, \gr}_{\cD_{X^M_{\DR}(f)}}$.
Thus $\int_M\cA_X$ yields a presheaf on~$X^M_{\DR}$:
\[
\begin{array}{ccccc}
\int_M \cA_X\colon & \bigl({\sf Aff}/X^M_{\DR}\bigr)^{\sf op} & \longrightarrow & \Mod_{\Bbbk}^{\epsilon, \gr} \\
& (\Spec(R) \xrightarrow{f}X^M_{\DR}) & \longmapsto & \int_M\cA_X(f)
\end{array}
\]
and we will now compare the presheaf $\int_M \cA_X$ with the presheaf of principal parts~$\cP_{X^M}$. 
\begin{prop}
\label{ppasfh}
There is an equivalence of graded mixed $\cD_{X^M_{\DR}}$-modules 
\[\int_M \cA_X \xrightarrow{\simeq} \cP_{X^M}~.\]
\end{prop}

\begin{proof}
We have a chain of equivalences of graded mixed $\cD_{X^M_{\DR}}$-modules
\[
\int_M\cA_X \simeq \underset{\Pi_\infty M}{\colim} \ \cA_X \simeq \cP_{X^M}
\]
where the first equivalence follows from Lemma~\ref{colim} and the second equivalence follows from Lemma~\ref{ppascolim}.
\end{proof}

\begin{remark}\label{failureforordinaryprincipalparts}
    It is essential to use {\em mixed} principal parts in the computation above. More precisely, ordinary principal parts is the $\cO_{X_{\DR}}$-algebra $q_*\cO_{X}$, where $q\colon X\ra X_{\DR}$ is the canonical map. Given a $d$-manifold $M$, we can similarly produce a sheaf of $\Disk_M$-algebras $\cA_X^{\sf ord}$. It is not true in general that 
    \[\int_M \cA_X^{\sf ord} \simeq q'_*\cO_{X^M}\]
    where $q' \colon X^M \ra X^M_{\DR}$ is the canonical map. For example, take $X = \AA^1$ and $M = \RR^{d} \sqcup \RR^d$. Note that the mapping stack $\Map(\RR^{d} \sqcup \RR^d, \AA^1) \simeq \AA^1 \times \AA^1 \simeq \AA^2$
    Consider the $\Bbbk$-point $f \colon \pt \ra \AA^2_{\DR}$ picking out the origin. We have 
    \[q'_*\cO_{X^M}(f) \simeq \Bbbk[\![x, y]\!]\]
    power series in two variables in cohomological degree $0$. The $\Disk_M$-algebra $\cA_X^{\sf ord}$ corresponds to the trivial local system on $\RR^{d} \sqcup \RR^d$ with fiber $\Bbbk[\![x]\!]$ where $x$ is of degree $0$. The factorization homology is given by
    \[\int_M \cA_X^{\sf ord} \simeq \Bbbk[\![x]\!]\otimes \Bbbk[\![y]\!]~,\]
    which is not equivalent to $\Bbbk[\![x, y]\!]$. The formalism of graded mixed complexes deals with issues involving completions.     
~\hfill$\Diamond$\end{remark}
\begin{remark}
    It is also essential that we work $\cD_{X^M_{\DR}}$-linearly. This is visible at non-reduced points: see the proof of Lemma~\ref{ppascolim}. 
~\hfill$\Diamond$\end{remark}

\section{Quantizing Mapping Stacks} \label{Quantizing Mapping Stacks}

We will now describe an approach to quantizing mapping stacks. 

Let $X$ be a derived Artin stack, locally of finite presentation over $\Bbbk$, $M$ a $d$-manifold, and $X^M = \Map(\Pi_\infty M, X)$ the associated mapping stack. Given an $n$-shifted Poisson structure on $X$, for $d \leq n+1$, and a framing on $M\times \RR^{n+1-d}$, we will produce an $\EE_{n+1-d}$-quantization of the algebra of principal parts of the mapping stack $X^M$. More precisely, given $\pi \in {\sf Poiss}(X,n)$, $n\geq 1$, we produce the $\Disk_M$-algebra of {\em quantized local principal parts}
\[\cA^{\pi}_X[\![\hbar]\!] \in \Alg_{\Disk_M}\Bigl(\Alg_{\EE_{n+1-d}}\bigl(\Mod^{\epsilon, \gr}_{\cD_{X^M_{\DR}}[\![\hbar]\!]}\bigr)\Bigr)~.\]
As Proposition~\ref{ppasfh} provides an equivalence of $\cD_{X^M_{\DR}}$-modules 
\[\cP_{X^M} \simeq \int_M \cA_X~,\]
a deformation of $\cA_X$ immediately determines a deformation of $\cP_{X^M}$: we have an $\EE_{n+1-d}$-quantization of principal parts on the mapping stack 
\[\cP^{\pi}_{X^M}[\![\hbar]\!] := \int_M\cA^{\hbar, \pi}_{X} \in \Alg_{\EE_{n+1-d}}\bigl(\Mod^{\epsilon, \gr}_{\cD_{X^M_{\DR}}[\![\hbar]\!]}\bigr)~.\]

\begin{remark}
    The algebra of principal parts $\cP_{X^M}$ only sees the underlying $\infty$-groupoid of $M$.  The algebra of quantized principal parts 
    \[\cP_{X^M}^{\pi}[\![\hbar]\!] := \int_M\cA_X^{\pi}[\![\hbar]\!]\]
    should be sensitive to a lot more information about $M$: the factorization homology $\int_M \cA_X^{\pi}[\![\hbar]\!]$ should be sensitive to the homotopy types of all of the configuration spaces of points ${\sf Conf}_k(M)$. 
~\hfill$\Diamond$\end{remark}

We will now define the $\Disk_M$-algebra $\cA_X^{\pi}[\![\hbar]\!]$. Recall the local system of commutative algebras
\[\cA_X\colon \Pi_{\infty}M \ra \CAlg\bigl(\Mod^{\epsilon, \gr}_{\cD_{X^M_{\DR}}}\bigr)\]
associated to the mapping stack $X^M$. 
Given a shifted Poisson structure $\pi \in {\sf Poiss}(X, n)$, by Lemma~\ref{poissbivecttodef}, we obtain a $\cD_{X_{\DR}}[\![\hbar]\!]$-linear $\PP_{n+1}$-algebra structure on $\cP_{X}[\![\hbar]\!]$. For $n\geq 1$, by formality of the $\EE_{n+1}$-operad, we obtain a functor 
\[\cP_X^{\pi}[\![\hbar]\!]\colon \bigl({\sf Aff}_{/X_{\DR}}\bigr)^{\op} \ra \Alg_{\EE_{n+1}}\bigl(\Mod_{\Bbbk[\![\hbar]\!]}^{\epsilon, \gr}\bigr)~.\]

\begin{definition}\label{localsystemfrompoiss}
Let $X$ be a derived Artin stack, locally of finite presentation over $\Bbbk$. Let $M$ be a $d$-manifold. Given a shifted Poisson structure $\pi \in {\sf Poiss}(X, n)$, the local system 
\[\cA_X^{\pi}[\![\hbar]\!]\colon \Pi_{\infty}M \ra \Alg_{\EE_{n+1}}\Bigl(\Fun\bigl(({\sf Aff}_{/X^M_{\DR}})^{\op}, \Mod_{\Bbbk[\![\hbar]\!]}^{\epsilon, \gr}\bigr)\Bigr)\]
is the adjoint of the composite 
\[M\times \bigl({\sf Aff}_{/X^M_{\DR}}\bigr)^{\op} \xra{\sf ev} \bigl({\sf Aff}_{/X_{\DR}}\bigr)^{\op} \xra{\cP_X^{\pi}[\![\hbar]\!]} \Alg_{\EE_{n+1}}\bigl(\Mod_{\Bbbk[\![\hbar]\!]}^{\epsilon, \gr}\bigr)~.\]
\end{definition}
Observe that the local system $\cA_X^{\pi}[\![\hbar]\!]$ is a module over $\cD_{X^M_{\DR}}$, more precisely, it is an object of $\Fun\Bigl(\Pi_{\infty}M, \Alg_{\EE_{n+1}}\bigl(\Mod^{\epsilon, \gr}_{\cD_{X^M_{\DR}}[\![\hbar]\!]}\bigr)\Bigr)$. Now given a framing of $M\times\RR^{n+1-d}$, we obtain a $\Disk_M$-algebra. 

\begin{definition} \label{diskalgfrompoiss}
    Let $X$ be a derived Artin stack, locally of finite presentation over $\Bbbk$. Let $M$ be a $d$-manifold. Given a shifted Poisson structure $\pi \in {\sf Poiss}(X, n)$, and a framing of $M\times \RR^{n+1-d}$, the algebra of {\em quantized local principal parts} is the object 
    \[\cA_X^{\pi}[\![\hbar]\!] \in \Alg_{\Disk_M}\Bigl(\Alg_{\EE_{n+1-d}}\bigl(\Mod^{\epsilon, \gr}_{\cD_{X^M_{\DR}}[\![\hbar]\!]}\bigr)\Bigr)\]
    corresponding to (Definition~\ref{localsystemfrompoiss})
    \[\cA_X^{\pi}[\![\hbar]\!] \in \Fun\Bigl(\Pi_{\infty}M, \Alg_{\EE_{n+1}}\bigl(\Mod^{\epsilon, \gr}_{\cD_{X^M_{\DR}}[\![\hbar]\!]}\bigr)\Bigr)\]
    under the equivalence (Lemma~\ref{diskasloc})
    \[\Alg_{\Disk_M}\Bigl(\Alg_{\EE_{n+1-d}}\bigl(\Mod^{\epsilon, \gr}_{\cD_{X^M_{\DR}}[\![\hbar]\!]}\bigr)\Bigr)\simeq \Fun\Bigl(\Pi_{\infty}M, \Alg_{\EE_{n+1}}\bigl(\Mod^{\epsilon, \gr}_{\cD_{X^M_{\DR}}[\![\hbar]\!]}\bigr)\Bigr)\]
\end{definition}

Note that by Lemma~\ref{poissbivecttodef}, we have an equivalence of $\cD_{X_{\DR}}$-linear $\EE_{n+1}$-algebras
\[\cP_X^{\pi}[\![\hbar]\!] \underset{\cD_{X_{\DR}}[\![\hbar]\!]}{\otimes}\cD_{X_{\DR}} \simeq \cP_X~.\]
Unwinding definitions, we also have 
\[\cA_X^{\pi}[\![\hbar]\!]\underset{\cD_{X^M_{\DR}}[\![\hbar]\!]}\otimes\cD_{X^M_{\DR}} \simeq \cA_X~.\]

Taking factorization homology, we obtain a quantum object on $X^M_{\DR}$ associated to an $n$-shifted Poisson structure on $X$. 

\begin{definition}\label{quantumprincipalparts}
    Let $X$ be a derived Artin stack, locally of finite presentation over $\Bbbk$. Let $M$ be a $d$-manifold. Given a shifted Poisson structure $\pi \in {\sf Poiss}(X, n)$, and a framing of $M\times \RR^{n+1-d}$, the algebra of {\em quantized principal parts on $X^M$} is 
    \[\cP_{X^M}^{\pi}[\![\hbar]\!]\colon = \int_M \cA_X^{\pi}[\![\hbar]\!] \]
an object in $\Alg_{\EE_{n+1-d}}\bigl(\Mod^{\epsilon, \gr}_{\cD_{X^M_{\DR}}[\![\hbar]\!]}\bigr)$.
\end{definition}

We view the presheaf  $\cP_{X^M}^{\pi}[\![\hbar]\!]$ of $\cD_{X^M_{\DR}}[\![\hbar]\!]$-modules as a quantization of the presheaf of principal parts~$\cP_{X^M}$. 
More punchily, it quantizes the mapping stack.

\begin{prop}
\label{defpp}
Let $X$ be a derived Artin stack, locally of finite presentation over $\Bbbk$. Let $M$ be a $d$-manifold. Given a shifted Poisson structure $\pi \in {\sf Poiss}(X, n)$, and a framing of $M\times \RR^{n+1-d}$, there is an equivalence 
\[
\cP_{X^M}^{\pi}[\![\hbar]\!] \underset{\cD_{X^M_{\DR}}[\![\hbar]\!]}\otimes\cD_{X^M_{\DR}} \simeq \cP_{X^M}
\]
in $\Alg_{\EE_{n+1-d}}\bigl(\Mod^{\epsilon, \gr}_{\cD_{X^M_{\DR}}}\bigr)$.
\end{prop}

\begin{proof}
We have $\cP_{X^M}^{\pi}[\![\hbar]\!]\colon = \int_M \cA_X^{\pi}[\![\hbar]\!]$. As reduction mod $\hbar$ is colimit-preserving and symmetric monoidal, we have 
\[\cP_{X^M}^{\pi}[\![\hbar]\!] \underset{\cD_{X^M_{\DR}}[\![\hbar]\!]}\otimes\cD_{X^M_{\DR}}  \simeq \int_M \cA_X^{\pi}[\![\hbar]\!]\underset{\cD_{X^M_{\DR}}[\![\hbar]\!]}\otimes\cD_{X^M_{\DR}} \simeq \int_M\cA_X \simeq \cP_{X^M}~.\]
where the last equivalence is Proposition~\ref{ppasfh}.
\end{proof}

\begin{remark}
    The deformation
     \[\cP_{X^M}^{\pi}[\![\hbar]\!]\colon = \int_M \cA_X^{\pi}[\![\hbar]\!] \in \Alg_{\EE_{n+1-d}}\bigl(\Mod^{\epsilon, \gr}_{\cD_{X^M_{\DR}}[\![\hbar]\!]}\bigr)\]
     can be thought of as a version of the AKSZ construction: given a shifted Poisson structure on the target $X$, we produce a deformation of the mapping stack $X^M$. \cite{Tom} produces shifted Poisson structures on mapping stacks using different methods. We conjecture that the semi-classical information associated to these deformations is compatible with the shifted Poisson structures produced in \cite{Tom}.
~\hfill$\Diamond$\end{remark}

\begin{conj}
     Let $X$ be a derived Artin stack, locally of finite presentation over $\Bbbk$. Let $M$ be a $d$-manifold. Given a shifted Poisson structure $\pi \in {\sf Poiss}(X, n)$, and a framing of $M\times \RR^{n+1-d}$, the deformation
     \[\cP_{X^M}^{\pi}[\![\hbar]\!]\colon = \int_M \cA_X^{\pi}[\![\hbar]\!] \in \Alg_{\EE_{n+1-d}}\bigl(\Mod^{\epsilon, \gr}_{\cD_{X^M_{\DR}}[\![\hbar]\!]}\bigr)\] 
     is a deformation quantization of the $(n-d)$-shifted Poisson structure on $X^M$ produced in \cite{Tom}. 
\end{conj}

We also immediately obtain deformations of various other objects associated to the mapping stack $X^M$. 
\begin{cor}
Let $X$ be an Artin derived stack locally of finite presentation over $\Bbbk$. 
Given an $n$-shifted Poisson structure $\pi \in {\sf Poiss}(X, n)$ with $n\geq1$, a $d$-manifold $M$ where $d \leq n+1$, and a framing of $M\times\RR^{n+1-d}$, we obtain 
\begin{itemize}
    \item an $\EE_{n+1-d}$-deformation of the algebra $\Gamma(X^M, \cO_{X^M})$ given by 
\[Z^\hbar_\pi(M):= \Gamma\Bigl(X^M_{\DR}, \vert\int_M\cA_X^{\hbar, \pi}\vert \Bigr)~.\]
\item an $\EE_{n+1-d}$-deformation of the site of formally \'etale derived stacks over $X^M$ given by 
\[
\begin{array}{ccccc}
\cO_{X^M}^{\pi, {\sf f\acute{e}t}}[\![\hbar]\!]:& X_{\sf f\acute{e}t}& \ra &\Alg_{\EE_{n+1}}(\Mod_{\Bbbk[\![\hbar]\!]})\\
&(Y \xra{f} X^M) &\mapsto &\Gamma(Y_{\DR}, \vert\cP^{\pi}_Y[\![\hbar]\!]\vert).
\end{array}
\]
\end{itemize}
\end{cor}
\begin{proof}
    This is immediate from Proposition \ref{defpp} and Lemmas \ref{quantizingglobalfunctions}, \ref{quantizingthestructuresheaf}. 
\end{proof}

If we have $d \leq n$, then $\Perf^{\epsilon, \gr}_{\cP^{\pi}_{X^M}}$ has the structure of an $\EE_{n-d}$ category (see Definition~\ref{quantumperfect}).  
 
\begin{cor}
Let $X$ be an Artin derived stack locally of finite presentation over $\Bbbk$. 
Given an $n$-shifted Poisson structure $\pi \in {\sf Poiss}(X, n)$ with $n\geq1$, a $d$-manifold $M$ where $d \leq n$, and a framing of $M\times\RR^{n+1-d}$, we obtain an $\EE_{n-d}$-category of quantum perfect complexes on the mapping stack 
\[
\Perf^{\pi}(X^M) := \Perf^{\epsilon, \gr}_{\cP^{\pi}_{X^M}[\![\hbar]\!]}
\]
\end{cor}

The procedure described above provides quantizations of several moduli spaces of interest. We list a few examples.

\begin{example}[Chern--Simons Theory]
    Let $\GG$ be an affine algebraic group over $\Bbbk$, and $M$ a $d$-manifold, with $d \leq 3$. Given a framing on $M\times\RR^{3-d}$, and a ``level" $\lambda\in \Sym^2(\gg)^{\GG} \simeq {\sf Poiss}(\sB\GG, 2)$, we obtain an $\EE_{3-d}$-quantization of the derived stack of $\GG$-local systems on $M$:
    \[\Map_{\dSt}(\Pi_{\infty}M, \sB\GG) \simeq \Loc_{\GG}(M)~.\] We will study this example in detail in the following section. 
    \begin{itemize}
        \item For $M = S^1$, we obtain an $\EE_2$-quantization of 
        \[ \Loc_{\GG}(S^1) \simeq [\GG/\GG]\]
        the adjoint quotient stack of $\GG$. 
        \item For $M$ a closed oriented surface $\Sigma$ (note that $\Sigma\times \RR$ always admits a framing), we obtain an $\EE_1$-quantization of local systems on a surface. These quantum character varieties have been extensively studied in the literature, and it would be interesting to compare the various quantizations produced. See \cite{Jordan} and \cite{BBJ}. 
        \item For $M$ a framed $3$-manifold, (note that every closed oriented $3$-manifold admits a framing), we obtain an $\EE_0$-quantization of $\Loc_{\GG}(M)$. This should be compared with the BV-BRST style quantizations of Chern-Simons theory. Our methods should glue the perturbative quantizations at each local system into a global object with a flat connection over $\Loc_{\GG}(M)$. 
    \end{itemize}
\end{example}

\begin{example}[BF Theory]
   Let $\GG$ be an affine algebraic group over $\Bbbk$, and consider the shifted cotangent bundle $T^{\ast}[n]\sB\GG$. The canonical shifted symplectic structure on $T^*[n]\sB\GG$ equips it with an $n$-shifted Poisson structure. Therefore, given a $d$-manifold $M$, with $d \leq n+1$, and a framing on $M\times\RR^{n+1-d}$, we obtain an $\EE_{n+1-d}$-quantization of the mapping stack $\Map_{\dSt}(\Pi_{\infty}M, T^{\ast}[n]\sB\GG)$. If $M$ is closed and oriented, this mapping stack identifies with a shifted cotangent bundle of the derived stack of $\GG$-local systems on $M$, i.e., we have 
   \[\Map_{\dSt}(\Pi_{\infty}M, T^{\ast}[n]\sB\GG) \simeq T^*[n-d]\Loc_{\GG}(M)~.\] 
   This stack is already interesting in low dimensions: for $n=1$ and $M = S^1$, we obtain an $\EE_1$-quantization of
       \[\Map_{\dSt}(\Pi_{\infty}S^1, T^{\ast}[1]\sB\GG) \simeq T^*[\GG/\GG]\]
    the cotangent bundle of the adjoint quotient $[\GG/\GG]$. 

\end{example}

\begin{example}[Rozansky--Witten Theory]
    Let $Y$ be a smooth complex algebraic variety, and consider the shifted cotangent bundle $T^{\ast}[2]Y$. The canonical symplectic structure on the shifted cotangent bundle gives $T^{\ast}[2]Y$ a $2$-shifted Poisson structure. Again, for a $d$-manifold $M$ with $d \leq 3$, we obtain quantizations of $\Map_{\dSt}(\Pi_{\infty}M, T^*[2]Y)$. For $M = S^1$, we obtain an $\EE_2$-quantization of $T^{\ast}\cL Y$; the shifted cotangent bundle of the free loop stack on $Y$.
\end{example}

\section{Quantizing $\sB\GG$ and $\Loc_{\GG}(M)$} \label{Quantizing BG and Loc}

The results of this section are similar to those of \cite{CPTVV} Section~3.6. Let $\GG$ be an affine algebraic group over $\Bbbk$, 
i.e., the affine group scheme of a finite type commutative Hopf  algebra over $\Bbbk$ 
(here ordinary or underived), 
and let $\gg$ denote its Lie algebra. 
In this section we will produce graded mixed $\EE_3$-deformations of the presheaf of principal parts $\cP_{\sB\GG}$ using methods different from those presented in the previous section. 

Recall that $\cP_{\sB\GG}$ is a prestack of graded mixed algebras over $\sB\GG_{\DR}$. We have a fiber sequence of derived stacks
\[\sB\gg \ra \sB\GG \ra \sB\GG_{\DR}~.\]
The action of $\GG_{\DR}$ on $\sB\gg$ makes $\DD(\sB\gg)$ a comodule over the Hopf algebra $\DD(\GG_{\DR})$.

Consider the universal $\GG_{\DR}$-bundle
\[\GG_{\DR} \ra \pt \ra \sB\GG_{\DR}~.\]
We will define a symmetric monoidal functor (Definition \ref{assocvbundle})
\[{\sf Assoc}\colon {\sf coMod}_{\DD(\GG_{\DR})}(\Mod^{\epsilon, \gr}_{\Bbbk}) \xrightarrow[]{} \Mod^{\epsilon, \gr}_{\cD_{\sB\GG_{\DR}}}\]
that takes a graded mixed $\DD(\GG_{\DR})$-comodule to the ``associated vector bundle" given by twisting the universal $\GG_{\DR}$-bundle by the comodule.  Lemma \ref{ppascomod} will show that there is an equivalence of $\cD_{\sB\GG_{\DR}}$-modules
\[{\sf Assoc}(\DD(\sB\gg)) \simeq \cP_{\sB\GG}~.\]
There is an identification of $\DD(\GG_{\DR})$-comodules 
\[\DD(\sB\gg) \simeq \sC^{\epsilon}_{\Lie}(\gg)~,\]
where $\sC^{\epsilon}_{\Lie}(\gg)$ is the graded mixed Chevalley--Eilenberg complex (See \cite{pav2} for definitions of graded mixed Chevalley--Eilenberg complexes.) 
Given a level $\lambda \in \Sym^2(\gg)^{\GG}$, 
we will produce an $\EE_3$-deformation $\sC^{\epsilon}_{\lambda}(\gg)[\![\hbar]\!]$ of $\sC^{\epsilon}_{\Lie}(\gg)$ as a $\DD(\GG_{\DR})[\![\hbar]\!]$-comodule (Lemma \ref{equivariance}). 
Consequently, we obtain a $\cD_{\sB\GG_{\DR}}[\![\hbar]\!]$-linear $\EE_3$-deformation $\cP^{\lambda}_{\sB\GG}[\![\hbar]\!]$ of the algebra of principal parts of $\sB\GG$. 
\begin{remark}
    The $\cD_{\sB\GG_{\DR}}[\![\hbar]\!]$-linear $\EE_3$-deformation $\cP^{\lambda}_{\sB\GG}[\![\hbar]\!]$ of the algebra of principal parts of $\sB\GG$ will be obtained by constructing a strict graded mixed $\PP_3$-algebra $\bigl(\sC^{\epsilon}_{\Lie}(\gg), \{-,-\}_{\lambda}\bigr)$ where the bracket is given by 
    \[\{-,-\}_{\lambda} = \hbar\cdot \lambda(-,-)\]
    on generators $\gg^{\vee}[-1](1) \subset \sC^{\epsilon}_{\Lie}(\gg)$. We will then use formality to obtain a $\EE_3$-algebra. Note that the level $\lambda$ only sits in one degree in $\hbar$, i.e., $\lambda \in \Sym^2(\gg)^{\GG} \simeq \pi_0{\sf Poiss}(\sB\GG, 2)$. We do not construct quantizations for general $\lambda\in \Sym^2(\gg)^{\GG}[\![\hbar]\!]$. 
~\hfill$\Diamond$\end{remark}
Further, given a $d$-manifold $M$ with $d\leq 3$, and a framing of $M\times\RR^{3-d}$,  we obtain a $\Disk_M$-algebra $\cA^{\lambda}_{\sB\GG}[\![\hbar]\!]$ quantizing the algebra of local principal parts $\cA_{\sB\GG}$ on the mapping stack $\sB\GG^M$. 
The algebra $\cA^{\lambda}_{\sB\GG}[\![\hbar]\!]$ can be thought of as the algebra of local observables for quantum Chern--Simons theory on the manifold~$M$. 
Taking factorization homology of $\cA^{\lambda}_{\sB\GG}[\![\hbar]\!]$ over $M$, we obtain a quantization of the derived stack of $\GG$-local systems on $M$ (Proposition~\ref{quantchernsimons}).

\subsection{A Künneth Formula}
We will make use of the following Künneth formula for the functor $\DD$ in the next few subsections. 

Recall the notion of an algebraisable affine formal derived stack from \cite{CPTVV} Definition 2.2.1. We will show that our examples of interest are algebraisable affine formal derived stacks, in order to obtain Künneth formulae for them. The reader can skip the discussion in the example below if they are willing to accept that the hypotheses of Lemma~\ref{kuenneth} apply to our examples. 

\begin{example}
    Let $\GG$ be an affine algebraic group over $\Bbbk$, and $\gg$ its Lie algebra.
    The following derived stacks of interest are algebraisable affine formal derived stacks: 
    \begin{itemize}
        \item $\GG_{\DR}$: it is the formal completion of $\pt$ along $(\GG_{\DR})_{\sf red} \simeq \GG \ra \pt$. Note that $\GG$ is reduced as $\Bbbk$ is of characteristic zero. 
        \item $\sB\gg$: it is the formal completion of $\sB\GG$ along $(\sB\gg)_{\sf red}\simeq \pt \mapsto \sB\GG$. 
        \item Take a map $f\colon \Spec(R) \ra \sB\GG_{\DR}$, classifying a $\GG_{\DR}$-torsor $P \ra \Spec(R)$. We have the corresponding map $f\colon \Spec(R_{\sf red}) \ra \sB\GG$ classfying a $\GG$-torsor $Q \ra \Spec(R_{\sf red})$. Observe that $P_{\sf red} \simeq Q$ and that $Q$ is affine. $P$ is an algebraizable affine formal derived stack: it is the formal completion of $\Spec(R)$ along  $Q \ra \Spec(R_{\sf red}) \ra \Spec(R)$. 
    \end{itemize}
\end{example}

\begin{lemma} \label{kuenneth}
   Assume we have a pullback square of algebraisable affine formal derived stacks 
   \[\begin{tikzcd}
	W & X \\
	Y & Z
	\arrow[from=1-1, to=1-2]
	\arrow[from=1-1, to=2-1]
	\arrow[from=1-2, to=2-2]
	\arrow[from=2-1, to=2-2]
\end{tikzcd}\]
    Assume that the natural map 
    \[W_{\sf red} \ra X_{\sf red} \underset{Z_{\sf red}}{\times}Y_{\sf red}\]
    is an equivalence. 
    Then the natural map of graded mixed algebras
    \[\DD(X)\otimes_{\DD(Z)}\DD(Y) \ra \DD(W)\]
    is an equivalence. 
\end{lemma}
\begin{proof}
    It suffices to check that the map is an equivalence of graded mixed modules. Further it suffices to check that the map is an equivalence after forgetting the mixed structures. As $X,Y,Z$ and $W$ are algebraisable affine formal derived stacks,  by \cite{CPTVV} Theorem 2.2.2 it suffices to show that 
    \[\Sym_{\cO(X_{\sf red})}(\LL_{X_{\sf red}\vert X}[-1](1))\underset{\Sym_{\cO(Z_{\sf red})}(\LL_{Z_{\sf red}\vert Z}[-1](1))}{\otimes} \Sym_{\cO(Y_{\sf red})}(\LL_{Y_{\sf red}\vert Y}[-1](1)) \ra \Sym_{\cO(W_{\sf red})}(\LL_{W_{\sf red}\vert W}[-1](1))\]
    is an equivalence. Note that all of the reduced stacks above are affine by the definition of an algebraisable affine formal derived stack. However, we have a pushout diagram in $\Mod_{\cO(W_{\sf red})}$
    \[\begin{tikzcd}
	{\cO(W_{\sf red})\otimes_{\cO(Z_{\sf red})}\LL_{Z_{\sf red}\vert Z}} & {\cO(W_{\sf red})\otimes_{\cO(Y_{\sf red})}\LL_{Y_{\sf red}\vert Y}} \\
	{\cO(W_{\sf red})\otimes_{\cO(X_{\sf red})}\LL_{X_{\sf red}\vert X}} & {\LL_{W_{\sf red}\vert W}}
	\arrow[from=1-1, to=1-2]
	\arrow[from=1-1, to=2-1]
	\arrow[from=1-2, to=2-2]
	\arrow[from=2-1, to=2-2]
\end{tikzcd}\]
As $\Sym_{\cO(W_{\sf red})}$ is a left adjoint, it preserves pushouts, and we have the result. 
\end{proof}
\subsection{The Associated Bundle Construction}
The group operation on $\GG_{\DR}$ induces the structure of a graded mixed Hopf algebra on $\DD(\GG_{\DR})$. Here $\DD(\GG_{\DR})$ should be thought of as a graded mixed version of the de Rham complex of $\GG$.

Consider the category of graded mixed comodules ${\sf coMod}_{\DD(\sB\GG)}$. We will construct a symmetric monoidal functor 
\[{\sf Assoc}\colon {\sf coMod}^{\epsilon, \gr}_{\DD(\GG_{\DR})} \ra \Mod^{\epsilon, \gr}_{\cD_{\sB\GG_{\DR}}}\]
which roughly takes a $\GG_{\DR}$-representation to its associated bundle over $\sB\GG_{\DR}$. 

Consider the adjunction 
\[{\sf triv}\colon \Mod_{\Bbbk}^{\epsilon, \gr} \rightleftarrows {\sf coMod}^{\epsilon, \gr}_{\DD(\GG_{\DR})} \colon \Bbbk\square(-)\]
where, ${\sf triv}$ denotes the trivial comodule functor, and $\Bbbk\square(-)$ denotes the cotensor product with the trivial comodule $\Bbbk$. Note that ${\sf triv}$ has the structure of a symmetric monoidal functor. For $R \in \CAlg^{\geq 0}$, this induces an adjunction 
\[{\sf triv}\colon \Mod_{\DD(R)}^{\epsilon, \gr} \rightleftarrows \Mod^{\epsilon, \gr}_{\DD(R)}({\sf coMod}^{\epsilon, \gr}_{\DD(\GG_{\DR})}) \colon \Bbbk\square(-)~.\]

Consider the universal $\GG_{\DR}$-bundle 
\[\GG_{\DR} \ra \pt \ra \sB\GG_{\DR}~.\]

Given a map $f\colon \Spec(R) \ra \sB\GG_{\DR}$, consider the corresponding $\GG_{\DR}$ torsor $P \ra \Spec(R)$. We obtain a map of graded mixed algebras $\DD(R) \ra \DD(P)$. The action of $\GG_{\DR}$ on $P$ makes $\DD(P)$ a graded mixed $\DD(\GG_{\DR})$-comodule. More precisely, it promotes $\DD(P)$ to an object of the category $\CAlg({\sf coMod}^{\epsilon, \gr}_{\DD(\GG_{\DR})})$. Give $\DD(R)$ the trivial $\DD(\GG_{\DR})$-comodule structure. The map of $\DD(\GG_{\DR})$-comodule algebras
\[\DD(R) \ra \DD(P)\]
is a Hopf--Galois context in the sense of \cite{Ergus2022Hopf}. Extension of scalars $\DD(P)\otimes-$ gives us an adjunction 
\[\Mod_{\DD(R)}^{\epsilon, \gr} \rightleftarrows \Mod^{\epsilon, \gr}_{\DD(R)}({\sf coMod}^{\epsilon, \gr}_{\DD(\GG_{\DR})}) \rightleftarrows \Mod^{\epsilon, \gr}_{\DD(P)}({\sf coMod}^{\epsilon, \gr}_{\DD(\GG_{\DR})})~.\]
The left adjoint 
\[\Mod_{\DD(R)}^{\epsilon, \gr} \xra{\sf triv} \Mod^{\epsilon, \gr}_{\DD(R)}({\sf coMod}^{\epsilon, \gr}_{\DD(\GG_{\DR})}) \xra{\DD(P)\otimes_{\DD(R)} -} \Mod^{\epsilon, \gr}_{\DD(P)}({\sf coMod}^{\epsilon, \gr}_{\DD(\GG_{\DR})})\]
is given by 
\[M \mapsto\DD(P)\otimes_{\DD(R)}M\]
where the coaction of $\DD(\GG_{\DR})$ is trivial. 
The right adjoint 
\[\Mod^{\epsilon, \gr}_{\DD(P)}({\sf coMod}^{\epsilon, \gr}_{\DD(\GG_{\DR})}) \xra{\sf fgt} \Mod^{\epsilon, \gr}_{\DD(R)}({\sf coMod}^{\epsilon, \gr}_{\DD(\GG_{\DR})}) \xra{\Bbbk\square(-)}\Mod_{\DD(R)}^{\epsilon, \gr}~\]
is given by 
\[N \mapsto \Bbbk\square N\]
where $\Bbbk\square N$ denotes the cotensor product of $\DD(\GG_{\DR})$-comodules and we remember the $\DD(R)$-action. Note that the left adjoint is symmetric monoidal. 

Given a Hopf--Galois context 
\[\DD(R) \ra \DD(P)\]
there is a shear map \cite{Ergus2022Hopf} Definition~3.1.7
\[{\sf shear} \colon \DD(P) \otimes_{\DD(R)}\DD(P) \ra \DD(\GG_{\DR})\otimes \DD(P)\otimes_{\DD(A)} \DD(P) \ra \DD(\GG_{\DR})\otimes\DD(P)~,\]
where the first map uses the coaction on $\DD(P)$ and the second map uses the algebra structure on $\DD(P)$. 
\begin{lemma} \label{shear}
    Given a map $\Spec(R)\rightarrow \sB\GG_{\DR}$ classifying a torsor $P \ra \Spec(R)$, the canonical shear map 
    \[{\sf shear} \colon \DD(P) \otimes_{\DD(R)}\DD(P) \ra \DD(\GG_{\DR})\otimes \DD(P)\otimes_{\DD(A)} \DD(P) \ra \DD(\GG_{\DR})\otimes\DD(P)~\]
    is an equivalence. 
\end{lemma}
\begin{proof}
     We have 
\[\GG_{\DR} \times P \xra{\simeq} P\underset{\Spec(R)}{\times}P\] and by \ref{kuenneth}, the following diagram commutes 
\[\begin{tikzcd}
	{\DD(P)\otimes_{\DD(R)}\DD(P)} & {\DD(\GG_{\DR})\otimes\DD(P)} \\
	{\DD\Bigl(P\underset{\Spec(R)}{\times}P\Bigr)} & {\DD(\GG_{\DR})\otimes\DD(P)}
	\arrow["{{\sf shear}}", from=1-1, to=1-2]
	\arrow["\simeq", from=1-1, to=2-1]
	\arrow["\simeq", from=1-2, to=2-2]
	\arrow["\simeq"', from=2-1, to=2-2]
\end{tikzcd}\]
Therefore, we have the result. 
\end{proof}

We have the following theorem of \cite{Ergus2022Hopf}, specialized to our context. 

\begin{theorem}[\cite{Ergus2022Hopf} Corollary 3.2.11]\label{comonad}
    Given a map $\Spec(R)\rightarrow \sB\GG_{\DR}$ classifying a torsor $P \ra \Spec(R)$, if the map  
    \[{\sf shear}\colon \DD(P) \otimes_{\DD(R)}\DD(P)\ra \DD(\GG_{\DR}) \otimes \DD(P)  \]
    is an equivalence, then the following are equivalent
    \begin{itemize}
        \item the functors of the adjunction 
    \[\Mod_{\DD(R)}^{\epsilon, \gr} \rightleftarrows \Mod^{\epsilon, \gr}_{\DD(P)}({\sf coMod}_{\DD(\GG_{\DR})})\]
    are equivalences of symmetric monoidal categories 
        \item the adjunction given by extension and restriction of scalars 
        \[\DD(P)\otimes_{\DD(R)}-\colon\Mod^{\epsilon, \gr}_{\DD(R)} \rightleftarrows \Mod^{\epsilon, \gr}_{\DD(P)} \colon{\sf fgt}\]
        is comonadic. 
    \end{itemize}
    
\end{theorem}

Now we will show that the adjunction 
\[\DD(P)\otimes_{\DD(R)}-\colon\Mod^{\epsilon, \gr}_{\DD(R)} \rightleftarrows \Mod^{\epsilon, \gr}_{\DD(P)} \colon{\sf fgt}\]
is comonadic. By \cite{Mathew2014TheGG} Proposition 3.22., it suffices to show that the map $\DD(R) \ra \DD(P)$ admits descent in the sense of \cite{Mathew2014TheGG} Definition~3.18.

We will need the following lemma. 
\begin{lemma}\label{descentforetalecovers}
Given an \'etale cover 
\[\Spec(R') \ra \Spec(R)\]
the map 
\[\DD(R) \ra \DD(R')\]
admits descent in the sense of \cite{Mathew2014TheGG} Definition~3.18.
\end{lemma}
\begin{proof}
    By \cite{Mathew2014TheGG} Corollary 3.33, the map $R \ra R'$ admits descent. Consider the symmetric monoidal functor
    \[\DD(R) \otimes_R- \colon \Mod_R \ra \Mod^{\epsilon, \gr}_{\DD(R)}.\]
    By \cite{Mathew2014TheGG} Corollary~3.21, the map $\DD(R) \ra \DD(R)\otimes_RR'$ admits descent as $R \ra R'$ does. Now observe that the natural map 
    \[\DD(R)\otimes_RR' \ra \DD(R')\]
    is an equivalence. Indeed, we have 
    \[R'_{\sf red} \simeq R_{\sf red}\otimes R'\]
    and by \cite{CPTVV}~Lemma 1.3.18
    \[\begin{split}
        \DD(R') &\simeq \DR(R'_{\sf red}/R') \\
        &\simeq \DR(R_{\sf red}\underset{R}{\otimes} R'/R') \\
        &\simeq \DR(R_{\sf red}\underset{R}{\otimes} R')\underset{\DR(R')}{\otimes} R' \\
        &\simeq \DR(R_{\sf red})\underset{\DR(R)}{\otimes} R' \\
        &\simeq \DD(R_{\sf red}/ R)\otimes_R R'~.
    \end{split}\]
\end{proof}

\begin{lemma}\label{descentalongatorsor}
   Given a map $f \colon \Spec(R) \ra \sB\GG_{\DR}$, classifing a torsor $P_f := \Spec(R) \times_{\sB\GG_{\DR}}\pt$, the map 
   \[\DD(R) \ra \DD(P)\]
   admits descent. 
\end{lemma}
\begin{proof}
    The $\GG_{\DR}$-torsor classified by $f\colon \Spec(R) \ra \sB\GG_{\DR}$ trivializes over an \'etale cover $\Spec(R')\ra \Spec(R)$, i.e. we have a section 
    \[\Spec(R') \ra P~.\]
    Consider the composite 
    \[\DD(R) \ra \DD(P) \ra \DD(R').\]
    The composite admits descent by Lemma~\ref{descentforetalecovers}. Therefore, by \cite{Mathew2014TheGG} Proposition~3.24, $\DD(R) \ra \DD(\Spec(R) \times_{\sB\GG_{\DR}}\pt)$ admits descent. 
\end{proof}
Therefore, by Theorem \ref{comonad}, we immediately have the following corollary. 
\begin{cor}
    Given a map $\Spec(R) \ra \sB\GG_{\DR}$ classifying a torsor $P$, the functor 
    \[ \Mod^{\epsilon, \gr}_{\DD(P)}({\sf coMod}^{\epsilon, \gr}_{\DD(\GG_{\DR})}) \ra \Mod_{\DD(R)}^{\epsilon, \gr}\]
    given by 
    \[M \mapsto \Bbbk\square M\]
    is an equivalence of symmetric monoidal categories. 
\end{cor}
\begin{definition}
Given a map $f\colon\Spec(R) \ra \sB\GG_{\DR}$ classifying a torsor $P$
    the functor
\[
{\sf Assoc}_f: {\sf coMod}_{\sC^{\epsilon}_{\Lie}(\gg, \cO(\GG))}(\Mod^{\epsilon, \gr}_{\Bbbk}) \xrightarrow[]{} \Mod^{\epsilon, \gr}_{\DD(R)}
\] 
is the composite of symmetric monoidal functors 
\[{\sf coMod}^{\epsilon, \gr}_{\DD(\sB\GG_{\DR})} \xra{\DD(P)\otimes-} \Mod^{\epsilon, \gr}_{\DD(P)}({\sf coMod}^{\epsilon, \gr}_{\DD(\GG_{\DR})}) \xra{\Bbbk\square(-)} \Mod^{\epsilon, \gr}_{\DD(R)}~.\]
\end{definition}

The construction above is natural in the map $\Spec(R) \ra \sB\GG_{\DR}$. Here's a sketch: given 
\[\Spec(S)\xra{g}\Spec(R) \xra{f} \sB\GG_{\DR}\] 
classifying torsors $P \ra \Spec(R)$ and $Q \ra \Spec(S)$, the following diagram commutes 
\[\begin{tikzcd}
	{} & {\Mod^{\epsilon, \gr}_{\DD(Q)}({\sf coMod}^{\epsilon, \gr}_{\DD(\GG_{\DR})})} & {\Mod^{\epsilon, \gr}_{\DD(S)}} \\
	& {\Mod^{\epsilon, \gr}_{\DD(P)}({\sf coMod}^{\epsilon, \gr}_{\DD(\GG_{\DR})})} & {\Mod^{\epsilon, \gr}_{\DD(R)}}
	\arrow["{\Bbbk\square-}", from=1-2, to=1-3]
	\arrow["{{\sf fgt}}", from=1-2, to=2-2]
	\arrow["{{\sf fgt}}", from=1-3, to=2-3]
	\arrow["{\Bbbk \square -}"', from=2-2, to=2-3]
\end{tikzcd}\]
Therefore, we obtain a natural transformation of functors 
\[{\sf Assoc}_f \ra {\sf Assoc}_{gf}~.\]

\begin{definition} \label{assocvbundle}
    The symmetric monoidal functor 
    \[{\sf Assoc}: {\sf coMod}_{\DD(\GG_{\DR})}(\Mod^{\epsilon, \gr}_{\Bbbk}) \xrightarrow[]{} \Mod^{\epsilon, \gr}_{\cD_{\sB\GG_{\DR}}}
    \] 
    takes an object $V\in {\sf coMod}_{\DD(\GG_{\DR})}(\Mod^{\epsilon, \gr}_{\Bbbk})$ to the object ${\sf Assoc}(V) \in \Mod^{\epsilon, \gr}_{\cD_{\sB\GG_{\DR}}}$ given by 
    \[{\sf Assoc}(V): ({\sf Aff}_{/\sB\GG_{\DR}})^{\op} \ra \Mod_{\Bbbk}^{\epsilon, \gr}\]
    \[(\Spec(R) \xra{f}\sB\GG_{\DR}) \mapsto {\sf Assoc}_f(V)~.\]
\end{definition}
As {\sf Assoc} is symmetric monoidal, we have
\[{\sf Assoc}(\Bbbk) \simeq \cD_{\sB\GG_{\DR}}~.\] 

Now consider the fiber sequence of derived stacks
\[\sB\gg \ra \sB\GG \ra \sB\GG_{\DR}~.\]
The action of $\GG_{\DR}$ on $\sB\gg$ makes $\DD(\sB\gg)$ a graded mixed $\DD(\GG_{\DR})$-comodule. The $\cD_{\sB\GG_{\DR}}$-algebra associated to $\DD(\sB\gg)$ is the algebra of principal parts on $\cP_{\sB\GG}$. 

\begin{lemma} \label{ppascomod}
    There is an equivalence of $\cD_{\sB\GG_{\DR}}$-algebras
    \[{\sf Assoc}(\DD(\sB\gg)) \simeq \cP_{\sB\GG} \in \Mod^{\epsilon, \gr}_{\cD_{\sB\GG_{\DR}}}~.\]
\end{lemma}
\begin{proof}
    We will construct a map of $\cD_{\sB\GG_{\DR}}$-algebras 
    \[\cP_{\sB\GG} \ra {\sf Assoc}(\DD(\sB\gg))~.\]
    Take a map $f \colon \Spec(R) \ra \sB\GG_{\DR}$, classifying a torsor $P \ra \Spec(R)$. We also have $\sB\GG_R = \Spec(R)\times_{\sB\GG_{\DR}}\sB\GG$. Pulling back $\sB\GG_R \ra \sB\GG_{\DR}$ along the universal torsor $\pt \ra \sB\GG_{\DR}$, we have a diagram where all of the squares are pullback diagrams
    \[\begin{tikzcd}
	{P\times \sB\gg} & P & \pt \\
	{\sB\GG_R} & {\Spec(R)} & {\sB\GG_{\DR}}
	\arrow[from=1-1, to=1-2]
	\arrow[from=1-1, to=2-1]
	\arrow[from=1-2, to=1-3]
	\arrow[from=1-2, to=2-2]
	\arrow[from=1-3, to=2-3]
	\arrow[from=2-1, to=2-2]
	\arrow["f"', from=2-2, to=2-3]
\end{tikzcd}\]
    Indeed, we have
    \[
    \begin{split}
        \sB\GG_R \times_{\sB\GG_{\DR}} \pt &\simeq (\Spec(R)\times_{\sB\GG_{\DR}}\sB\GG)\times_{\sB\GG_{\DR}}\pt \\
        &\simeq (\Spec(R) \times_{\sB\GG_{\DR}}\pt)\times (\sB\GG \times_{\sB\GG_{\DR}}\pt
        \simeq P \times \sB\gg~.
    \end{split}
    \]
    Moreover, this equivalence is $\GG_{\DR}$-equivariant. Therefore we have a $\GG_{\DR}$-torsor 
    \[P \times \sB\gg \ra \sB\GG_{R}\]
    over $\Spec(R)$. Therefore, we have a map of $\DD(\GG_{\DR})$-comodules (where the action on $\DD(\sB\GG_R)$ is trivial)
    \[\cP_{\sB\GG}(f) \simeq\DD(\sB\GG_R) \ra \DD(P\times \sB\gg) \simeq \DD(P)\otimes\DD(\sB\gg)~.\]
    Therefore by adjunction, we have a map of $\DD(R)$-algebras 
    \[\cP_{\sB\GG}(f) \ra \Bbbk\square(\DD(P)\otimes\DD(\sB\gg)) \simeq {\sf Assoc}(\DD(\sB\gg)(f)\]
    This construction is natural in the map $f\colon\Spec(R) \ra \sB\GG_{\DR}$, and we have a map of $\cD_{\sB\GG_{\DR}}$-algebras
    \[\cP_{\sB\GG} \ra {\sf Assoc}(\DD(\sB\gg))~.\]
    This map is an equivalence. Indeed, we have 
    \[\sB\GG_R \simeq \colim(\dots(P\times \sB\gg)\times\GG_{\DR}\rightrightarrows (P\times \sB\gg))~,\]
    and as $\DD$ takes colimits of derived stacks to limits, we have 
    \[\DD(\sB\GG_{R}) \simeq \lim((\DD(P)\otimes\DD(\sB\gg)) \rightrightarrows (\DD(P)\otimes\DD(\sB\gg))\otimes\DD(\GG_{\DR}) \dots)~.\]
    However, note that the limit above computes $\Bbbk\square(\DD(P)\otimes\DD(\sB\gg))$, and we have the result. 
\end{proof}
\subsection{An explicit presentation of the Hopf Algebra $\DD(\GG_{\DR})$} 
By \cite{saf} Lemma 1.14, identifying $\gg$ with right-invariant vector fields, we have an identification of graded mixed algebras (see \cite{pav2} for definitions of graded mixed Chevalley--Eilenberg complexes)
\[
\DD(\GG_{\DR}) \simeq \sC^{\epsilon}_{\Lie}(\gg, \cO(\GG))
\]
where $\gg$ acts on $\cO(\GG)$ by infinitesimal left translations. 
The algebra $\sC^{\epsilon}_{\Lie}(\gg, \cO(\GG))$ obtains a Hopf algebra structure  transferred from $\DD(\GG_{\DR})$. \cite{saf} Section 1.5 describes the Hopf algebra structure under the identification above. We include it here for convenience.  All tensor products are over $\Bbbk$ unless indicated otherwise. 

The coproduct on $\sC^{\epsilon}_{\Lie}(\gg, \cO(\GG))$ is determined by the following properties. 
\begin{itemize}
    \item The diagram 
    \[
    \begin{tikzcd}
        \cO(\GG) \arrow[r, "\Delta"] \arrow[d] & \arrow[d] \cO(\GG) \otimes \cO(\GG)\\
        {\sf fgt}_{\epsilon}\bigl(\sC^{\epsilon}_{\Lie}(\gg, \cO(\GG))\bigr) \arrow[r, "\Delta"] &{\sf fgt}_{\epsilon}\bigl(\sC^{\epsilon}_{\Lie}(\gg, \cO(\GG))\bigr)\otimes {\sf fgt}_{\epsilon}\bigl(\sC^{\epsilon}_{\Lie}(\gg, \cO(\GG))\bigr)
    \end{tikzcd}
    \]
    commutes. 
    \item For $\xi \in \gg^{\vee}[-1] \subset \sC^{\epsilon}_{\Lie}(\gg, \cO(\GG))$, we have 
    \[\Delta(\xi) = \xi\otimes 1 + {\sf coAd}_{g_{(1)}}(1\otimes \xi)~,\]
    where $g_{(1)}$ is the coordinate of the first factor of $\cO(\GG)$ on $\sC^{\epsilon}_{\Lie}(\gg, \cO(\GG))\otimes \sC^{\epsilon}_{\Lie}(\gg, \cO(\GG))$.
\end{itemize}

\begin{remark}
    As (non-mixed) graded algebras, one could see the identification 
    \[\Sym_{\cO(\GG)}\bigl(\LL_{\GG/\Bbbk}[-1](1)\bigr)\simeq {\sf fgt}_{\epsilon}\bigl(\DD(\GG_{\DR})\bigr) \simeq \bigl(\sC^{\epsilon}_{\Lie}(\gg, \cO(\GG))\bigr) \simeq \Sym_{\Bbbk}(\gg^{\vee}[-1](1))\otimes\cO(\GG)\]
    as coming from the trivialization of the tangent bundle provided by the  Maurer--Cartan form
    \[\theta: \TT_{\GG/\Bbbk} \xra{\simeq} \cO(\GG)\otimes\gg~.\]
    \cite{saf} further proves that the equivalence is compatible with the mixed structures. 
    This is analogous to what happens in ordinary Lie theory: given a Lie group $G$ with Lie algebra $\gg$, the Maurer--Cartan form $\theta \in \Omega^1(G, \gg) $ trivializes the tangent bundle 
    \[\theta: TG \xra{\cong} G \times \gg~,\]
    providing the map 
    \[\Omega^{\bullet}(G) \cong C^{\infty}(G)\otimes_{\RR}\bigwedge^{\bullet}\gg^{\vee}\]
    a well-known isomorphism.~\hfill$\Diamond$
\end{remark}

Consider the coadjoint action 
\[{\sf coAd} \colon \GG \ra {\sf GL}(\gg^{\vee})~.\]
At the universal element ${\sf id}_{\cO(\GG)} \in \GG(\cO(\GG))$, we have 
\[{\sf coAd}: \cO(\GG)\otimes\gg^{\vee} \ra \cO(\GG)\otimes \gg^{\vee}\]
Choose a basis $\{\xi^i\}_i$ of $\gg^{\vee}$, and define elements $\alpha^i_j \in \cO(\GG)$ by the formula 
\[
{\sf coAd}(1 \otimes \xi^i) = \sum_j\alpha^i_j\otimes \xi^j.
\]
Then the coproduct 
\[\Delta: \sC^{\epsilon}_{\Lie}(\gg, \cO(\GG)) \rightarrow \sC^{\epsilon}_{\Lie}(\gg, \cO(\GG))\otimes\sC^{\epsilon}_{\Lie}(\gg, \cO(\GG))\]
is given on generators of $\gg^{\vee}[-1]$ by 
\[
\Delta(\xi^i) = \xi^i\otimes 1 + \sum_j \alpha^i_j\otimes \xi^j~.
\]
As $\Delta$ is a map of commutative algebras, given $f \in \cO(\GG)$ and $\xi^{i_1},\dots,\xi^{i_m}$ we have 
\[\Delta(f\xi^{i_1}\dots\xi^{i_m}) = \Delta(f)\Delta(\xi^{i_1})\dots\Delta(\xi^{i_m})\]
where $\Delta(f)$ denotes the coproduct on the coalgebra $\cO(\GG)$. As $\sC^{\epsilon}_{\Lie}(\gg, \cO(\GG))$ is generated by $\gg^{\vee}[-1]$ and $\cO(\GG)$ as a commutative algebra, the formula above computes the coproduct on all of~$\sC^{\epsilon}_{\Lie}(\gg, \cO(\GG))$.
 
\subsection{$\sC^{\epsilon}_{\Lie}(\gg)$ as a $\sC^{\epsilon}_{\Lie}(\gg, \cO(\GG))$-comodule}
By \cite{saf} Theorem 1.15, we have an identification of $\DD(\GG_{\DR})$-comodules
\[\DD(\sB\gg)\simeq \sC^{\epsilon}_{\Lie}(\gg)~.\]
We will describe $\sC^{\epsilon}_{\Lie}(\gg)$ as a $\sC^{\epsilon}_{\Lie}(\gg, \cO(\GG))$-comodule. Note that the quotient map 
\[\GG_{\DR} \rightarrow \GG_{\DR}/\GG \simeq \sB\gg\]
induces a map
\[\sC^{\epsilon}_{\Lie}(\gg) \simeq \DD(\sB\gg) \rightarrow \DD(\GG_{\DR})\simeq \sC^{\epsilon}_{\Lie}(\gg, \cO(\GG)).\]
This can be viewed as the inclusion of $\GG$-invariant differential forms on $\GG$ into differential forms on $\GG$. We will use the same notation to denote elements of $\sC^{\epsilon}_{\Lie}(\gg)$ and elements of $\sC^{\epsilon}_{\Lie}(\gg, \cO(\GG))$. Note that $\GG_{\DR} \rightarrow \GG_{\DR}/\GG \simeq \sB\gg$ is equivariant for the translation action of $\GG_{\DR}$. Therefore $\sC^{\epsilon}_{\Lie}(\gg) \rightarrow \sC^{\epsilon}_{\Lie}(\gg, \cO(\GG))$ is a $\sC^{\epsilon}_{\Lie}(\gg, \cO(\GG))$-comodule map. Therefore the coaction \[\rho: \sC^{\epsilon}_{\Lie}(\gg) \rightarrow \sC^{\epsilon}_{\Lie}(\gg, \cO(\GG))\otimes \sC^{\epsilon}_{\Lie}(\gg)\] is determined by the commutative diagram 
\[\begin{tikzcd}
	{\sC^{\epsilon}_{\Lie}(\gg)} & {} & {\sC^{\epsilon}_{\Lie}(\gg, \cO(\GG))} \\
	{\sC^{\epsilon}_{\Lie}(\gg)\otimes \sC^{\epsilon}_{\Lie}(\gg, \cO(\GG))} && {\sC^{\epsilon}_{\Lie}(\gg, \cO(\GG))\otimes \sC^{\epsilon}_{\Lie}(\gg, \cO(\GG))}
	\arrow["i", from=1-1, to=1-3]
	\arrow["\rho"', from=1-1, to=2-1]
	\arrow["\Delta", from=1-3, to=2-3]
	\arrow["{i\otimes {\sf id}}"', from=2-1, to=2-3]
\end{tikzcd}\]
Let $\{\xi^i\}$ be a basis for $\gg^{\vee}[-1]$. Again define functions $\alpha^i_j \in \cO(\GG)$ by the formula 
\[{\sf coAd}_g(\xi^i) = \sum_j\alpha^i_j(g)\xi^j.\]
Then the commutativity of the above diagram gives us a formula for the coaction
\[\rho: \sC^{\epsilon}_{\Lie}(\gg, \cO(\GG)) \rightarrow \sC^{\epsilon}_{\Lie}(\gg, \cO(\GG))\otimes \sC^{\epsilon}_{\Lie}(\gg)\]
\[\xi^i \mapsto \xi^i\otimes 1 + \sum_j \alpha^i_j\otimes \xi^j.\]
The explicit formula for the coaction will be useful to us in the proof of Lemma~\ref{equivariance}. 
In the next subsection, we will describe some $\EE_3$-deformations of $\sC^{\epsilon}_{\Lie}(\gg)$ as a $\sC^{\epsilon}_{\Lie}(\gg, \cO(\GG))$-comodule. As ${\sf Assoc}$ is symmetric monoidal, by Lemma \ref{ppascomod}, any such deformation will give us a deformation of the presheaf of principal parts.  

\subsection{$\EE_3$-deformations of $\sC^{\epsilon}_{\Lie}(\gg)$}

Let $\hbar$ be a formal variable of weight $2$ and cohomological degree $0$. Let $\sC^{\epsilon}_{\sf Lie}(\gg)[\![\hbar]\!]:= \sC^{\epsilon}_{\sf Lie}(\gg)\otimes \Bbbk[\![\hbar]\!]$ and $\sC^{\epsilon}_{\sf Lie}(\gg, \cO(\GG))[\![\hbar]\!] := \sC^{\epsilon}_{\sf Lie}(\gg, \cO(\GG)) \otimes\Bbbk[\![\hbar]\!]$. 
As extension of scalars 
\[\Mod_{\Bbbk}^{\epsilon, \gr} \xrightarrow[]{-\otimes \Bbbk[\![\hbar]\!]} \Mod_{\Bbbk[\![\hbar]\!]}^{\epsilon, \gr}\]
is symmetric monoidal, $\sC^{\epsilon}_{\sf Lie}(\gg, \cO(\GG))[\![\hbar]\!]$ is a Hopf algebra in $\Mod^{\epsilon, \gr}_{\Bbbk[\![\hbar]\!]}$ and $\sC^{\epsilon}_{\sf Lie}(\gg)[\![\hbar]\!]$ is a $\sC^{\epsilon}_{\sf Lie}(\gg, \cO(\GG))[\![\hbar]\!]$-comodule. In order to produce $\EE_3$-deformations of $\sC^{\epsilon}_{\sf Lie}(\gg)[\![\hbar]\!]$, we will produce $\PP_3$-deformations and then use formality of $\EE_3$. 

\begin{definition}
    Given $\lambda \in \Sym^2(\gg)^{\GG}$, let 
    \[\{-, -\}_{\lambda}: \sC^{\epsilon}_{\Lie}(\gg)[\![\hbar]\!]\underset{\Bbbk[\![\hbar]\!]}\otimes\sC^{\epsilon}_{\Lie}(\gg)[\![\hbar]\!] \rightarrow \sC^{\epsilon}_{\Lie}(\gg)[\![\hbar]\!]\]
    be the graded bi-derivation determined by 
    \[\hbar\cdot \lambda : \gg^{\vee}[-1]\otimes \gg^{\vee}[-1] \rightarrow \Bbbk[\![\hbar]\!].\]
\end{definition}

We have the following lemma. 

\begin{lemma}
Given $\lambda \in \Sym^2(\gg)^{\GG}$, the pair $(\sC^{\epsilon}_{\Lie}(\gg)[\![\hbar]\!], \{-, -\}_{\lambda})$ defines an object in $\Alg_{\PP_3}(\Mod_{\Bbbk[\![\hbar]\!]}^{\epsilon, \gr})$. 
\end{lemma}

\begin{proof}
Note that ${\sf fgt}_{\epsilon}(\sC^{\epsilon}_{\Lie}(\gg)[\![\hbar]\!]]) \simeq \Sym_{\Bbbk}(\gg^{\vee}[-1](1))[\![\hbar]\!]$. On generators $\xi^1, \xi^2 \in \gg^{\vee}[-1](1)$, we have 
\[\{\xi, \eta\}_{\lambda} = \hbar\cdot\lambda(\xi, \eta) = \hbar\cdot \lambda(\eta, \xi) = \{\eta, \xi\}_{\lambda}~.\]
Therefore the bracket is graded anti-symmetric. 
    
Given generators $\xi^1, \xi^2, \xi^3 \in \gg^{\vee}[-1](1)$, note that 
\[\{\xi^i, \{\xi^j, \xi^k\}_{\lambda}\}_{\lambda} = 0, \{i, j, k\} = \{1, 2, 3\}\]
for degree reasons. 
Therefore, the Jacobi identity is satisfied. 

Now we will check that $\{-, -\}_{\lambda}$ is compatible with the Chevalley--Eilenberg differential. 
Let $\{\xi^i\}$ be a basis for $\gg^{\vee}$. Write the Chevalley--Eilenberg differential as 
\[d_{\sf CE}(\xi^i) = \frac{-1}{2}c^{i}_{jk} \xi^j\xi^k\]
where $\{c^{i}_{jk}$\} are the structure constants for the dual basis $\{\xi_i\}$ of $\gg$. 
Also write the pairing $\lambda$ as 
\[\lambda = \lambda^{ij}\xi_i\xi_j\]
As the pairing $\lambda$ is $\GG$-invariant, for any $\xi_k$ we must have 
\[\begin{split}
        0 &= \lambda^{ij}{\sf ad}_{\xi_k}(\xi_i)\xi_j + \lambda^{ij}\xi_i{\sf ad}_{\xi_k}(\xi_j) \\
        &= \lambda^{ij}c^l_{ki}\xi_l\xi_j + \lambda^{ij}c^r_{kj}\xi_i\xi_r\
\end{split}\]
Extracting the coefficient of $\xi_a\xi_b$, we have 
\[0 = \lambda^{ib}c^a_{ki} + \lambda^{aj}c^{b}_{kj}~,\]
which we rewrite as 
\[0 = \lambda^{bj}c^a_{kj} + \lambda^{aj}c^b_{kj}\]
We will now check that 
\[0 = d_{\sf CE}(\{\xi^a, \xi^b\}_{\lambda}) = \{d_{\sf CE}(\xi^a), \xi^b\}_{\lambda} - \{\xi^a, d_{\sf CE}(\xi^b)\}_{\lambda}~.\]
We have 
\[\begin{split}
        \{d_{\sf CE}(\xi^a), \xi^b\}_{\lambda} &= \left\{\frac{-1}{2}c^a_{ij}\xi^i\xi^j, \xi^b\right\}_{\lambda} \\
        &= \frac{-\hbar}{2}\bigl(c^a_{ij}\lambda^{jb}\xi^i - c^a_{ij}\lambda^{ib}\xi^j\bigr) \\
        &= \frac{-\hbar}{2}\bigl(c^a_{ij}\lambda^{jb}\xi^i - c^a_{ji}\lambda^{jb}\xi^i\bigr) \\
        &= \frac{-\hbar}{2}\bigl(c^a_{ij}\lambda^{jb}\xi^i + c^a_{ij}\lambda^{jb}\xi^i\bigr) \\
        &= -\hbar(c^a_{ij}\lambda^{jb}\xi^i)
\end{split}\]
as we have $c^a_{ij} = - c^a_{ji}$. 
Similarly, we have 
\[
\begin{split}
        \{\xi^a, d_{\sf CE}(\xi^b)\}_{\lambda} &= \left\{\xi^a, \frac{-1}{2}c^b_{ij}\xi^i\xi^j\right\}_{\lambda} \\
        &= \frac{-\hbar}{2}\bigl(c^b_{ij}\lambda^{ai}\xi^j - c^b_{ij}\lambda^{aj}\xi^i\bigr) \\
        &= \frac{-\hbar}{2}\bigl(c^b_{ji}\lambda^{aj}\xi^i - c^b_{ij}\lambda^{aj}\xi^i\bigr) \\
        &= \frac{-\hbar}{2}\bigl(-c^b_{ij}\lambda^{aj}\xi^i - c^b_{ij}\lambda^{aj}\xi^i\bigr) \\
        &= \hbar(c^b_{ij}\lambda^{aj}\xi^i)~.
\end{split}
\]
Therefore, we have 
\[\{d_{\sf CE}(\xi^a), \xi^b\}_{\lambda} - \{\xi^a, d_{\sf CE}(\xi^b)\}_{\lambda} = -\hbar(c^a_{ij}\lambda^{jb}+c^b_{ij}\lambda^{aj})\xi^i = 0\]
by the $\GG$-invariance formula calculated above (note that $\lambda^{\mu \nu} =\lambda^{\nu \mu}$). 
Therefore we have the result. 
\end{proof}

We will now show that $(\sC^{\epsilon}_{\Lie}(\gg), \{-, -\}_{\lambda})$ is in fact a $\PP_3$-algebra in $\sC^{\epsilon}_{\Lie}(\gg, \cO(\GG))[\![\hbar]\!]$-comodules. 

\begin{lemma} \label{equivariance}
The $\PP_3$-bracket 
\[\{-, -\}_{\lambda}: \sC^{\epsilon}_{\Lie}(\gg)[\![\hbar]\!]\underset{\Bbbk[\![\hbar]\!]}\otimes\sC^{\epsilon}_{\Lie}(\gg)[\![\hbar]\!] \rightarrow \sC^{\epsilon}_{\Lie}(\gg)[\![\hbar]\!]\]
is a map of $\sC^{\epsilon}_{\Lie}(\gg, \cO(\GG))[\![\hbar]\!]$-comodules.
That is, 
\[
(\sC^{\epsilon}_{\Lie}(\gg)[\![\hbar]\!], \{-, -\}_{\lambda})
\in \Alg_{\PP_3}({\sf coMod}_{\sC^{\epsilon}_{\Lie}(\gg, \cO(\GG))[\![\hbar]\!]}(\Mod^{\epsilon, \gr}_{\Bbbk[\![\hbar]\!]})).
\]
\end{lemma}

\begin{proof}
The proof is by direct computation. Let $\rho'$ denote the coaction on $\sC^{\epsilon}_{\Lie}(\gg)[\![\hbar]\!]\underset{\Bbbk[\![\hbar]\!]}\otimes\sC^{\epsilon}_{\Lie}(\gg)[\![\hbar]\!]$. We will check that 
\[\rho \circ \{\cdot, \cdot\}_{\lambda} = ({\sf id} \otimes \{\cdot, \cdot\}_{\lambda})\circ \rho'.\]
Let $\{\xi^i\}$ be a basis for $\gg^{\vee}$. Define functions $\alpha^i_j \in \cO(\GG)$ by the formula 
\[{\sf coAd}_g(\xi^i) = \sum_j\alpha^i_j(g)\xi^j.\]
The coaction of $\sC^{\epsilon}_{\Lie}(\gg, \cO(\GG))[\![\hbar]\!]$ on $\sC^{\ast}_{\Lie}(\gg)[\![\hbar]\!]$
\[
\rho: \sC^{\epsilon}_{\Lie}(\gg)[\![\hbar]\!] \rightarrow \sC^{\epsilon}_{\Lie}(\gg, \cO(\GG))[\![\hbar]\!]\underset{\Bbbk[\![\hbar]\!]}\otimes\sC^{\epsilon}_{\Lie}(\gg)[\![\hbar]\!]
\]
satisfies
\[
\rho(\xi^i) = \xi^i\otimes 1 + \sum_j \alpha^i_j\otimes \xi^j.
\]
We have 
\[
(\rho \circ \{\cdot, \cdot\}_{\lambda})(\xi^i \otimes\xi^j) = \rho(\hbar \cdot \lambda(\xi^i, \xi^j)) = 1 \otimes \hbar \cdot \lambda(\xi^i, \xi^j).
\]
We will now compute $({\sf id} \otimes \{\cdot, \cdot\}_{\lambda})\circ \rho')(\xi^i\otimes \xi^j)$. 
Let $\tau$ denote the braiding on $\Mod^{\epsilon, \gr}_{\Bbbk[\![\hbar]\!]}$. Recall that
\[\rho' = (m\otimes \sf{id} \otimes \sf{id})\circ({\sf id} \otimes \tau \otimes {\sf id})\circ(\rho \otimes \rho)\]
where $m$ denotes the product on $\sC^{\epsilon}_{\Lie}(\gg, \cO(\GG))$. 
We have 
\begin{align*}
(\rho \otimes \rho)(\xi^i\otimes \xi^j) &= \xi^i\otimes 1 \otimes\xi^j \otimes 1 \\
&+ \sum \xi^i\otimes 1 \otimes \alpha^j_l \otimes \xi^l \\
&+ \sum \alpha^i_k\otimes \xi^k \otimes \xi^j \otimes1 \\
&+ \sum \alpha^i_k \otimes \xi^k \otimes \alpha^j_l \otimes \xi^l.
\end{align*}
Swapping the middle two factors introduces a sign (note that only the cohomological grading introduces signs).  
\begin{align*}
\tau \circ (\rho \otimes \rho)(\xi^i\otimes \xi^j) &= \xi^i\otimes\xi^j\otimes 1  \otimes 1 \\
&+ \sum \xi^i\otimes \alpha^j_l\otimes 1  \otimes \xi^l \\
&- \sum \alpha^i_k\otimes \xi^j\otimes \xi^k  \otimes1 \\
&+ \sum \alpha^i_k \otimes \alpha^j_l\otimes \xi^k  \otimes \xi^l,
\end{align*}
and therefore 
\begin{align*}
\rho'(\xi^i\otimes \xi^j) &= \xi^i \xi^j\otimes 1  \otimes 1 \\
&+ \sum \xi^i \alpha^j_l\otimes 1  \otimes \xi^l \\
&- \sum \alpha^i_k \xi^j\otimes \xi^k  \otimes1 \\
&+ \sum \alpha^i_k \alpha^j_l\otimes \xi^k  \otimes \xi^l.
\end{align*}
The composition $({\sf id} \otimes \{\cdot, \cdot\}_{\lambda})\circ \rho')(\xi^i\otimes \xi^j)$ is therefore given by 
\begin{align*}
({\sf id} \otimes \{\cdot, \cdot\}_{\lambda})\circ \rho')(\xi^i\otimes \xi^j) &= \sum \alpha^i_k \alpha^j_l\otimes \{\xi^k,\xi^l\}_{\lambda} \\
&= \sum \alpha^i_k \alpha^j_l\otimes \hbar \cdot \lambda(\xi^k,\xi^l) \\
&= 1 \otimes \hbar\cdot \lambda(\xi^i, \xi^j),
\end{align*}
where the last equality is given by $\GG$-invariance of $\lambda$. 
Therefore we have checked the desired equality on generators of $\gg^{\vee}$. 
Observe that the difference 
\[\rho \circ \{\cdot, \cdot\}_{\lambda} - ({\sf id} \otimes \{\cdot, \cdot\}_{\lambda})\circ \rho'\]
is a $\Bbbk[\![\hbar]\!]$-linear biderivation on $\sC^{\epsilon}_{\Lie}(\gg)[\![\hbar]\!]\underset{\Bbbk[\![\hbar]\!]}\otimes\sC^{\epsilon}_{\Lie}(\gg)[\![\hbar]\!]$ with coefficients in $\sC^{\epsilon}_{\Lie}(\gg)[\![\hbar]\!]$. 
This bi-derivation vanishes on generators so it is identically zero, and we have the result. 
\end{proof}

By formality of $\EE_3$, the following is immediate. 

\begin{cor} \label{e3def}
Given $\lambda \in \Sym^2(\gg)^{\GG}$, we have a $\GG_{\DR}$-equivariant $\EE_3$- deformation quantization of the graded mixed Chevalley--Eilenberg complex:  \[
\sC^{\epsilon}_{\lambda}(\gg)[\![\hbar]\!] \in \Alg_{\EE_3}({\sf coMod}_{\sC^{\epsilon}_{\Lie}(\gg, \cO(\GG))[\![\hbar]\!]}(\Mod^{\epsilon, \gr}_{\Bbbk[\![\hbar]\!]}))
\]
such that $\sC^{\epsilon}_{\lambda}(\gg)[\![\hbar]\!]\underset{\Bbbk[\![\hbar]\!]}\otimes\Bbbk_{\hbar = 0} \simeq \sC^{\epsilon}_{\Lie}(\gg)$. 
\end{cor} 

The following lemma is just a base-changed version of Lemma~\ref{ppascomod}.

\begin{lemma} \label{comp}
Consider the symmetric monoidal functor
\[ {\sf Assoc}^{\hbar}\colon {\sf coMod}_{\sC^{\epsilon}_{\Lie}(\gg, \cO(\GG))[\![\hbar]\!]}(\Mod^{\epsilon, \gr}_{\Bbbk[\![\hbar]\!]}) \xrightarrow[]{} \Mod^{\epsilon, \gr}_{\cD_{\sB\GG_{\DR}}[\![\hbar]\!]}\]
\[V \mapsto {\sf Assoc}^{\hbar}(V): (\Spec(R)\xrightarrow[]{f}\sB\GG_{\DR})\mapsto \Bbbk[\![\hbar]\!]\square (\DD(P)\otimes V)\]
where $P \ra \Spec(R)$ is the torsor classified by $\Spec(R)\xrightarrow[]{f}\sB\GG_{\DR}$.
This functor sends $\sC^{\epsilon}_{\Lie}(\gg)[\![\hbar]\!]$ to the presheaf of principal parts $\cP_{\sB\GG}[\![\hbar]\!]$. 
\end{lemma}

Therefore, we have the following result. 

\begin{prop} \label{quantum}
Given $\lambda \in \Sym^2(\gg)^{\GG}$, we obtain a $\cD_{\sB\GG_{\DR}}[\![\hbar]\!]$-linear $\EE_3$-deformation of the prestack of principal parts of $\sB\GG$: an object $\cP^{\lambda}_{\sB\GG}[\![\hbar]\!] \in \Alg_{\EE_3}(\Mod^{\epsilon, \gr}_{\cD_{\sB\GG_{\DR}}[\![\hbar]\!]})$ together with an equivalence \[
\cP^{\lambda}_{\sB\GG}[\![\hbar]\!] \underset{\cD_{\sB\GG_{\DR}}[\![\hbar]\!]}\otimes(\cD_{\sB\GG_{\DR}})_{\hbar =0} \simeq \cP_{\sB\GG}~.
\]
\end{prop}

\begin{proof}
Combine Corollary~\ref{e3def} and Lemma~\ref{comp}. 
\end{proof}
    
\subsection{Quantizing $\Loc_{\GG}(M)$}
Let $M$ be a smooth $d$-manifold with the homotopy type of a finite CW complex. For $d \leq 3$, given a framing of $M \times \RR^{3-d}$, we will construct quantizations of $\Loc_{\GG}(M) \simeq \Map(\Pi_\infty M, \sB\GG)$. 

Given $\lambda\in \Sym^2(\gg)^{\GG}$,  
we obtain a graded mixed $\EE_3$-algebra 
\[
\cP^{\lambda}_{\sB\GG}[\![\hbar]\!] \in \Alg_{\EE_3}(\Mod^{\epsilon, \gr}_{\cD_{\sB\GG_{\DR}}[\![\hbar]\!]})
\]
deforming principal parts of $\sB\GG$, by Proposition~\ref{quantum}. 
Definition~\ref{diskalgfrompoiss} constructs a $\Disk_M$-algebra 
\[
\cA_{\sB\GG}^{\lambda}[\![\hbar]\!] \in \Alg_{\Disk_M}(\Mod^{\epsilon, \gr}_{\cD_{\Loc_{\GG}(M)_{\DR}}[\![\hbar]\!]})
\]
deforming $\cA_{\sB\GG}$. 
Proposition~\ref{defpp} then gives us a quantization of $\Loc_{\GG}(M)$. 
Therefore we have the following proposition. 

\begin{prop} \label{quantchernsimons}
Given $\lambda\in \Sym^2(\gg)^{\GG}$, and a framing of $M \times \RR^{3-d}$, we obtain  an $\EE_{3-d}$-deformation of the presheaf of principal parts of $\Loc_{\GG}(M)$:
\[\cP_{\Loc_{\GG}(M)}^{\lambda}[\![\hbar]\!]:= \int_{M} \cA^{\lambda}_{\sB\GG}[\![\hbar]\!]~.\]
\end{prop}

\begin{remark}
We describe $\cP_{\Loc_{\GG}(M)}^{\lambda}[\![\hbar]\!]$ at a $\Bbbk$-point $\Spec(\Bbbk)\xrightarrow[]{f}\Loc_{\GG}(M)_{\DR}$. 
Note that $f$ picks out a $\GG$-local system $L_{f}$ on~$M$. 
Reducing mod $\hbar$, we have 
\begin{align*}
\cP_{\Loc_{\GG}(M)}^{\lambda}[\![\hbar]\!](f)\underset{\Bbbk[\![\hbar]\!]}\otimes\Bbbk_{\hbar = 0} 
&\simeq \cP_{\Loc_{\GG}(M)}(f)\\ 
&\simeq \DD(\Spec(\Bbbk)\times_{\Loc_{\GG}(M)_{\DR}}\Loc_{\GG}(M))\\ 
&\simeq \sC^{\epsilon}_{\Lie}(\Omega_c^{\bullet}(M, \gg_{L_f}))
\end{align*}
where $\gg_{L_f}$ is the adjoint bundle of the local system $L_f$. At the trivial local system, we obtain a deformation of $\sC^{\epsilon}_{\Lie}(\Omega_c^{\bullet}(M, \gg))$. This recovers a quantization described in~\cite{CFG}.  
~\hfill$\Diamond$\end{remark}

If we have $d \leq 2$, then $\Perf_{\cP^{\lambda}_{\Loc_{\GG}(M)}}$ has the structure of an $\EE_{2-d}$-category. The formal localization result of \cite{CPTVV} also gives the following. 
 
\begin{prop} \label{perfchernsimons}
For $d \leq 2$, given $\lambda \in \Sym^2(\gg)^{\GG}$, and a framing of $M \times \RR^{3-d}$, we obtain an $\EE_{2-d}$-category: 
\[
\Perf^{\lambda}(\Loc_{\GG}(M)) := \Perf^{\epsilon, \gr}_{\cP_{\Loc_{\GG}(M)}^{\lambda}}.
\]
\end{prop}

\begin{remark}
If $\dim(M)\leq 2$, \cite{BBJ} construct quantizations of $\qc(\Loc_{\GG}(M))$ using different methods. There is an equivalence of symmetric monoidal categories \cite{BZFN}
\[\qc(\Loc_{\GG}(M)) \simeq \int_M \Rep(\GG)~,\]
where the right-hand side is factorization homology of categories. Therefore given an $\EE_2$-deformation $\Rep_q(\GG)$ of $\Rep(\GG)$, one obtains a deformation 
\[\qc_q(\Loc_{\GG}(M)) \simeq \int_M \Rep_q(\GG).\]
 It is interesting to note that our quantizations of $\Loc_{\GG}(M)$ are obtained by computing factorization homology of an $\EE_3$-algebra instead of an $\EE_2$-category. We intend to return to the relation of these two quantizations in future work.
~\hfill$\Diamond$
\end{remark}

\subsection{Quantization yields a nontrivial deformation}
In this subsection, we will show that the quantizations of $\Loc_{\GG}(M)$ produced above are in fact nontrivial. We begin with a computation. 

Let $V$ be a finite dimensional $\Bbbk$-vector space concentrated in degree $0$, viewed as an Abelian Lie algebra. Pick a pairing $\lambda\in \Sym^2(V)$, and consider the graded mixed $\Bbbk[\![\hbar]\!]$-linear $\PP_3$-algebra 
\[\Bigl(\Sym_{\Bbbk}(V^{\vee}[-1](1))\otimes \Bbbk[\![\hbar]\!], \{-,-\}_{\lambda}\Bigr)~.\]
Note that the internal and mixed differentials are zero. Choose a basis $(\xi^i)_i$ of $V^{\vee}[-1]$ that diagonalizes the pairing
\[\lambda(\xi^i, \xi^j) = c_i\delta_{ij}~.\]

We have 
\[\Sym_{\Bbbk}(V^{\vee}[-1](1))\otimes \Bbbk[\![\hbar]\!] \simeq \Bbbk[\![\xi^1,\dots \xi^k, \hbar]\!]\]
where $\xi^i$ have weight 1 and cohomological degree 1, and $\hbar$ has weight 2 and cohomological degree 0. 
The $\PP_3$-bracket on generators is given by 
\[\{\xi^i, \xi^j\}_{\lambda} = \delta_{ij}c_i\hbar, \]
\[\{\xi^i, \hbar\}_{\lambda} = 0 = \{\hbar, \hbar\}_{\lambda}~.\]
The constructions of \cite{Markarian2017Weyl} Proposition 9, and \cite{CPTVV} Section 3.6 provide an explicit $\EE_3$-quantization of $\PP_3$-algebras of the form above. Let $FM_3$ denote the Fulton-MacPherson model for the $\EE_3$-operad. Pick an $\sS\sO(3)$-invariant volume form $\omega$ on $FM_3(2) \cong S^2$. Let $\omega_{ij}$ denote the pullback of $\omega$ along the map $FM_3(I) \ra FM_3(2)$ that forgets all points except $i$ and $j$. 
The $\EE_3$-quantization of $A^{\lambda} = \Bigl(\Sym_{\Bbbk}(V^{\vee}[-1](1))\otimes \Bbbk[\![\hbar]\!], \{-,-\}_{\lambda}\Bigr)$ can be described by the structure maps 
\[(A^{\lambda})^{I} \xra{{\sf exp}(a)} (A^{\lambda})^I\otimes \sC^{\ast}(FM_3(I), \Bbbk) \xra{m\otimes {\sf id}}A^{\lambda}\otimes \sC^{\ast}(FM_3(I), \Bbbk)\]
where 
\[a = \sum_{i\neq j}\partial_{\lambda}^{i,j}\otimes \omega_{ij}~,\]
where $\partial_{\lambda}^{i,j}$ is $\{-, -\}_{\lambda}$ applied to the $i$-th and $j$-th factors. The exression ${\sf exp}(a)$ is a formal sum, which only makes sense upon evaluating on chains. In particular, the binary operation is given by 
\[\Bbbk[\![\xi^1,\dots \xi^k, \hbar]\!] \otimes_{\Bbbk[\![\hbar]\!]} \Bbbk[\![\xi^1,\dots \xi^k, \hbar]\!] \ra \Bbbk[\![\xi^1,\dots \xi^k, \hbar]\!]\otimes\sC^{\ast}(S^2, \Bbbk)\]
\[\xi^i\otimes\xi^j \mapsto \xi^i\xi^j\otimes1 + \delta_{ij}c_{i}\hbar \otimes \omega~,\]
i.e. the fundamental class of $S^2$ acts via the Poisson bracket. 

Consider the graded mixed $\Bbbk[\![\hbar]\!]$-linear $\EE_0$-algebra 
\[P^{\lambda} := \int_{S^3}\Bbbk[\![\xi^1,\dots \xi^k, \hbar]\!]_{\lambda}\]
where we view $\Bbbk[\![\xi^1,\dots \xi^k, \hbar]\!]_{\lambda}$ as the $\Disk_{S^3}$-algebra corresponding to the trivial local system of $\EE_3$-algebras on $S^3$ with fiber $\Bbbk[\![\xi^1,\dots \xi^k, \hbar]\!]_{\lambda}$. 

Also consider the object 
\[P^{\sf triv}:= \int_{S^3} \Bbbk[\![\xi^1,\dots \xi^k, \hbar]\!]~,\]
where $\Bbbk[\![\xi^1,\dots \xi^k, \hbar]\!]$ is given the commutative algebra structure. 
\begin{lemma}
    We have an equivalence of graded mixed $\Bbbk[\![\hbar]\!]$-modules 
    \[P^{\sf triv} \simeq \Bbbk[\![p_1, \dots, p_k, q_1, \dots, q_k, \hbar]\!]\]
    where 
    $p_i$ has weight $1$ and cohomological degree $1$, $q_i$ has weight $1$ and cohomological degree $-2$, and $\hbar$ has weight 2 and cohomological degree $0$. 
\end{lemma}
\begin{proof}
    By Lemma~\ref{colim}, we have 
    \[P^{\sf triv} := \int_{S^3} \Bbbk[\![\xi^1,\dots \xi^k, \hbar]\!] \simeq \Sym(H_{\ast}(S^3, \Bbbk)[-1](1))\otimes_{\Bbbk}\Bbbk[\![\hbar]\!] \simeq \Bbbk[\![p_1, \dots, p_k, q_1, \dots, q_k, \hbar]\!]\]
    where $p_i:= [{\sf pt}]\otimes x_i$, and $q_i:= [S^3]\otimes x_i$. Note that both the internal and mixed differentials are 0.
\end{proof}

We will now show that $P^{\lambda}$ is not equivalent to $P^{\sf triv}$ viewed as deformations of 
\[P = \int_{S^3}\Bbbk[\![\xi^1,\dots \xi^k]\!] \simeq  \Bbbk[\![p_1, \dots, p_k, q_1, \dots, q_k]\!]\in \CAlg(\Mod_{\Bbbk}^{\epsilon, \gr})\]
This will follow from a spectral sequence computation. The following argument is very similar to \cite{Markarian2017Weyl} Theorem 1. It is also similar to the argument in \cite{Costello2016FactorizationAI} Section 5.1. 
\begin{lemma} \label{nontrivialdeformation}
    For $\lambda\neq 0$, the object $P^{\lambda}$ is not equivalent to $P^{\triv}$ as objects of $\Alg_{\EE_0}\Bigl(\Mod^{\epsilon, \gr}_{\Bbbk[\![\hbar]\!]}\Bigr)$~.
\end{lemma}
\begin{proof}
    We will verify that the semiclassical structure associated to the $\EE_0$-deformation above is nontrivial. Consider the filtration on $P^{\lambda}$ by powers of $\hbar$:
    \[F^p \colon = P^{\lambda} \otimes_{\Bbbk[\![\hbar]\!]}\hbar^p\Bbbk[\![\hbar]\!].\]
    Its associated graded is given by 
    \[\gr^{\bullet} = P[\![\hbar]\!] \simeq \Bbbk[\![p_1, \dots, p_k, q_1, \dots, q_k, \hbar]\!]\]
    following the notation of the previous lemma: $p_i$ has weight $1$ and cohomological degree $1$, $q_i$ has weight $1$ and cohomological degree $-2$, and $\hbar$ has weight 2 and cohomological degree $0$. Note that there are three gradings involved: the cohomological grading, weight grading, and powers of $\hbar$ gradings. The first page of the spectral sequence associated to the filtration is given by 
    \[E_1^{p,q} = H^{p+q}(\gr^p)  = H^{p+q}(\hbar^p\Bbbk[\![p_1, \dots, p_k, q_1, \dots, q_k]\!]) = (\hbar^p\Bbbk[\![p_1, \dots, p_k, q_1, \dots, q_k]\!])^{p+q}~.\]
    Note that the superscript refers to the {\em cohomological} grading. The differential on the first page is given by 
    \[d_1= \hbar\sum_i c_i\partial_{p_i}\partial_{q_i}~.\]
    See the proof of \cite{Markarian2017Weyl} Theorem 1 for a calculation of the differential: Consider the weight $2$ piece of the $E_1$-page. It is generated as a $\Bbbk$-vector space by degree 2 monomials in $p_i$ and $q_j$ (note that they have weight 1), and $\hbar$ (which has weight 2). Therefore the only potentially nonzero differential in weight 2 is given by 
    \[d_1: \bigoplus_{i, j}\Bbbk\cdot p_iq_j=\bigl(E^{0,-1}_1\bigr)_2 \ra \bigl(E_1^{1,-1}\bigr)_2 = \Bbbk[\![\hbar]\!]\]
    \[p_iq_j \mapsto \delta_{ij}c_i\hbar~.\]
    As we assumed $\lambda \neq 0$, we must have $c_i \neq 0$ for some $i$, and therefore the class $\hbar$ in the weight $2$ piece of the $E_1$-page is killed. Observe that no other classes can contribute to $H^0(P^\lambda_2)$, therefore we must have 
    \[H^0(P^\lambda_2) = 0~.\]
    However, for the trivial deformation $P^{\sf triv}$, we have 
    \[H^0(P^{\sf triv}_2) = P_2^{\sf triv} = \Bbbk\cdot\hbar \neq 0.\]
    Therefore we have the result. 
\end{proof}

Now let $\GG$ be an affine algebraic group over $\Bbbk$, and consider the derived stack of $\GG$-local systems $\Loc_{\GG}(S^3)$. Given a level $\lambda\in \Sym^2(\gg)^{\GG}$, consider the $\EE_0$-quantization 
\[\cP^{\lambda}_{\Loc_{\GG}(S^3)}:= \int_{S^3}\cA_{\sB\GG}^{\lambda}[\![\hbar]\!]\]
constructed in Proposition~\ref{quantchernsimons}. 

\begin{prop}
    The quantization
    \[\cP^{\lambda}_{\Loc_{\GG}(S^3)}:= \int_{S^3}\cA_{\sB\GG}^{\lambda}[\![\hbar]\!]\]
    is not the trivial deformation of principal parts on $\Loc_{\GG}(S^3)$.
\end{prop}
\begin{proof}
    At the trivial local system 
    \[f\colon \Spec(\Bbbk) \ra \Loc_{\GG}(S^3)\]
    we have an $\EE_0$-algebra in $\Mod^{\epsilon, \gr}_{\Bbbk[\![\hbar]\!]}$
    \[\cP^{\lambda}_{\Loc_{\GG}(S^3)}(f) \simeq \int_{S^3}\cA_{\sB\GG}^{\lambda}[\![\hbar]\!](f)\simeq \int_{S^3}\sC^{\lambda}_{\Lie}(\gg)[\![\hbar]\!]~,\]
    where the expression on the right denotes the $\Disk_{S^3}$-algebra corresponding to the trivial local system on $S^3$ with fiber $\sC^{\lambda}_{\Lie}(\gg)[\![\hbar]\!] \in \Alg_{\EE_3}(\Mod_{\Bbbk[\![\hbar]\!]}^{\epsilon, \gr})$. Upon applying ${\sf fgt_{\epsilon}}$, we are reduced to the case studied in Lemma~\ref{nontrivialdeformation}, and therefore we have the result. 
\end{proof}

\bibliographystyle{alphaurl}
\bibliography{ref}
\end{document}